\documentstyle[12pt]{article}
\begin{document}
\title{Differentiable functions of Cayley-Dickson numbers.}
\author{S.V. Ludkovsky.}
\date{27.05.2004}
\maketitle

\begin{abstract}
We investigate superdifferentiability of functions defined on regions of the
real octonion (Cayley) algebra and obtain a noncommutative version
of the Cauchy-Riemann conditions. Then we study
the noncommutative analog of the Cauchy integral
as well as criteria for functions of an octonion variable
to be analytic. In particular, the octonion
exponential and logarithmic functions are being considered.
Moreover, superdifferentiable functions of variables belonging to
Cayley-Dickson algebras (containing the octonion algebra as the proper
subalgebra) finite and infinite dimensional are investigated.
Among main results there are the Cayley-Dickson algebras
analogs of Caychy's theorem, Hurtwitz', argument principle,
Mittag-Leffler's, Rouche's and Weierstrass' theorems.
\end{abstract}

\section{Introduction}
Functions of real variables with values in Clifford algebras
were investigated, for example, in \cite{brdeso}.
In this article we continue our investigations of functions
of variables belonging to noncommutative superalgebras
\cite{luoyst} considering here functions of octonion variables
and also of more general Cayley-Dickson algebras ${\cal A}_r$,
$r\ge 4$, containing the octonion algebra as the proper subalgebra.
The algebra of octonions is alternative, that gives possibility
to define residues of functions in a reasonable way.
The Cayley-Dickson algebras of larger dimensions are not already
alternative and proceedings for them are heavier.
Nevertheless, using their power-associativity and distributivity
it is possible to define differentiability of functions of variables
belonging to the Cayley-Dickson algebra such that the differential
has sufficiently well properties, that to define subsequently
the line integral over such variable. This integral is extended
on spaces of continuous functions along rectifiable paths.
It is necessary to note, that the graded structure of the Cayley-Dickson
algebra ${\cal A}_r$ over $\bf R$ and its noncommutativity for the real
dimension not less, than $4$ causes well properites of
superdifferentiable functions $f: U\to {\cal A}_r$, where $U$ is
open in ${\cal A}_r^n$. Apart from the complex case a derivative
$f'$ of $f\in C^1(U,{\cal A}_r)$ is an operator and not  in general
a number even for $U\subset {\cal A}_r$, $n=1$.
\par The theory of ${\cal A}_r$-holomorphic functions
given below can not be reduced to the theory of holomorphic functions
of several complex variables. Moreover, ${\cal A}_r$-holomorphic
functions have many specific features in comparison with real
locally analytic functions as, for example, the argument principle,
homotopy theorem $2.15$, theorem about representations of multiples of
functions with the help of line integral along a loop, etc. show.
\par Dirac had used biquaternions (complexified quaternions)
in his investigations of quantum mechanics. The Dirac
equations $D_zf_1 - D_{\tilde z}f_2=m(f_1+f_2)$ and
$D_{\tilde z}f_1 + D_zf_2=m(f_1-f_2)$ on the space
of right superlinearly $(z,{\tilde z})$-superdifferentiable
functions of quaternion variable
can be extended on the space of $(z,{\tilde z})$-superdifferentiable
functions $f_1(z,{\tilde z})$ and $f_2(z,{\tilde z})$
(see our definition in \S 2.2 below and in \cite{luoyst}) gives
evident physical interpretation of a solution $(f_1,f_2)$, $r\ge 2$, as
spinors, where $m$ is a mass of an elementary particle.
We extend the operator $(D_z^2+D_{\tilde z}^2)$ from
the space of right superlinearly $(z,{\tilde z})$-superdifferentiable
functions on the space of $(z,{\tilde z})$-superdifferentiable functions $f$, 
hence we get the Klein-Gordon equation $(D_z^2+D_{\tilde z}^2)f=2m^2f$
in the particular case of quaternions and ${\cal A}_r$-algebras, $r\ge 3$.
It is necessary to note that previous authors
have used right (or left) superlinearly superdifferentiable functions,
that does not form an algebra of functions and they have used
multiple and iterated integrals and the Gauss-Ostrogradskii-Green
formula, but they have not used line integrals over ${\cal A}_r$
(see, for example, \cite{guetze} and references therein).
While development their theory Yang and Mills known in theoretical
and mathematical physics had actively worked with quaternions, but
they have felt lack of the available theory of quaternion functions
existing in their time. Yang also have expressed the idea, that
possibly in quantum field theory it is worthwhile to use quaternion time
(see page 198 \cite{guetze}). It is known also the
use of complex time through the Wick rotation in quantum mechanics
for getting solutions
of problems, where the imaginary time is used for interpretations
of probabilities of tunneling under energy barriers (walls).
Using the special unitary group
embedded into the quaternion skew field $\bf H$ it makes equivalent
under isomoprphism with $SO(3)$ all spatial axes. On the other hand,
the major instrument for measurement is the spectrum.
When there are deep energy wells or high energy walls, then
it makes obstacles for penetrating electromagnetic waves and radiation,
that is well known also in astronomy, where black holes are actively
studied (see page 199 and \S 3.b \cite{guetze} and references therein).
W. Hamilton in his lectures on quaternions also tackled a question of
events in $\bf H$ and had thought about use of quaternions in astronomy and
celestial mechanics (see \cite{hamilt,rothe} and references therein).
Therefore, in general to compare the sequence of events it may be
necessary in definite situations to have the same dimensional time space
as the coordinate space. On the other hand, spatial isotropy at least local
in definite domains makes from each axis under rotations and dilatations
$SU(2)\times {\bf R}$ isomorphic with $\bf H$. Therefeore, it appears
that in definite situations it would be sufficent to use ${\bf H}^4$
instead of the Minkowski space-time ${\bf R}^{1,3}$,
where ${\bf R}^{1,3}$ has the embedding into ${\bf H}$.
Since $\bf H$ as the $\bf R$-linear space is isomorphic with
$\bf R^4$, then there exists the embedding $\zeta $ of $\bf R^{1,3}$
into $\bf H$ such that $\zeta (x_1,x_i,x_j,x_k)=
x_1+x_ii+x_jj+x_kk$, where the $\bf R^{1,3}$-norm is given by
the equation $|x|_{1,3}=(x^2+{\tilde x}^2)/2=Re (x^2)=
x_1^2-x_i^2-x_j^2-x_k^2$ and the $\bf R^{1,3}$ scalar product
is given by the equality
$(x,y)_{1,3}:=(xy+{\tilde y}{\tilde x})/2=Re (xy)$
$=x_1y_1-x_iy_i-x_jy_j-x_ky_k$, where $x=x_1+x_ii+x_jj+x_kk$,
$x_1,...,x_k\in \bf R$ (see \S 2.1). This also can be used
for embeddings of hyperbolic manifolds into quaternion manifolds.
Then ${\bf H}^4$ can be embedded into the algebra ${\cal A}_4$
of sedenions. It is also natural for describing systems with spin,
isospin, flavor, color and their interactions.
The enlargement of the space-time
also is dictated in some situations by symmetry properties of
differential equations or a set of operators describing a system.
For example, special unitary groups $SU(n)$, for $n=3, 5-8, 11$, etc.,
exceptional Lie groups,
are actively used in theory of elementary particles \cite{guetze}, but these
groups can be embedded into the corresponding
Cayley-Dickson algebra ${\cal A}_r$ \cite{baez}.
Indeed, $U(m)\subset GL(n,{\bf C})\subset {\bf C}^{n^2}$ while
${\bf C}^m$ has the embedding into ${\cal A}_r$, where
$m=2^{r-1}$, such that ${\bf C}^m \ni (x^1+iy^1,...,
x^m+iy^m)=:\xi \mapsto z:= (x^1+i_1y^1)+i_2(x^2+i_2^*i_3y^2)+...+
i_{2m-2}(x^m+i_{2m-2}^*i_{2m-1}y^m)\in {\cal A}_r$, since
$(i_l^*i_k)^2=-1$ for each $l\ne k\ge 1$, where $ \{ i_0,i_1,...,
i_{2m-1} \} $ is the basis of generators of ${\cal A}_r$, $i_0=1$,
$i_k^2=-1$, $i_0i_k=i_ki_0$, $i_li_k=-i_ki_l$ for each $k\ne l\ge 1$,
$i=(-1)^{1/2}$, $z^*=x^1-i_1y^1-i_2x^2-i_3y^2-...-i_{2m-2}x^m-
i_{2m-1}y^m$, the norm $(zz^*)^{1/2}=:|z|=(\sum_{k=1}^m
|x^k+iy^k|^2)^{1/2}=:|\xi |$ satisfies the parallelogramm
identity and induces the scalar product.
\par This paper as the previous one \cite{luoyst} is devoted to
the solution of the W. Hamilton problem of developing line integral and
holomorphic function theory of quaternion variables, but now
we consider general case of Cayley-Dickson algebras variables.
\par In this paper we investigate differentiability of functions
defined on regions of the real octonion (Cayley) algebra.
For this we consider specific definition of superdifferentiability
and obtain a noncommutative version of the Cauchy-Riemann conditions
in the particular case of right superlinear superdifferentiability.
Then we study the noncommutative analog of the Cauchy integral
as well as criteria for functions of an octonion variable
to be analytic. In particular, the octonion
exponential and logarithmic functions are being considered.
Moreover, superdifferentiable functions of variables belonging to
Cayley-Dickson algebras (containing the octonion algebra as the proper
subalgebra) finite and infinite dimensional are investigated.
Among main results there are the Cayley-Dickson algebras analogs
of Caychy's theorem,
Hurtwitz', argument principle, Mittag-Leffler's, Rouche's
and Weierstrass' theorems.
\par The results of this paper can serve for subsequent investigations
of special functions of Cayley-Dickson algebra variables,
noncommutative sheaf theory, manifolds of noncommutative geometry
over Cayley-Dickson algebras, their groups of loops and diffeomorphisms
(see also \cite{lulgcm,lulsqm,lustptg,lupm,oystaey}).
\section{Differentiability of functions of octonion variables}
\par To avoid misunderstandings we first introduce notations.
\par {\bf 2.1.} We write $\bf H$ for the skewfield of quaternions over
the real field $\bf R$ with the classical quaternion
basis $1$, $i$, $j$, $k$, satisfying relations of ${\cal A}_2$
(see Introduction).
The quaternion skewfield $\bf H$ has an anti-automorphism $\eta $
of order two $\eta: z\mapsto {\tilde z}$, where ${\tilde z}=
w_1-w_ii-w_jj-w_kk$, $z=w_1+w_ii+w_jj+w_kk$; $w_1,...,w_k\in \bf R$.
There is a norm in $\bf H$
such that $|z|=|z{\tilde z}|^{1/2}$, hence ${\tilde z}=|z|^2z^{-1}$.
\par The algebra $\bf K$ of octonions (octaves, the Cayley algebra)
is defined as an eight-dimensional algebra over $\bf R$ with a basis,
for example,
\par $(1)$ ${\bf b}_3:={\bf b}:= \{ 1, i, j, k, l, il, jl, kl \} $ such that
\par $(2)$ $i^2=j^2=k^2=l^2=-1$, $ij=k$, $ji=-k$,
$jk=i$, $kj=-i$, $ki=j$, $ik=-j$, $li=-il$, $jl=-lj$, $kl=-lk$;
\par $(3)$ $(\alpha +\beta l)(\gamma +\delta l)=(\alpha \gamma
-{\tilde {\delta }}\beta )+(\delta \alpha +\beta {\tilde {\gamma }})l$ \\
is the multiplication law in $\bf  K$
for each $\alpha $, $\beta $, $\gamma $, $\delta \in \bf H$,
$\xi :=\alpha +\beta l\in \bf K$, $\eta :=\gamma +\delta l\in \bf K$,
${\tilde z}:=v-wi-xj-yk$ for a quaternion $z=v+wi+xj+yk
\in \bf H$ with $v, w, x, y\in \bf R$.
\par The octonion algebra is neither commutative, nor associative,
since $(ij)l=kl$, $i(jl)=-kl$, but it is distributive and
${\bf R}1$ is its center. If $\xi :=\alpha +\beta l\in \bf K$,
then
\par $(4)$ ${\tilde {\xi }}:={\tilde {\alpha }}-\beta l$ is called
the adjoint element of $\xi $, where $\alpha , \beta \in \bf H$.  Then
\par $(5)$ $(\xi \eta )^{\tilde .}={\tilde {\eta }}{\tilde {\xi }}$,
${\tilde {\xi }} + {\tilde {\eta }}= (\xi +\eta )^{\tilde .}$ and
$\xi {\tilde {\xi }}=|\alpha |^2+|\beta |^2$, \\
where $|\alpha |^2=\alpha {\tilde {\alpha }}$ such that
\par $(6)$ $\xi {\tilde {\xi }}=:|\xi |^2$ and $|\xi |$
is the norm in $\bf K$. Therefore,
\par $(7)$ $|\xi \eta |=|\xi | |\eta |$, \\
consequently, $\bf K$ does not contain divisors of zero
(see also \cite{kansol,kurosh,ward}). The multiplication
of octonions satsifies equations:
\par $(8)$ $(\xi \eta )\eta =\xi (\eta \eta )$,
\par $(9)$ $\xi (\xi \eta )=(\xi \xi )\eta $, \\
that forms the alternative system.
In particular, $(\xi \xi )\xi =\xi (\xi \xi )$. Put
${\tilde {\xi }}=2a-\xi $, where $a=Re (\xi ):=(\xi + {\tilde {\xi }})/2
\in \bf R$. Since ${\bf R}1$ is the center of $\bf K$ and
${\tilde {\xi }}\xi =\xi {\tilde {\xi }}=|\xi |^2$, then
from $(8,9)$ by induction it follows, that for each
$\xi  \in \bf K$ and each $n$-tuplet (product), $n\in \bf N$,
$\xi (\xi (... \xi \xi )... )=(...(\xi \xi )\xi ...)\xi $
the result does not depend on an order of brackets (order
of consequtive multiplications), hence the definition of $\xi ^n:=
\xi (\xi (... \xi \xi )...)$ does not depend on the order of brackets.
This also shows that $\xi ^m\xi ^n=\xi ^n\xi ^m$,
$\xi ^m{\tilde {\xi }}^m={\tilde {\xi }}^m\xi ^n$ for each
$n, m\in \bf N$ and $\xi \in \bf K$.
Apart from the quaternions, the octonion algebra can not be
realized as the subalgebra of the algebra ${\bf M}_8({\bf R})$
of all $8\times 8$-matrices over $\bf R$, since $\bf K$ is not associative,
but ${\bf M}_8({\bf R})$ is associative.
The noncommutative nonassociative octonion algebra $\bf K$
is the $\bf Z_2$-graded $\bf R$-algebra ${\bf K}={\bf K}_0+{\bf K}_1$, where
elements of ${\bf K}_0$ are $\underline {even}$ 
and elements of ${\bf K}_1$ are $\underline {odd}$
(see, for example, \cite{kosshaf,kurosh,waerd}).
There are the natural embeddings ${\bf C}\hookrightarrow \bf K$
and ${\bf H}\hookrightarrow \bf K$, but neither $\bf K$ over $\bf C$,
nor $\bf K$ over $\bf H$, nor $\bf H$ over $\bf C$ are algebras, since
the centres of them are $Z({\bf H})=Z({\bf K})=\bf R$.
\par We consider also the Cayley-Dickson algebras ${\cal A}_n$
over $\bf R$, where $2^n$ is its dimension over $\bf R$.
They are constructed by induction starting from $\bf R$
such that ${\cal A}_{n+1}$ is obtained from ${\cal A}_n$
with the help of the doubling procedure, in particular,
${\cal A}_0:=\bf R$, ${\cal A}_1=\bf C$, ${\cal A}_2=\bf H$,
${\cal A}_3=\bf K$ and ${\cal A}_4$ is known as the sedenion algebra
\cite{baez}.
The Cayley-Dickson algebras are $*$-algebras, that is, there
is a real-linear mapping ${\cal A}_n\ni a\mapsto a^*\in {\cal A}_n$
such that
\par $(10)$ $a^{**}=a$,
\par $(11)$ $(ab)^*=b^*a^*$ for each $a, b\in {\cal A}_n$.
Then they are nicely normed, that is,
\par $(12)$ $a+a^*=:2 Re (a) \in \bf R$ and
\par $(13)$ $aa^*=a^*a>0$ for each $0\ne a\in {\cal A}_n$.
The norm in it is defined by
\par $(14)$ $|a|^2:=aa^*$.
We also denote $a^*$ by $\tilde a$.
Each $0\ne a\in {\cal A}_n$ has a multiplicative inverse
given by $a^{-1}=a^*/|a|^2$.
\par The doubling procedure is as follows. Each $z\in {\cal A}_{n+1}$
is written in the form $z=a+bl$, where $l^2=-1$, $l\notin {\cal A}_n$,
$a, b \in {\cal A}_n$. The addition is componentwise. The conjugate is
\par $(15)$ $z^*:=a^* -bl$. \\
The multiplication is given by Equation $(3)$.
\par The basis of ${\cal A}_n$ over $\bf R$ is denoted by
${\bf b}_n:= {\bf b}:= \{ 1, i_1,...,i_{2^n-1} \} $, where
$i_s^2=-1$ for each $1\le s \le 2^n-1$, $i_{2^{n-1}}:=l$
is the additional element of the doubling procedure of
${\cal A}_n$ from ${\cal A}_{n-1}$, choose $i_{2^{n-1}+m}=
i_ml$ for each $m=1,...,2^{n-1}-1$, $i_0:=1$.
\par An algebra is called alternative, if each its
subalgebra generated by two elements is associative.
An algebra is called power-associative, if its any subalgebra
generated by one element is associative.
Only for $n=0,...,3$ the Cayley-Dickson
algebras ${\cal A}_n$ are division alternative algebras.
For $n\ge 4$ the Cayley-Dickson algebras ${\cal A}_n$ are not
division algebras, but they are power-associative.
To verify the latter property consider $z\in {\cal A}_n$ written
in the form $z=v+M$, where $v=Re (z)$, $M:=(z-z^*)/2=:Im (z)$.
Then $v$ and $M$ commute and they are orthogonal,
$M^*=-M$. Therefore, the subalgebra generated by $z$ is associative
if and only if the algebra generated by $M$ is associative.
Since $M^*M=MM^*=|M|^2$ and $M^*=-M$, then the subalgebra generated
by $M$ is associative.
\par Verify that $z$ and $\tilde z$ in ${\cal A}_r$
are independent variables.
Suppose contrary that there exists $\gamma \in {\cal A}_r$ such that
$z+\gamma \tilde z=0$ for each $z\in {\cal A}_r$.
Write $z=a+bl$, $\gamma =\alpha +\beta l$,
where $a, b, \alpha , \beta \in {\cal A}_{r-1}$.
Then $z+\gamma \tilde z=0$ is equivalent to $\alpha a^*+b^*\beta =-a$
and $-b\alpha +\beta a=-b$. Consider $z$ with $a\ne 0$ and $b\ne 0$.
Then from the latter two equations we get:
$\alpha =-(a+b^*\beta )a|a|^{-2}=
1+b^*\beta a|b|^{-2}$, since $a^{-1}=a^*|a|^{-2}$.
This gives $\beta =-2b Re (a)|z|^{-2}$ and
$\alpha =1-2aRe (a)|z|^{-2}$. Taking in particular
$|z|=1$ and $Re (a)\ne 0$ and varying $z$ we come to the contardiction,
since $\gamma $ is not a constant. Therefore, $z$ and $z^*$
are two variables in ${\cal A}_r$ left-linearly
(or right-linearly) independent over ${\cal A}_r$.
\par For each ${\cal A}_r$ with $r\ge 2$, $r\in \bf N$,
there are the identities:
$z+z^*=2w_1$, $s(zs^*)=z^* + 2w_ss$ for each $s\in \hat b$, where
$z=\sum_{s\in \bf b}w_ss$, $w_s\in \bf R$ for each $s\in {\bf b}:=
\{ 1,i_1,...,i_{2^r-1} \} $, ${\hat b}:={\bf b}\setminus \{ 1 \} $,
hence $z^*=(2^r-2)^{-1} \{ -z +\sum_{s\in \hat b}s(zs^*) \} $
for each $r\ge 2$, $r\in \bf N$.
Therefore, $z^*$ does not play so special
role in ${\cal A}_r$, $r\ge 2$, as for $\bf C$.
\par In view of noncommutativity of ${\cal A}_r$, $r\ge 3$,
and the identities in it caused by the conjugation $(15)$, multiplication
$(3)$ and addition laws, for example, Equations $(7,8,9)$
for octonions, also Equations $(10,11,14)$ and Conditions $(12,13)$
in the general case of ${\cal A}_r$
a polynomial function $P: U\to {\cal A}_r$ in variables $z$ and
$z^{-1}$ may have several different representations
$$\mbox{\v{P}}(z)=\sum_{k,q(m)} \{
b_{k,1}z^{k_1}...b_{k,m}z^{k_m} \} _{q(m)},$$
where $b_{k,j}\in {\cal A}_r$ are constants, $k=(k_1,...,k_m)$,
$m\in \bf N$, $k_j=(k_{j,1},...,k_{j,n})$, $k_{j,l}\in \bf Z$,
$z^{k_j}:=\mbox{ }^1z^{k_{j,1}}... \mbox{ }^nz^{k_{j,n}},$
$\mbox{ }^lz^0:=1$, $U$ is an open subset of ${\cal A}_r^n$.
Certainly, we can consider $z-z_0$ instead of $z$ in the formula
of $\mbox{\v{P}}(z)$ on the right side, when a marked point $z_0$
is given.
In view of the nonassociativity of ${\cal A}_r$
here is used the notation $\{ a_1...a_m \} _{q(m)}$ of the product
of elements $a_1,...,a_m\in {\cal A}_r$ corresponding to the order
of products in this term defined by the position of brackets
$q(m):= (q_m,...,q_3)$, where $a_v:=(b_{k,v}z^{k_v})$
for each $v=1,...,m$, $q_m\in \bf N$ means that the first
(the most internal bracket) corresponds to the multiplication
$a_{q_m}a_{q_m+1}$ such that to the situation $(a_1...(a_ta_{t+1})...
(a_wa_{w+1})...a_m)$ with formally two simultaneous independent
multiplications, but $t<w$ by our definition of ordering there
corresponds $q_m=t$. After the first multiplication we get the
product of ${a'}_1,...,{a'}_{m-1}\in {\cal A}_r$, where not less, than
$m-2$ of these elements are the same as in the preceding term,
then $q_{m-1}$ corresponds to the first multiplication in this new term.
We omit $q_2$ and $q_1$, since they are unique.
Each term  $\{ b_{k,1}z^{k_1}...b_{k,m}z^{k_m} \} _{q(m)}=:
\omega (b_k,z)\ne 0$ we consider as a word of length
$\xi (\omega )=\sum_{j,l}\delta (k_{j,l}) +\sum_j\kappa (b_{k,j})$,
where $\delta (k_{j,l})=0$ for $k_{j,l}=0$ and $\delta (k_{j,l})=1$
for $k_{j,l}\ne 0$, $\kappa (b_{k,j})=j$ for $b_{k,j}=1$,
$\kappa (b_{k,j})=j+1$ for $b_{k,j}\in {\cal A}_r\setminus \{ 0,1 \} $.
A polynomial $P$ is considered as a phrase $\mbox{\v{P}}$ of a length
$\xi (\mbox{\v{P}}):=\sum_k\xi (\omega (b_k,z))$.
Using multiplication of constants in ${\cal A}_r$, commutativity
of $v\in \bf R$ with each $\mbox{ }^lz$ and $\mbox{ }^l\tilde z$,
and $\mbox{ }^lz^a\mbox{ }^lz^b=\mbox{ }^lz^{a+b}$
and $\mbox{ }^l{\tilde z}^a\mbox{ }^l{\tilde z}^b=
\mbox{ }^l{\tilde z}^{a+b}$, $\mbox{ }^lz\mbox{ }^l{\tilde z}=
\mbox{ }^l{\tilde z}\mbox{ }^lz$, it is possible to consider
representations of $P$ as phrases $\mbox{\v{P}}$ of a minimal lenght
$\xi (\mbox{\v{P}})$, then order them lexicographically by vectors
$q(m)$. We choose one such $\mbox{\v{P}}$ of a minimal lenght
and then minimal with respect to the lexicographic ordering of $q(m)$.
In view of the commutativity of the addition for terms $\{
a_1...a_m \} _{q(m)}$ and $\{ {a'}_1...{a'}_m \} _{q'(m)} $ of equal length
and having different vectors $q(m)$ and $q'(m)$ the order
of $q(m)$ and $q'(m)$ for $\mbox{\v{P}}$
is not important such that the minimality
is tested by all orderings of $q(m)$ among all terms of a given lenght
in $\mbox{\v{P}}$.
\par If $f: U\to {\cal A}_r$ is a function presented by a convergent by
$z$ series $f(z)=\sum_nP_n(z)$,
where $P_n(vz)=v^nP_n(z)$ for each $v\in \bf R$
is a $\bf R$-homogeneous polynomial, $n\in \bf Z$,
then we consider among all representations of $f$
such for which $\xi (\mbox{\v{P}}_n)$ is minimal for each $n\in \bf Z$.
This serves us to find representatives in classes of equivalent elements
of the $\bf R$-algebra of all polynomials on $U$ and $z$-analytic
functions. The corresponding family of locally $z$-analytic functions on $U$
is denoted by $C^{\omega }_z(U,{\cal A}_r)$ or $C^{\omega }(U,{\cal A}_r)$.
Each element of $C^{\omega }_z(U,{\cal A}_r)$ by our definition
is a unique phrase which may be infinite. We do not exclude
a possibility that two different phrases $f$ and $g$ may have the same
set-theoretical graph $\Gamma _f:= \{ (z,f(z)): z\in U \} $ as mappings
from $U$ into ${\cal A}_r$, for example, when a class of equivalence
defined by a graph has nonunique element of minimal length.
If each $P_n$ for $f$ has a decomposition of a particular left type
$$\mbox{\v{P}}(z)=\sum_kb_k(z^k),$$
where $0\le k\in \bf Z$, $b_k\in {\cal A}_r$, then
the space of all such locally analytic functions on $U$
is denoted by $\mbox{ }_lC^{\omega }(U,{\cal A}_r).$
The space of locally analytic functions $f$ having right type
decompositions for each $P_n$
$$\mbox{\v{P}}(z)=\sum_k(z^k)b_k$$
is denoted by $\mbox{ }_rC^{\omega }(U,{\cal A}_r)$.
The corresponding space in variables $(z,{\tilde z})$ is denoted by
$C^{\omega }_{z,\tilde z}(U,{\cal A}_r)$ and in variables
$\tilde z$ by $C^{\omega }_{\tilde z}(U,{\cal A}_r)$,
where $C^{\omega }_{z,\tilde z}:=C^{\omega }_{\mbox{ }_1z,\mbox{ }_2z}
(U^2,{\cal A}_r)|_{\mbox{ }_1z=z,\mbox{ }_2z=\tilde z} $,
$\mbox{ }_1z$ and $\mbox{ }_2z\in U$.
By our definition each element of $C^{\omega }_{z,\tilde z}(U,{\cal A}_r)$
is a unique phrase which may be infinite.
\par The $\bf R$-linear space $C^{\omega }_{z,\tilde z}(U,{\cal A}_r)$
is dense in the $\bf R$-linear space $C^0(U,{\cal A}_r)$
of all continuous functions $f: U\to {\cal A}_r$. We denote by
$C^0_z(U,{\cal A}_r)$ the $\bf R$-linear space of all equivalence classes
of Cauchy sequences from $C^{\omega }_z(U,{\cal A}_r)$ converging
relative to the $C^0$-uniformity. Analogously we define $C^0_{\tilde z}
(U,{\cal A}_r)$ and $C^0_{z,\tilde z}(U,{\cal A}_r)$.
\par {\bf 2.1.1. Definition.} If ${\cal G}\in C^0_z(U,{\cal A}_r)$,
then we say that $\cal G$ is $z$-represented.
Elements of $C^0_{\tilde z}$ (or
$C^0_{z,\tilde z}$) we call $\tilde z$- (or $(z,{\tilde z})-$ respectively)
represented functions. If $f: U\to {\cal A}_r$ is a (set-theoretic)
continuous function and ${\cal G}\subset {\cal F}\in C^0_{z,\tilde z}
(U,{\cal A}_r)$ such that each Cauchy sequence $\{ \zeta _n:
n\in {\bf N} \} $ from a family $\cal G$ (of converging Cauchy sequences)
converges to $f$ relative to the $C^0$-uniformity, then we call each
$g\in {\cal G}={\cal G}(f)$ an algebraic continuous function $g$.
(Since $\cal F$ is an equivalence class, then each $\{ \zeta _n: n \}
\in \cal G$ converges to the same limit $f$).
\par We say that $\cal G$ posses property $A$, if each $\{ \zeta _n:
n \} \in \cal G$ posseses property $A$. We say that $f$ posseses
property $A$, if there exists ${\cal G}(f)$ possesing
property $A$ and such that ${\cal G}(f)=\cal H$, where either
${\cal H}\in C^0_z(U,{\cal A}_r)$ or ${\cal H}\in C^0_{\tilde z}
(U,{\cal A}_r)$, or ${\cal H}\in C^0_{z,\tilde z}(U,{\cal A}_r)$.
If ${\cal G}\ne \emptyset $, ${\cal G}\subset {\cal H}$ and
${\cal G}\ne \cal H$ we will talk about $(A,{\cal G})$ property.
If $f$ is of higher class of smoothness $C^n$, $C^{\infty }$
or $C^{\omega }$, etc., then we take intersections of ${\cal G}(f)$
and $\cal H$ with $C^n$ or $C^{\infty }$, or $C^{\omega }$, etc.
supposing convergence relative to the respective uniformity
and such that ${\cal G}(f)\cap C^n$ or ${\cal G}(f)\cap C^{\infty }$,
or ${\cal G}(f)\cap C^{\omega }$, etc. is nonvoid.
Writing arguments $(\mbox{ }_1z,...,\mbox{ }_nz)$ of $f$  we will
outline the situation in each case indicating a subset of variables
by which property $A$ is accomplished. We may write for short
$f(z)$ or $f$ instead of $f(z,{\tilde z})$ in situations, in which
it can not cause a confusion. Our general supposition is that a
continuous function has a $(z,{\tilde z})$-representation, if another
will not be specified.
\par {\bf 2.1.2. Proposition.} {\it Let $A$ be an ${\cal A}_r$-additive,
${\bf R}$-homogeneous operator $A: {\cal A}_r^n\to {\cal A}_r^n$,
$r\ge 2$, $n\ge 1$, 
then $A$ is $\bf R$-linear and there exists a finite family
$A_j$ of right-${\cal A}_r$-linear and $B_j$ of left-${\cal A}_r$-linear
operators $j\in \{1,2,...,2^r \} $ independent of $h\in {\cal A}_r^n$
such that $A(h)=\sum_jA_j(hB_j)$ for each $h\in {\cal A}_r^n$, where
we write $A_j(h)=:A_jh$ and $B_j(h)=:hB_j$.}
\par {\bf Proof.} The first statement is evident.
To prove the second mention that $A(h)=\sum_{j=0}^{2^r-1}
A(i_j)h_j$, where $h_j\in \bf R^{n2^r}$, $h=\sum_{j=0}^{2^r-1}h_ji_j$,
$A(i_j)$ are $\bf R$-linear operators independent of $h_j$.
On the other hand, \\
$h_j=(-hi_j+ i_j(2^r-2)^{-1} \{ -h
+\sum_{j=1}^{2^r-1}i_j(hi_j^*) \} )/2$ for each $j=1,2,...,2^r-1$, \\
$h_0=(h+ (2^r-2)^{-1} \{ -h +\sum_{j=1}^{2^r-1}i_j(hi_j^*) \} )/2$. \\
Substituting these expressions of $h_j$, $j=0,1,...,2^r-1$, into
each term $A(i_j)h_j$ we get the second statement.
\par {\bf 2.2. Definition.} Consider
an open region $U$ in ${\cal A}_r^n$, $r\ge 3$,
the $n$-fold product of copies of ${\cal A}_r$, and let $f: U\to {\cal A}_r$
be a function. Then $f$ is said to be $z$-superdifferentiable
at a point $(\mbox{ }^1z,...,\mbox{ }^nz)=z\in  U$,
$\mbox{ }^1z,...,\mbox{ }^nz\in {\cal A}_r$, if it satisfies
Conditions $(2-7)$ below and if it can be written in the form
$$(1)\quad f(z+h)=f(z)+\sum_{j=1}^nA_j\mbox{ }^jh + \epsilon (h)|h|$$
for each $h\in {\cal A}_r^n$ such that $z+h\in U$, where
$A_j$ is an ${\cal A}_r$-valued ${\cal A}_r$-additive
$\bf R$-homogeneous operator of $h$-variable, in general
it is non-linear for each $j=1,...,n$ and $A_j$ is denoted
by $(Df(z)).e_j$ and there exists a derivative $f'(z)$ such that
a differential is given by
$$(2)\quad Df(z).h:=f'(z).h:=
\sum_{j=1}^n (\partial f(z)/\partial \mbox{ }^jz)\mbox{ }^jh,$$
where $\epsilon (h)$, $\epsilon : {\cal A}_r^n\to {\cal A}_r$,
is a function continuous at zero such that $\epsilon (0)=0$, \\
$e_j=(0,...,0,1,0,...,0)$ is the vector
in ${\cal A}_r^n$ with $1$ on $j$-th place,
$$(3)\quad Df(z).h=:(Df)(z;h)$$
such that $(Df)(z;h)$ is additive in
$h$ and $\bf R$-homogeneous, that is,
$$(4)\quad (Df)(z;h_1+h_2)=(Df)(h_1)+(Df)(h_2)\mbox{ and }
(Df)(z;vh)=v(Df)(z;h)$$
for each $h_1$, $h_2$ and $h\in {\cal A}_r^n$, $v\in \bf R$.
There are imposed conditions:
$$(5)\quad (\partial _zz).h=h, \quad \partial _z1=0,\quad
\partial _z{\tilde z}=0,\quad \partial _{\tilde z}z=0,\quad
(\partial _{\tilde z}{\tilde z}).h=\tilde h$$
$$\mbox{also } D=\partial _z + \partial _{\tilde z} , \quad
(D(fg)).h=((Df).h)g+f(Dg).h$$
for a product of two supedifferentiable functions $f$ and $g$
and each $h\in {\cal A}_r^n$, where the notation $\partial _z$
corresponds to $\partial /\partial z$ and $\partial _{\tilde z}$
corresponds to $\partial /\partial {\tilde z}$.
We also have distributivity laws relative to multiplication from the right
by elements $\lambda \in {\cal A}_r$:
$$(6)\quad (D(f+g))(z;h\lambda )=(Df)(z;h\lambda )+(Dg)(z;h\lambda ),$$
$$(Df)(z;h(\lambda _1+\lambda _2))=(Df)(z;h\lambda _1)+
(Df)(z;h\lambda _2)$$
for each superdifferentiable functions
$f$ and $g$ at $z$ and each $\lambda $, $\lambda _1$ and
$\lambda _2\in {\cal A}_r$.  There are also left distributive laws:
$$(7)\quad (D\lambda (f+g))(z;h)=\lambda (Df)(z;h)+\lambda (Dg)(z;h),$$
$$(D(\lambda _1+\lambda _2)f)(z;h)=\lambda _1(Df)(z;h)+
\lambda _2(Df)(z;h).$$
If use $(z,{\tilde z})$-representation of polynomials and functions, then
we define $(z,{\tilde z})$-superdifferentiability by the pair
$(z,{\tilde z})$, by $z$ and by $\tilde z$ such that \\
$(8)\quad D_z{\tilde z}=0$, $D_{\tilde z}z=0$, $(D_zz).h=h$,
$(D_{\tilde z}{\tilde z}).h=\tilde h$,\\
$(D_{z,\tilde z}(fg)).h=((D_{z,\tilde z}f).h)g+f(D_{z,\tilde z}g).h$ 
and $(D_{z,\tilde z}f).h=(D_zf).h+(D_{\tilde z}f).h$ \\
for each two $(z,{\tilde z})$-superdifferentiable functions
$f$ and $g$, each $h\in {\cal A}_r^n$ (see also \cite{luoyst}).
We take a function $g(\mbox{ }_1z,\mbox{ }_2z)$
in the $z$-representation by $\mbox{ }_1z$ and $\mbox{ }_2z$,
then consider the operator $D$ by the variable $(\mbox{ }_1z,\mbox{ }_2z)$
and in the expression
$(Dg(\mbox{ }_1z,\mbox{ }_2z)).h$, put for the components
$\mbox{ }_1z=z$, $\mbox{ }_2z=\tilde z$, $\mbox{ }_1h=\mbox{ }_2h=:
\alpha \in {\cal A}_r^n$ and consider the function $g(z,{\tilde z})=:f$,
where $z=(\mbox{ }^1z,...,\mbox{ }^nz)\in U\subset {\cal A}_r^n$,
${\tilde z}=(\mbox{ }^1{\tilde z},...,\mbox{ }^n{\tilde z})$,
$az:=(a\mbox{ }^1z,...,a\mbox{ }^nz)$,
$zb:=(\mbox{ }^1zb,...,\mbox{ }^nzb)$ for each $a, b\in {\cal A}_r$.
\par  If there is a function $g(\mbox{ }_1z,\mbox{ }_2z)$
on an open subset $W$ in ${\cal A}_r^{2n}$ with values in ${\cal A}_r$ 
$(\mbox{ }_1z,\mbox{ }_2z)$-superdifferentiable at a point
$(\mbox{ }_1z,\mbox{ }_2z)$,
$\mbox{ }_1z$ and $\mbox{ }_2z\in {\cal A}_r^n$, also
\par $g(\mbox{ }_1z,\mbox{ }_2z)|_{\mbox{ }_1z=z, \mbox{ }_2z=\tilde z}
=:f(z,{\tilde z})$, $z=\xi $, \\
then we say that $f$ is $(z,{\tilde z})$-superdifferentiable at a point
$\xi $ and
\par $(9)\quad (D_zf(z,{\tilde z})).h=(\partial f(z,{\tilde z})/\partial z).h:=
\{ (Dg(\mbox{ }_1z,\mbox{ }_2z)).(h,0) \}
|_{\mbox{ }_1z=z, \mbox{ }_2z=\tilde z}$,
\par $(D_{\tilde z}f(z,{\tilde z})).h=
(\partial f(z,{\tilde z})/\partial {\tilde z}).h:=
\{ (Dg(\mbox{ }_1z,\mbox{ }_2z)).(0,h) \}
|_{\mbox{ }_1z=z, \mbox{ }_2z=\tilde z}$, \\
where $h\in {\cal A}_r^n$ and $f$ is supposed to be defined
by $g$ and its restriction on $\{ (\mbox{ }_1z, \mbox{ }_2z)\in W:$
$\mbox{ }_2z=(\mbox{ }_1z)^{\tilde .} \} $, \\
$(10)\quad D_{\mbox{ }_1z}g(\mbox{ }_1z, \mbox{ }_2z).h:=
D_{(\mbox{ }_1z, \mbox{ }_2z)}g(\mbox{ }_1z, \mbox{ }_2z).(h,0),$ \\ 
$\quad D_{\mbox{ }_2z}g(\mbox{ }_1z, \mbox{ }_2z).h:=
D_{(\mbox{ }_1z, \mbox{ }_2z)}g(\mbox{ }_1z, \mbox{ }_2z).(0,h),$ \\
$(D_zf(z,{\tilde z})).e_j=:\partial f(z,{\tilde z})/\partial \mbox{ }^jz,$
$(D_{\tilde z}f(z,{\tilde z})).e_j=:\partial f(z,{\tilde z})/
\partial \mbox{ }^j{\tilde z}.$ \\
Since the Cayley-Dickson
algebras are over $\bf R$ and Fr\'echet differentials are unique,
then for functions $g: W\to {\cal A}_r$ 
and $f: U\to {\cal A}_r$, their superdifferentials
$(Dg).h$ and $(D_{(z,{\tilde z})}f).\alpha $
are unique, so we have $D_z=D_{\mbox{ }_1z}|_{\mbox{ }_1z=z}$,
$D_{\tilde z}=D_{\mbox{ }_2z=\tilde z}$ in the
$(z,{\tilde z})$-representation, where
$U$ is open in ${\cal A}_r^n$ such that
$\{ (\mbox{ }_1z=z, \mbox{ }_2z={\tilde z}): z\in U \} \subset W$.
In particular, if there are functions $f_1, f_2, f_3$ such that
$f_3=f_1(z,{\tilde z})f_2(z,{\tilde z})$, $f_j=g_j(\mbox{ }_1z,
\mbox{ }_2z)|_{\mbox{ }_1z=z, \mbox{ }_2z=\tilde z}$, $j=1,2,3$,
where either $g_j$ for each $j$ is presented by a minimized series
of \S 2.1, or while multiplication $(f_1,f_2)\mapsto f_1f_2$
no any reorganization of their series, for example, by minimality
is made (this is the case of our definition above), then
\par $(11)\quad (D_zf_1f_2).h=((D_zf_1).h)f_2+f_1(D_zf_2).h$ and
\par $(D_{\tilde z}f_1f_2).h=((D_{\tilde z}f_1).h)f_2+f_1(D_{\tilde z}f_2).h$
\\
for each $h\in {\cal A}_r^n$, since $D_{z,\tilde z}=D_z+D_{\tilde z}$.
Generally, \\
$(D_{z,\tilde z}f_1f_2).h=((D_{z,\tilde z}f_1).h)f_2+
f_1(D_{z,\tilde z}f_2).h$ for each $(z,{\tilde z})$-superdifferentiable
functions $f_1$ and $f_2$ on $U$ and each $h\in {\cal A}_r^n$.
\par A function $f: U\to {\cal A}_r$ is called $\tilde z$-superdifferentiable
at a point $\xi $, if there exists a function $g: U\to {\cal A}_r$ such that
$g({\tilde z})=f(z)$ and $g(z)$ is $z$-superdifferentiable at $\xi $.
\par {\bf Notation.} We may write a function $f(z)$ with $z\in {\cal A}_r$,
$r\ge 2$, in variables $((\mbox{ }^jw_s: s\in {\bf b}): j=1,...,n )$,
${\bf b}:={\bf b}_r$, as $F((\mbox{ }^jw_s: s\in {\bf b}): j=1,...,n)=
f\circ \sigma ((\mbox{ }^jw_s: s\in {\bf b}): j=1,...,n)$, where
$\sigma ((\mbox{ }^jw_s: s\in {\bf b}): j=1,...,n)=
(\mbox{ }^jz: j=1,...,n)$ is a bijective mapping.
For $U$ open in ${\cal A}_r^n$ and $F: U\to {\cal A}_r$
we can write $F$ in the form $F=\sum_{s\in {\bf b}}F_ss$,
where $F_s\in \bf R$ for each $s\in \bf b$,
$F_{vs}:=vF_s$ for each $v\in \bf R$.
\par {\bf 2.2.1. Proposition.} {\it Let $g: U\to {\cal A}_r^m$, $r\ge 3$,
and $f: W\to {\cal A}_r^n$ be two superdifferentiable functions
on $U$ and $W$ respectively such that $g(U)\supset W$,
$U$ is open in ${\cal A}_r^k$, $W$ is open in ${\cal A}_r^m$,
$k, n, m\in \bf N$. Then the composite function $f\circ g(z):=
f(g(z))$ is superdifferentiable on $V:=g^{-1}(W)$ and
\par $(Df\circ g(z)).h=(Df(g)).((Dg(z)).h)$ \\
for each $z\in V$ and each $h\in {\cal A}_r^k$, where
$f$ and $g$ are simultaneously $(z,{\tilde z})$, or $z$, or
$\tilde z$-superdifferentiable and hence $f\circ g$ is of the same
type of superdifferentiability.}
\par {\bf Proof.} Since $g$ is superdifferentiable, then
$g$ is continuous and $g^{-1}(W)$ is open in ${\cal A}_r^k$.
Then $f\circ g(z+h)-f\circ g(z)=(Df(g))|_{g=g(z)}.(g(z+h)-g(z))+
\epsilon _f(\eta )|\eta |$, where $\eta =g(z+h)-g(z)$,
$g(z+h)-g(z)=(Dg(z)).h+\epsilon _g(h)|h|$ (see \S 2.2).
Since $Df$ is ${\cal A}_r^m$-additive and $\bf R$-homogeneous
(and continuous) operator on ${\cal A}_r^m$, then \\
$f\circ g(z+h)-f\circ g(z)=(Df(g))|_{g=g(z)}.((Dg(z)).h)
+\epsilon _{f\circ g}(h)|h|$, where \\
$\epsilon _{f\circ g}(h)|h|:=\epsilon _f((Dg(z)).h+\epsilon _g(h)|h|)
|(Dg(z)).h+\epsilon _g(h)|h||$ \\
$+[(Df(g))|_{g=g(z)}.(\epsilon _g(h))]|h|)$, \\
$|(Dg(z)).h+\epsilon _g(h)|h||\le 
[\| Dg(z) \| + |\epsilon _g(h)|] |h|$, hence \\
$|\epsilon _{f\circ g}(h)|\le |\epsilon _f((Dg(z)).h+\epsilon _g(h)|h|)|
[\| Dg(z) \| $  $+ |\epsilon _g(h)|] +
\| (Df(g))|_{g=g(z)} \| |\epsilon _g(h)|$ \\
and inevitably $\lim_{h\to 0} \epsilon _{f\circ g}(h)=0$.
Moreover, $\epsilon _{f\circ g}(h)$ is continuous in $h$, since
$\epsilon _g$ and $\epsilon _f$ are continuous functions,
$Df$ and $Dg$ are continuous operators. Evidently, if
$\partial _{\tilde z}f=0$ and $\partial _{\tilde z}g=0$ on
domains of $f$ and $g$ respectively, then $\partial _{\tilde z}f\circ g=0$
on $V$, since $D=\partial _z+\partial _{\tilde z}$.
\par {\bf 2.3. Proposition.} {\it A function $f: U\to {\cal A}_r$ is
$z$-superdifferentiable at a point $a\in U$ if and only if
$F$ is Frech\'et differentiable at $a$ and
$\partial _{\tilde z}f(z)|_{z=a}=0$. If $f$ is $z$-superdifferentiable
on $U$, then $f$ is $z$-represented on $U$.
If $f'(a)$ is right superlinear on the superalgebra
${\cal A}_r^n$, then $f$ is $z$-superdifferentiable at $a\in U$
if and only if $F$ is Frech\'et differentiable at $a\in U$ and satisfies
the following equations:
$$(1)\quad (\partial F_{ps}/\partial \mbox{ }^jw_p)=
((ps)p^*)^*((qs)q^*) (\partial F_{qs}/\partial \mbox{ }^jw_q),
\quad \mbox{for each } p, q, s \in {\bf b} $$
or shortly:
$$(2)\quad \partial F/\partial \mbox{ }^jw_1=
(\partial F/\partial \mbox{ }^jw_q)q^*$$
for each $q\in {\hat b}_r$ and each $j=1,...,n$.
A $(z,{\tilde z})$-superdifferentiable function $f$ at $a\in U$
is $z$-superdifferentiable at $a\in U$ if and only if
$D_{\tilde z}f(z,{\tilde z})|_{z=a}=0$.}
\par {\bf Proof.} For each canonical closed compact set $U$ in ${\cal A}_r$
the set of all polynomial by $z$
functions is dense in the space of all continuous on $U$
Frech\'et differentiable functions on $Int (U)$.
\par As usually a set $A$ having structure of an additive group
and having distributive multiplications of its elements
on Cayley-Dickson numbers $z\in {\cal A}_v$ from the left
and from the right is called a vector space over ${\cal A}_v$.
In such sence it is the $\bf R$-linear space and also
the left and right module over ${\cal A}_v$.
For two vector spaces $A$ and $B$ over ${\cal A}_v$
consider their ordered tensor product $A\otimes B$ over ${\cal A}_v$
consisting of elements $a\otimes b := (a,b)$ such that $a\in A$ and
$b\in B$, $\alpha (a,b)=(\alpha a,b)$ and $(a,b)\beta = (a,b\beta )$
for each $\alpha , \beta \in {\cal A}_v$, $(a_1\otimes b_1)(a_2\otimes b_2)=
a_1a_2\otimes b_1b_2$ for each $a_1, a_2\in A$ and $b_1, b_2\in B$.
In the aforementioned respect
$A\otimes B$ is the $\bf R$-linear space
and at the same time left and right module over ${\cal A}_v$.
Then $A\otimes B$ has the structure of the
vector space over ${\cal A}_v$. By induction consider
tensor products $\{ C_1\otimes C_2\otimes ... \otimes C_n \} _{q(n)}$,
where $C_1,...,C_n \in \{ A, B \} $, $q(n)$ indicates on the order
of tensor multiplications in $ \{ * \} $.
For two ${\cal A}_v$-vector spaces $V$ and $W$ their direct sum
$V\oplus W$ is the ${\cal A}_v$-vector space consisting of all
elements $(a,b)$ with $a\in V$ and $b\in W$ such that $\alpha (a,b)=
(\alpha a,\alpha b)$ and $(a,b)\beta =(a\beta ,b\beta )$ for each
$\alpha $ and $\beta \in {\cal A}_v$.
Therefore, the direct sum of all different
tensor products $\{ C_1\otimes C_2\otimes ... \otimes C_n \} _{q(n)}$,
which are $\bf R$-linear spaces and left and right modules over
${\cal A}_v$, provides the minimal tensor space $T(A,B)$ generated
by $A$ and $B$.
\par Operators $\partial _z$ and $\partial _{\tilde z}$ are uniquely defined
on $C^{\omega }_z(U,{\cal A}_r)$ and $C^{\omega }_{\tilde z}
(U,{\cal A}_r)$, hence they are unique on the tensor space
$T(C^{\omega }_z(U,{\cal A}_r),C^{\omega }_{\tilde z}(U,{\cal A}_r))$,
which is dense in $C^{\omega }_{z,\tilde z}(U,{\cal A}_r)$,
since $C^{\omega }_{z,\tilde z}(U,{\cal A}_r):=
C^{\omega }_{\mbox{ }_1z,\mbox{ }_2z}(U^2,{\cal A}_r)
|_{\mbox{ }_1z=z,\mbox{ }_2z=\tilde z}$. Therefore, operators
$\partial _z$ and $\partial _{\tilde z}$ are uniquely defined on
$C^{\omega }_{z,\tilde z}(U,{\cal A}_r)$.
\par If there is a product $fg$ of two phrases $f$ and $g$ from
$C^{\omega }_{z,\tilde z}(U,{\cal A}_r)$, then if it is reduced
to a minimal phrase $\xi $, then it is made with the help of
$z^nz^m=z^{n+m}$ and ${\tilde z}^n{\tilde z}^m={\tilde z}^{n+m}$
and identities for constants in ${\cal A}_r$, since no any
shortening related with their permutation $z{\tilde z}={\tilde z}z$
or substitution of $z$ on $\tilde z$ or $\tilde z$ on $z$, for example,
using identity ${\tilde z}=l(zl^*)$ is not allowed in $C^{\omega }_{z,
\tilde z}(U,{\cal A}_r)$ in accordance with our convention in \S 2.1, since
$C^{\omega }_{z,\tilde z}(U,{\cal A}_r):=C^{\omega }_{\mbox{ }_1z,
\mbox{ }_2z}(U^2,{\cal A}_r)|_{\mbox{ }_1z=z,\mbox{ }_2z=\tilde z}$
and in $C^{\omega }_{\mbox{ }_1z,\mbox{ }_2z}(U^2,{\cal A}_r)$ variables
$\mbox{ }_1z$ and $\mbox{ }_2z$ do not commute, $\mbox{ }_1z$ and
$\mbox{ }_2z$ are different variables which are not related.
Therefeore, $\partial _z\xi .h=(\partial _zf.h)g+f(\partial _zg.h)$
and $\partial _{\tilde z}\xi .h=(\partial _{\tilde z}f.h)g+
f(\partial _{\tilde z}g.h)$, hence $\partial _z$ and
$\partial _{\tilde z}$ are correctly defined.
In particular functions of the form of series $f=\sum
\{ \mbox{ }_{l_1}f... \mbox{ }_{l_t}f \} _{q(t)}$
converging on $U$ together with its superdifferential on $Int (U)$
such that each $\mbox{ }_lf$ is superlinearly
$z$-superdifferentiable on $Int (U)$ relative to the superalgebra
${\cal A}_r$ is dense in the
$\bf R$-linear space of $z$-superdifferentiable functions $g$ on $U$,
since $(Dg(z)).h$ is continuous by $(z,h)$.
We can use $\delta $-approximation for each $\delta >0$
of $Dg(z).h$ on a sufficiently small open subset $V$ in $U$
such that $z\in V$ by functions $\zeta _n$ polynomial in $z$
and $\bf R$-homogeneous ${\cal A}_r^n$-additive in $h$ and partition of
unity in $U$ by $C^{\omega }_z$-functions and then consider
functions $\xi _n$ with ${\xi _n}'$ corresponding to $\zeta _n$, since
$\epsilon (h)$ is continuous at $0$ for each $z\in U$ and for each
canonical closed compact subset $W$ in $U$ from each open covering
we can choose a finite subcovering of $W$.
\par From Conditions $2.2.(2-7)$ it follows, that the $z$-superdifferentiability
conditions are defined uniquely on space of polynomials.
In view of Conditions $2.2.(1-7)$ the $z$-superdifferentiability of a
polynomial or a converging series $P$ on $U$
means that it is expressible through a sum or a converging series
of products of $\mbox{ }^jz$ and constants from ${\cal A}_r$.
Therefore, each $z$-superdifferentiable function $f$ on $U$
is the equivalence class of all Cauchy sequences from $C^{\omega }_z(U,
{\cal A}_r)$ converging to $f$ relative to the $C^1$-uniformity, since
$(D*).h: C^1\to C^0$ is continuous for each $h\in {\cal A}_r^n$.
\par Suppose that $f$ is $z$-superdifferentiable at a point $a$.
To each $f'(z)$ there corresponds a $\bf R$-linear
operator on the Euclidean space $\bf R^{2^rn}$. Moreover,
we have the distributivity and associativity
laws for $(Df)(z;h)$ relative to the right multiplication
on elements $\lambda \in {\cal A}_r$ (see \S 2.1, 2.2).
Then $f(a+h)-f(a)=\partial _af(a).h+ \epsilon (h)|h|$ and
$\partial _{\tilde z}f(z)|_{z=a}=0$, since generally
$f(a+h)-f(a)=(\partial _af(a)).h+(\partial _{\tilde a}f(a)).h
+\epsilon (h)|h|$,
where $\epsilon (h)$ is continuous by $h$ and $\epsilon (0)=0$.
Vice versa, if $F$ is Frech\'et differentiable and
$\partial _{\tilde z}f(z)|_{z=a}=0$, then
expressing $w_ss$ for each $s\in {\bf b}_r$ through linear combinations
of $z$ (with multiplication on constants from ${\cal A}_r$ on the left
and on the right) with constant coefficients we get
the increment of $f$ as above.
\par Consider now the particular case, when
$f'$ is right superlinear on the superalgebra ${\cal A}_r^n$
and $\partial _{\tilde z}f(z)|_{z=a}=0$.
In this case $f'(a)$ is right ${\cal A}_r$-linear.
Using the definition of the $z$-superderivative and that there
is a bijective correspondence
between $z$ and $((\mbox{ }^jw_s: s\in {\bf b}_r): j=1,...,n) $,
$\mbox{ }^jw_s\in \bf R$ for each $s\in {\bf b}_r$, $j=1,...,n$,
we consider a function $f=f(z)$ and $F((\mbox{ }^jw_s: s\in {\bf b}_r):
j=1,...,n)=f\circ \sigma $, where $f$ is $z$-superdifferentiable by
$z$, hence $F$ is Frech\'et differentiable by
$((\mbox{ }^jw_s: s\in {\bf b}_r): j=1,...,n)$
and we obtain the expressions:
$$\partial F/\partial \mbox{ }^jw_s=(\partial F/\partial \mbox{ }^jz).
(\partial \mbox{ }^jz/\partial \mbox{ }^jw_s),$$
since $\partial \mbox{ }^jz/\partial \mbox{ }^kw_s=0$
and $\partial \mbox{ }^j{\tilde z} / \partial \mbox{ }^kw_s=0$
and $\partial f(z)/\partial \mbox{ }^j{\tilde z}|_{z=a}=0$
for each $k\ne j$. From $\partial \mbox{ }^jz/\partial \mbox{ }^jw_s=s$
for each $s\in {\bf b}_r$ we get Equations $(2)$, since
$ps=-sp$, $F_{ps}=-F_{sp}$ and $pp^*=1$ for each $p\ne s\in \hat b$.
Using the equality $F=\sum_{s\in \bf b}F_ss$ we get
Equations $(1)$ from the latter equations, since
$qf$ is right superlinearly superdifferentiable together
with $f$ for each $q\in \bf b$, $((ps)p^*)^*((qs)q^*)\in \{ -1, 1 \}
\subset \bf R$ for each $p, q, s\in \bf b$.
\par Let now $F$ be Frech\'et differentiable at $a$ and
let $F$ be satisfying Conditions $(1)$.
Then $f(z)-f(a)=\sum_{j=1}^n \sum_{s\in \bf b}
(\partial F/\partial \mbox{ }^jw_s)\Delta \mbox{ }^jw_s$
$+\epsilon (z-a)|z-a| $,
where $\Delta (\mbox{ }^jw_s: s\in {\bf b})
=\sigma ^{-1}(\mbox{ }^jz)-\sigma
(\mbox{ }^ja)$ for each $j=1,...,n$.
From Conditions $(2)$ equivalent to $(1)$ we get  \\
$f(z)-f(a)= \sum_{j=1}^n \sum_{s\in \bf b}
(\partial F/\partial \mbox{ }^jw_1)s\Delta \mbox{ }^jw_s$
$+\epsilon (z-a)|z-a| =$ \\
$\sum_{j=1}^n(\partial F/\partial \mbox{ }^jw_1) \Delta \mbox{ }^jz$
$+\epsilon (z-a)|z-a| ,$ \\
where $\epsilon $ is a function continuous at $0$ and $\epsilon (0)=0$.
Therefore, $f$ is superdifferentiable by $z$ at $a$ such that
$f'(a)$ is right superlinear, since $\partial F/\partial \mbox{ }^jw_s$
are real matrices and hence $f'(a).(\mbox{ }_1h\lambda _1+
\mbox{ }_2h\lambda _2)=(f'(a)\mbox{ }_1h)\lambda _1+
(f'(a)\mbox{ }_2h)\lambda _2$ for each $\lambda _1$ and $\lambda _2\in
{\cal A}_r$ and each $\mbox{ }_1h$ and $\mbox{ }_2h\in {\cal A}_r^n$.
\par The last statement of this proposition follows from Definition 2.2.
\par {\bf 2.3.1. Notation.} If $f: U\to {\cal A}_r$ is either
$z$-superdifferentiable or $\tilde z$-superdifferentiable at $a\in U$
or on $U$, then we can write also $D_{\tilde z}$ instead of
$\partial _{\tilde z}$ and $D_z$ instead of $\partial _z$
at $a\in U$ or on $U$ respectively in situations, when it can not
cause a confusion, where $U$ is open in ${\cal A}_r^n$.
\par {\bf 2.4. Corollary.} {\it Let $f$ be a continuously
superdifferentiable function by $z$
with a right superlinear superdifferential on the superalgebra
${\cal A}_r^n$, $r\ge 3$, in an open subset $U$ in ${\cal A}_r^n$
and let $F$ be twice continuously differentiable by $((\mbox{ }^jw_s:
s\in {\bf b}_r): j=1,...,n)$
in $U$, then each component $F_s$ of $F$ is the harmonic functions
by pairs of variables $(\mbox{ }^jw_p,\mbox{ }^jw_q)$
for each $p\ne q\in {\bf b}_r$ namely:
$$(1)\quad \bigtriangleup _{\mbox{ }^jw_p,\mbox{ }^jw_q}F_s=0,$$
for each $j=1,..,n,$ where
$\bigtriangleup _{\mbox{ }^jw_p,\mbox{ }^jw_q}F_s:=
\partial ^2F_s/\partial \mbox{ }^jw_p^2+
\partial ^2F_s/\partial \mbox{ }^jw_q^2$.}
\par {\bf Proof.} From Equations $2.3.(1)$
and in view of the twice continuous differentiability of $F$
it follows, that $(\partial ^2F_s/\partial \mbox{ }^jw_p^2)=
(\partial ^2F_{(sp^*)q}/\partial \mbox{ }^jw_p\partial \mbox{ }^jw_q)=$
$\partial ^2F_{(((sp^*)q)p^*)q}/\partial \mbox{ }^jw_q^2)=$
$-(\partial ^2F_s/\partial \mbox{ }^jw_q^2)$,
since $F_{vs}=vF_s$ for each $v\in \bf R$ and each $s\in {\bf b}_r$,
$p\ne q\in {\bf b}_r$ and hence $p^*q\in {\hat b}_r$,
$t^2=-1$ for each $t\in {\hat b}_r$, $p^*=-p$ for each $p\in {\hat b}_r$,
$pq=-qp$ for each $p\ne q\in {\hat b}_r$.
\par {\bf 2.5. Note and Definition.} Let $U$ be an open
subset in ${\cal A}_r$ and let $f: U\to {\cal A}_r$, $r\ge 3$,
be a function defined on $U$ such that
$$(i)\quad f(z,{\tilde z})= \{ f^1(z,{\tilde z})...
f^j(z,{\tilde z}) \} _{q(j)},$$
where each function $f^s(z,{\tilde z})$ is presented by a Laurent series
$$(ii)\quad f^s(z,{\tilde z})=\sum_{n=n_0}^{\infty } \sum_{m=m_0}^{\infty }
(f^s_{n,m} (z-\zeta )^n)({\tilde z}-{\tilde {\zeta }})^m$$
converging on $U$, where $f^s_{n,m}\in {\cal A}_r$, $z\in U$, $\zeta \in
{\cal A}_r$ is a marked point, $n$ and $m\in \bf Z$,
if $n_0<0$ or $m_0<0$, then $\zeta \notin U$. Consider the case
$f^s_{-1,m}=0$ for each $s$ and $m$. The case with terms
$f^s_{-1,m}\ne 0$ will be considered later.
\par Let $[a,b]$ be a segment in $\bf R$ and $\gamma : [a,b]\to {\cal A}_r$
be a continuous function. Consider a partitioning $P$ of
$[a,b]$, that is, $P$ is a finite subset of $[a,b]$
consisting of an increasing sequence of points
$a=c_0<...<c_k<c_{k+1}<...<c_t=b$, then the norm of $P$
is defined as $|P|:=\max _k(x_{k+1}-x_k)$ and the $P$-variation
of $\gamma $ as $v(\gamma ;P):=\sum_{k=0}^{t-1}
|\gamma (c_{k+1})-\gamma (c_k)|,$ where $t=t(P)\in \bf N$.
The total variation (or the length) of $\gamma $ is defined as
$V(\gamma )=\sup_Pv(\gamma ;P)$. Suppose that $\gamma $ is rectifiable,
that is, $V(\gamma )<\infty $. For $f$  having decomposition
$(2.5.i,ii)$ with
$f^s_{-1,m}=0$ for each $s$ and $m$ and a rectifiable path
$\gamma : [a,b]\to U$
we define a (noncommutative) Cayley-Dickson algebra line integral.
Consider more general situation.
\par Let $f: U\to {\cal A}_r$ be a continuous function, where $U$
is open in ${\cal A}_r$, $f$ is defined by a continuous function
$\xi : U^2\to {\cal A}_r$ such that
\par $(1)$ $\xi (\mbox{ }_1z,\mbox{ }_2z)|_{\mbox{ }_1z=z,
\mbox{ }_2z=\tilde z}=f(z,{\tilde z})$ \\
or shortly $f(z)$ instead of $f(z,{\tilde z})$, where
$\mbox{ }_1z$ and $\mbox{ }_2z\in U$. Let also $g: U^2\to {\cal A}_r$
be a continuous function, which is $\mbox{ }_1z$-superdifferentiable
such that
\par $(2)$ $(\partial g(\mbox{ }_1z,\mbox{ }_2z)/\partial \mbox{ }_1z).1=
\xi (\mbox{ }_1z,\mbox{ }_2z)$ on $U^2$. Then put
\par $(3)$ ${\hat f}(z,{\tilde z}).h:={\hat f}_z(z,{\tilde z}).h:=
[(\partial g(\mbox{ }_1z,\mbox{ }_2z)/
\partial \mbox{ }_1z).h]|_{\mbox{ }_1z=z,\mbox{ }_2z=\tilde z}$
for each $h\in {\cal A}_r$. Shortly we can write these as
$(\partial g(z,{\tilde z})/ \partial z).1=f(z,{\tilde z})$ and
${\hat f}_z(z,{\tilde z}).h:={\hat f}(z).h
:=(\partial g(z,{\tilde z})/ \partial z).h$.
If the following limit exists
$$(4)\quad \int_{\gamma }f(z,{\tilde z})dz:=
\lim_P I(f,\gamma ;P),\mbox{ where}$$
$$(5)\quad I(f,\gamma ;P):=\sum_{k=0}^{t-1}
{\hat f}(z_{k+1},{\tilde z}_{k+1}).(\Delta z_k),$$
where $\Delta z_k:=z_{k+1}-z_k$, $z_k:=\gamma (c_k)$ for each
$k=0,...,t$, then we say that $f$ is line integrable along $\gamma $
by the variable $z$. Analogously we define
$\int_{\gamma }f(z,{\tilde z})d{\tilde z}$
with $(\partial g(z,{\tilde z})/ \partial {\tilde z}).1=f(z,{\tilde z})$,
${\hat f}_{\tilde z}(z,{\tilde z}):=
(\partial g(z,{\tilde z})/\partial {\tilde z}).h$,
where $g(\mbox{ }_1z,\mbox{ }_2z)$ is $\mbox{ }_2z$-superdifferentiable.
\par {\bf Remark.} In view of Definitions $2.1, 2.2$ and Proposition
$2.3$ Conditions $2.5.(1-3)$ are correct, for example, the taking of
functions $\xi $ and $g$ in the $(\mbox{ }_1z,\mbox{ }_2z)$-representation
is sufficient.
\par This definition is justified by the following lemma and proposition.
\par {\bf 2.5.1. Lemma.} {\it Let $f: U\to {\cal A}_r$ satisfy Conditions
$2.2.(1,2,5)$, where $U$ is an open subset in ${\cal A}_r$. Then conditions
\par $(1)$ $\partial f(z)/\partial {\bar z}_{2j,2j+1}=0$ for each
$j=0,1,...,2^{r-1}-1$ and $z\in U$, \\
where $z_{2j,2j+1}:=w_s+w_ps^*p$ with $s=i_{2j}$, $p=i_{2j+1}$,
are equivalent with $\partial _{\tilde z}f(z)=0$ for each $z\in U$.}
\par {\bf Proof.}  Since $(z+h)^* - z^*=h^*$, $((\lambda h)^*)^*=\lambda h$
for each $z$ and $h\in {\cal A}_r^n$ and each $\lambda \in {\cal A}_r$,
$(h_1+h_pp)^*=h_1-h_pp$ for each $p\in \hat b$, $(s(h_s+h_ps^*p))^*=
(h_s-h_ps^*p)s^*$ for each $s\ne p\in \hat b$, $(\partial z/\partial z).h
=h$, $(\partial {\tilde z}/ \partial {\tilde z}).h={\tilde h}$,
$\partial z/\partial {\tilde z}=0$ and  $\partial {\tilde z}
/\partial z=0$, where $h_p\in \bf R$ for each $p\in \bf b$, hence \\
$(\partial f(z)/\partial z_{2j,2j+1}).h=
(\partial f(z)/\partial z).(sh)|_{h\in {\bf R}\oplus s^*p\bf R}$ and  \\
$(\partial f(z)/\partial {\bar z}_{2j,2j+1}).h=
(\partial f(z)/\partial {\tilde z}).(sh)|_{h\in {\bf R}\oplus s^*p\bf R}$,
where $s=i_{2j}$ and $p=i_{2j+1}$, then \\
$\partial z/\partial {\bar z}_{2j,2j+1}=0$,
$\partial {\tilde z}/\partial z_{2j,2j+1}=0$ and \\
$(\partial {\bar z}_{2j,2j+1}/\partial {\bar z}_{2j,2j+1}).h=\tilde h$,
$(\partial z_{2j,2j+1}/\partial z_{2j,2j+1}).h=h$ for each
$h\in {\bf R}\oplus s^*p\bf R$. In view of Proposition 2.2.1 from
$\partial _{\tilde z}f(z)=0$ on $U$ and Conditions $2.2.(1,2,5)$
it follows, that for each $j$: \\
$\partial f(z)/\partial {\bar z}_{2j,2j+1}=(\partial f(z)/\partial
{\tilde z}).(\partial {\tilde z}/\partial {\bar z}_{2j,2j+1})+
(\partial f(z)/\partial z).(\partial z/\partial {\bar z}_{2j,2j+1})=0$, \\
since $\partial _{\tilde z}f(z)=0$ and $\partial z/
\partial {\bar z}_{2j,2j+1}=0$.  Generally for $f$ its derivative
$f'(z)$ need not be right
(or left) superlinear on ${\cal A}_r$. Then $\partial f(z)/\partial
{\bar z}_{2j,2j+1}$ are $\bf R$-homogeneous additive operators
on the $\bf R$-linear subspace $({\bf R}\oplus s^*p{\bf R})$ of
${\cal A}_r$.
\par Let $f$ satisfy $2.2.(1,2,5)$ and $(1)$, then
\par $(Df(z)).h=\sum_{s\in \bf b}(Df(z)).h_ss=\sum_{s\in \bf b}
(\partial f(z)/\partial w_s)h_s$, \\
since $\partial z/\partial w_s=s$ for each $s\in \bf b$,
where $h=\sum_{s\in \bf b}h_ss\in {\cal A}_r$, $h_s\in \bf R$
$\forall $ $s\in \bf b$. From \\
$(\partial f(z)/\partial w_s)h_s+ (\partial f(z)/\partial w_p)h_p
=(\partial f(z)/\partial z_{2j,2j+1}).(h_s +s^*ph_p)$ \\
for each $s=i_{2p}$ and $p=i_{2j+1}$, since \\
$(\partial f(z)/\partial w_s)=(\partial f(z)/\partial z_{2j,2j+1}).1
+(\partial f(z)/\partial {\bar z}_{2j,2j+1}).1$, \\
$(\partial f(z)/\partial w_p)=(\partial f(z)/\partial z_{2j,2j+1}
- \partial f(z)/\partial {\bar z}_{2j,2j+1}).(s^*p)$ \\
and $\partial f(z)/\partial {\bar z}_{2j,2j+1}=0$,
it follows that $(Df(z)).h=(\partial f(z)/\partial z).h$,
since by Proposition 2.2.1 and Conditions $(1)$ \\
$(\partial f(z)/\partial {\tilde z}).(h_ss+h_pp)
=(\partial f(z)/\partial {\bar z}_{2j,2j+1}).(h_s+h_ps^*p)=0$
for each $j$ \\
and $(\partial f/\partial {\tilde z}).h=\sum_{j=0}^{2^{r-1}-1}
(\partial f/\partial {\tilde z}).(i_{2j}(h_{2j}+h_{2j+1}i_{2j}^*i_{2j+1}))$,
since $h_s\in \bf R$ for each $s\in \bf b$.
\par {\bf 2.6. Proposition.} {\it  Let $f$ be a function
as in \S 2.5 and suppose that there are two constants $r$ and $R$ such
that the Laurent series $(2.5.i,ii)$ converges in the set
$B(a,r,R,{\cal A}_m):=\{ z\in {\cal A}_m: r\le |z-a|\le R \} $ for each
$s=1,...,j$, let also $\gamma $ be a rectifiable path contained in
$U\cap B(a,r',R',{\bf H})$, where $r<r'<R'<R$, $m\ge 3$.
Then the Cayley-Dickson algebra line integral exists.}
\par {\bf Proof.} At first mention that ${\cal A}_m$ is the normed
algebra such that $|\xi \eta |\le |\xi | |\eta |$ for each
$\xi $ and $\eta \in {\cal A}_m$. This can be proved by induction
starting from $\bf C$ and using the doubling procedure.
Suppose $m\ge 2$ and ${\cal A}_{m-1}$ is the normed algebra,
then $|(a,b)(c,d)|=(|ac-d^*b|^2+|da+bc^*|^2)^{1/2}\le $
$(|ac|^2+|d^*b|^2+|da|^2+|bc^*|^2)^{1/2}\le (|a|^2+|b|^2)^{1/2}
(|c|^2+|d|^2)^{1/2}=|(a,b)| |(c,d)|$, where $a, b, c, d\in
{\cal A}_{m-1}$, $(a,b)$ and $(c,d)\in {\cal A}_m$.
Since each $f^s$ converges in $B(a,r,R,{\cal A}_m)$, then
$${\overline {\lim }}_{n+m>0}|f^s_{n,m}|^{1/(n+m)}R\le 1,\mbox{ hence}$$
$$\| f \|_{\omega }:=\prod_{s=1}^j (\sup_{n+m<0}|f^s_{n,m}| r^{n+m},
\sup_{n+m\ge 0} |f^s_{n,m}| R^{n+m})<\infty $$
 and inevitably
$$\| f \| _{1,\omega ,B(a,r',R',{\cal A}_m)}
:=\prod_{s=1}^j [(\sum_{n+m<0}|f^s_{n,m}|{r'}^{n+m})
+(\sum_{n+m>0}|f^s_{n,m}| {R'}^{n+m})] <\infty .$$
For each locally $z$-analytic
function $f$ in $U$ and each $z_0$ in $U$ there exists a ball of radius
$r>0$ with center $z_0$ such that $f$ has a decomposition
analogous to $(2.5.i,ii)$ in this ball with all $n$ nonnegative.
Consider two $z$-locally analytic
functions $f$ and $q$ on $U$ such that $f$ and $q$ noncommute.
Let $f^0:=f$, $q^0:=q$, $q^{-n}:=q^{(n)}$, $(\partial (q^n)/\partial z).1 =:
q^{n-1}$ and $q^{-k-1}=0$ for some $k\in \bf N$, then
\par $(i)$ $(fq)^1=f^1q-f^2q^{-1}+f^3q^{-2}+...+(-1)^kf^{k+1}q^{-k}$.
In particular, if $f=az^n$, $q=bz^k$, with $n>0$, $k>0$,
$b\in {\cal A}_m\setminus {\bf R}I$, then $f^p=[(n+1)...(n+p)]^{-1}az^{n+p}$
for each $p\in {\bf N}$, $q^{-s}=k(k-1)...(k-s+1)bz^{k-s}$
for each $s\in \bf N$. Also
\par $(ii)$ $(fq)^1=fq^1-f^{-1}q^2+f^{-2}q^3+...+(-1)^nf^{-n}q^{n+1}$,
when $f^{-n-1}=0$ for some $n\in \bf N$.
Apply $(i)$ for $n\ge m$ and $(ii)$ for $n<k$
to solve the equation $(\partial g(z,{\tilde z})/\partial z).1=
f(z,{\tilde z})$ for each $z\in U$. If $f$ and $q$
have series converging in $Int (B(z_0,r,{\cal A}_m))$,
then these formulas show that there exists a $z$-analytic
function $(fq)^1$ with series converging in $Int (B(0,r,{\cal A}_m))$,
since $\lim_{n\to \infty }(nr^n)^{1/n}=r$, where $0<r<\infty $. 
Applying this formula by induction to products of polynomials
$\{ P_1...p_n \} _{q(n)}$ and converging series, we get
$g$. Since $f$ is locally analytic, then $g$ is also locally analytic.
Therefore, for each locally $z$-analytic function
$f$ there exists the operator $\hat f$. Considering a function
$G$ of real variables corresponding to $g$ we get that
due to Lemma 2.5.1 all solutions
$g$ differ on constants in ${\cal A}_m$,
since $\partial g/\partial w_s+(\partial g/\partial w_p).(s^*p)=0$
for each  $s=i_{2j}$, $p=i_{2j+1}$, $j=0,1,...,2^{r-1}-1$
and $\partial g/\partial w_1$ is unique,
hence $\hat f$ is unique for $f$. Therefore,
$$(1)\quad | \{ f^s_{n,m} (z_{j+1}-a)^k(\Delta z_j)(z_{j+1}-a)^{n-k}
({\tilde z}_{j+1}-{\tilde a})^m \} _{q(n+m+1)} | \le $$
$$|f^s_{n,m}| |(z_{j+1}-a)^n
({\tilde z}_{j+1}-{\tilde a})^m| |\Delta z_j|.$$
From Equation $(1)$ it follows, that 
$|I(f,\gamma ;P)| \le \| f \| _{1,\omega ,B(a,r',R',{\bf H})} v(\gamma ;P),$
for each $P$, and inevitably
$$(2)\quad |I(f,\gamma ;P)-I(f,\gamma ;Q)|\le w(\hat f;P)V(\gamma )$$
for each $Q\supset P$, where
$$(3)\quad w(\hat f;P):=\max_{(z, \zeta \in \gamma ([c_j,c_{j+1}]))}
\{ \| {\hat f}(z)-{\hat f}(\zeta ) \|:
\quad z_j=\gamma (c_j), c_j\in P \} ,$$
$ \| {\hat f}(z) -{\hat f}(\zeta ) \|
:= \sup_{h\ne 0}|{\hat f}(z).h -{\hat f}(\zeta ).h |/|h|$.
Since $\lim_{n\to \infty } (n)^{1/n}=1,$ then
$\lim_P \omega ({\hat f},P)=0.$
From $\lim_Pw({\hat f};P)=0$ the existence of $\lim_PI(f,\gamma ;P)$
now follows.
\par {\bf 2.7. Theorem.} {\it Let $\gamma $ be a rectifiable
path in $U$, then the Cayley-Dickson algebra ${\cal A}_r$, $r\ge 3$,
line integral has a continuous extension on the space $C^0_b(U,{\cal A}_r)$
of bounded continuous functions $f: U \to {\cal A}_r$.
This integral is an $\bf R$-linear and left-${\cal A}_r$-linear
and right-${\cal A}_r$-linear functional on $C^0_b(U,{\cal A}_r)$.}
\par {\bf Proof.} Since $\gamma $ is continuous and $[a,b]$ is compact,
then there exists a compact canonical closed subset $V$ in ${\cal A}_r$,
that is, $cl (Int(V))=V$, such that $\gamma ([a,b])\subset V\subset U$.
Let $f\in C^0_b(U,{\cal A}_r)$,
then in view of the Stone-Weierstrass theorem
for a function $F(w_s: s\in {\bf b})=f\circ \sigma (w_s: s\in {\bf b})$
and each $\delta >0$ there exists a polynomial $T$ such that
$\| F- T\|_0<\delta $,
where $\| f\|_0:=\sup_{z\in U}|f(z)|$. This polynomial takes values in
${\cal A}_r$, hence it has the form: $T=\sum_{s\in \bf b}T_ss$,
where $T_s: U\to \bf R$. There are relations $zp=\sum_{s\in \bf b}w_ssp$,
$sp=-ps$ for each $s\ne p\in {\hat b}$, consequently,
$zp=-w_p+\sum_{s\in {\bf b}; s\ne p}w_ssp$,
$(zp)^*=p^*z^*= -w_p + \sum_{s\in {\bf b}; s\ne p}w_s(sp)^*$
$=-w_p - \sum_{s\in {\bf b}, s\ne p}w_ssp$, since $p^*=-p$
for each $p\in {\hat b}$. Then $w_1=(z+{\tilde z})/2$ and
$w_p=(p{\tilde z}-zp)/2$ for each $p\in \hat b$, where
we use the identity ${\tilde z}=(2^r-2)^{-1} \{ -z +\sum_{s\in \hat b}
s(zs^*) \} $ (see \S 2.1).
Write $F=\sum_{s\in \bf b}F_ss$, where $F_s\in \bf R$ for each $s\in \bf b$,
then the application of the Stone-Weierstrass theorem by real variables
$(w_s: s\in {\bf b})$ expressed through $z$ and $s\in \bf b$ with real
constant multipliers gives that the $\bf R$-linear space of functions
given by Equations $(2.5.i,ii)$ is dense in $C^0_b(U,{\cal A}_r)$.
\par Consider a function $g(z)$ on $U$. Let $g(z)$ be superdifferentiable
by $z$. Consider a space of all such that $g$ on $U$ for which
$(Dg(z)).s$ is a bounded continuous function on $U$ for each $s\in \bf b$,
it is denoted by $C^1_b(U,{\cal A}_r)$
and it is supplied with the norm $\| g \| _{C^1_b}:=\| g\|_{C^0_b}+
\sum_{s\in \bf b} \| (Dg(z).s \|_{C^0_b}$, where
$\| g \|_{C^0_b}:=\sup_{z\in U}|g(z)|$,
hence $(Dg(z)).h\in C^0_b(U\times B(0,0,1,{\cal A}_r),{\cal A}_r),$
where $h\in B(0,0,1,{\cal A}_r)$. Therefore, there exists
a positive constant $C$ such that
$$(1)\quad \sup_{h\ne 0} | (Dg(z)).h|/|h| \le C
\sum_{s\in \bf b} \| (Dg(z)).s \|_{C^0_b} ,$$
since $h=\sum_{s\in \bf b}h_ss$ for each $h\in {\cal A}_r$
and $Dg(z)$ is $\bf R$-linear and
$(Dg(z)).(\mbox{ }^1h+\mbox{ }^2h)=(Dg(z)).\mbox{ }^1h+(Dg(z)).\mbox{ }^2h$
for each $\mbox{ }^1h$ and $\mbox{ }^2h\in {\cal A}_r$,
where $h_s$ is a real number for each $s\in \bf b$,
$G(w_s: s\in {\bf b}):=g\circ \sigma (w_s: s\in {\bf b})$
is Frech\'et differentiable on an open subset $U_{\sigma }\subset
\bf R^{2^r}$ such that $\sigma (U_{\sigma })=U$.
\par In \S 2.6 it was shown that the equation
$(\partial g(z,{\tilde z})/\partial z).1=f(z,{\tilde z})$
has a solution in a class of locally $z$-analytic
functions on $U$. The subset $C^{\omega }(U,{\cal A}_r)$ is dense in
the uniform space $C^0_b(U,{\cal A}_r)$.
\par If $g=\{ g^1...g^j \} _{q(j)}$ is a product of functions
$g^s\in C^1_b(U,{\cal A}_r)$, then $(Dg(z)).h=
\sum_{v=1}^j\{ g^1(z)...g^{v-1}(z)[(Dg^v(z)).h]g^{v+1}(z)...g^j(z)
\} _{q(j)}$ for each $h\in {\cal A}_r$.
Consider the space ${\hat C}^0_b(U,{\cal A}_r):= \{
(Dg(z)).s: s\in {\bf b} \} $.
It has an embedding $\xi $ into $C^0_b(U,{\cal A}_r)$ and
$\| g \| _{C^1_b}\ge \sum_{s\in \bf b} \| (Dg(z)).s \| _{C^0_b} $.
In view of Inequality $(1)$ the completion of
${\hat C}^0_b(U,{\cal A}_r)$ relative to $\| * \|_{C^0_b(U,{\cal A}_r)}$
coincides with $C^0_b(U,{\cal A}_r)$.
\par Let $\{ f^v: \quad v\in {\bf N} \} $ be a sequence of functions
having decomposition $(2.5.i,ii)$ and converging to $f$ in
$C^0_b(U,{\cal A}_r)$ relative to the metric $\rho (f,q):=\sup_{z\in U}
|f(z)-q(z)|$ such that $f^v=\xi ((Dg^v(z)).s: s\in {\bf b})$
for some $g^v\in C^1_b(U,{\cal A}_r)$. Relative to this metric
$C^0_b(U,{\cal A}_r)$ is complete. We have the equality
$$\partial (\int_0^q F (\phi h_s: s\in {\bf b}))/\partial q =
F(w_s: s\in {\bf b})$$
for each continuous function $F$ on $U_{\sigma }$, where
$w_s=w_{0,s}+qh_s$ for each $s\in \bf b$,
$(w_{0,s}: s\in {\bf b})+\phi (h_s: s\in {\bf b})
\in U_{\sigma }$ for each $\phi \in \bf R$ with $0\le \phi < q+\epsilon $,
$0<\epsilon <\infty $, $h_s\in \bf R$ for each $s\in \bf b$.
Let $z_0$ be a marked point in $V$.
There exists $R>0$ such that $\gamma $ is contained in the
interior of the parallelepiped $V:=
\{ z\in {\cal A}_r: z=\sum_{s\in \bf b}w_ss;$ $|w_s-w_{0,s}|\le R $
$\mbox{for each}$ $s\in {\bf b} \} $.
\par If $V$ is not contained in
$U$ consider a continuous extension of a continuous function $F$
from $V\cap U_0$ on $V$, where $U_0$ is a closed subset in $U$
such that $Int (U_0)\supset \gamma $ (about the theorem
of a continuous extension see \cite{eng}).
Therefore, suppose that $F$ is given on $V$. Then the function
$F_1(w_s: s \in {\bf b}):=\int_{w_{0,1}}^{w_1}...\int_{w_{0,t}}^{w_t}
F(w_s: s\in {\bf b})dw_1...dw_t$ is in $C^1(V,{\cal A}_r)$
(with one sided derivatives on $\partial V$ from inside $V$),
where $t:=i_{2^r-1}$.
Consider a foliation of $V$ by $(2^r-1)$-dimensional
$C^0$-manifolds $\Upsilon _z$ such that $\Upsilon _z\cap \Upsilon _{z_1}=
\emptyset $ for each $z\ne z_1$, where $z, z_1\in \gamma $,
$\bigcup_{z\in \gamma }\Upsilon _z=V_1$, $V_1$ is a canonical closed
subset in ${\cal A}_r$ such that $\gamma \subset V_1\subset V$.
Choose this foliation such that to have decomposition
of a Lebesgue measure $dV$ into the product of measures $d\nu (z) $
along $\gamma $ and $d\Upsilon _z$ for each $z\in \gamma $.
In view of the Fubini theorem there exists $\int_Vf(w_s: s\in {\bf b})dV=
\int_{\gamma }(\int_{\Upsilon _z}f(z)d\Upsilon _z)d\nu (z)$.
If $\gamma $ is a straight line segment then
$\int_{\gamma } f(z)dz$ is in $L^1(\Upsilon ,{\cal A}_r)$.
Let $U_{\bf R}$ be a real region in $\bf R^{2^r}$ corresponding to
$U$ in ${\cal A}_r$.
\par Consider the Sobolev space $W^q_2(U_{\bf R},{\bf R^{2^r}})$
of functions $h: U_{\bf R}\to \bf R^{2^r}$ for which $D^{\alpha }h\in
L^2(U_{\bf R},{\bf R^{2^r}})$ for each $|\alpha |\le q$, where
$0\le q\in \bf Z$. In view of Theorem 18.1.24
\cite{hoermpd} (see also the notation there)
if $A\in \Psi ^m$ is a properly supported
pseudodifferential elliptic operator of order $m$ in the sence that
the principal symbol $a\in S^m(T^*(X))/S^{m-1}(T^*(X))$ has an inverse
in $S^{-m}(T^*(X))/S^{-m-1}(T^*(X))$, then one can find
$B\in \Psi ^{-m}$ properly supported such that $BA-I\in \Psi ^{-\infty }$,
$AB-I\in \Psi ^{-\infty }$. One calls $B$ a parametrix for $A$.
In view of Proposition 18.1.21 \cite{hoermpd} each $A\in \Psi ^m$
can be written as a sum $A=A_1+A_0$, where $A_1\in \Psi ^m$ is properly
supported and the kernel of $A_0$ is in $C^{\infty }$.
In particular we can take a pseudodifferential operator
with the principal symbol $a(x,\xi )=(b+|\xi |^2)^{s/2}$,
where $b>0$ is a constant and $s\in \bf Z$, which corresponds
to $b+\Delta $ for $s=1$ up to minor terms,
where $\Delta =\nabla ^2$ is the Laplacian
(see also Theorem 3.2.13 \cite{grubb} about its parametrix family).
For estimates of a solution there may be also applied Theorem 3.3.2
and Corollary 3.3.3 \cite{grubb} concerning parabolic
pseudodifferential equations for our particular case
corresponding to $(\partial g(z,{\tilde z})/\partial z).1=f$
rewritten in real variables.
\par Due to the Sobolev theorem (see \cite{stein,trieb})
there exists an embedding of the Sobolev space $W^{2^{r-1}+1}_2
(V,{\cal A}_r)$ into $C^0(V,{\cal A}_r)$ such that
\par $(2)$ $\| g \|_{C^0}\le C\| g \|_{W^{2^{r-1}+1}_2}$
for each $g\in W^{2^{r-1}+1}_2$, where $C$ is a positive constant
independent of $g$.
If $h\in W^{k+1}_2(V,{\cal A}_r)$, then $\partial h/\partial w_s
\in W^k_2(V,{\cal A}_r)$ for each $k\in \bf N$ and in particular for
$k=2^{r-1}+1$ and each $s\in \bf b$ (see \cite{stein}).
On the other hand, $\| h \|_{L^2(V,{\cal A}_r)}\le
\| h \|_{C^0(V,{\cal A}_r)}(2R)^{2^{r-1}}$
for each $h\in L^2(V,{\cal A}_r)$. Therefore,
\par $(3)$ $\| A^{-k}h \| _{W^k_2(V,{\cal A}_r)}
\le C \| h \| _{C^0(V,{\cal A}_r)}(2R)^{k+2^{r-1}}$ for each $k\in \bf N$,
where $C=const >0$, $A$ is an elliptic pseudodifferential
operator such that $A^2$ corresponds to $(1+\Delta )$.
For the estimate below the Gronwall Lemma is used
(see, for example, Section 3.3.1 \cite{beldal}), which reads
as follows. Let $\phi (t)$ and $\psi (t)$ be measurable
bounded functions, and $\eta (t)$ be a continuous nonnegative
function such that
\par $\phi (t)\le X+\psi (t)+\int_0^t\eta (\tau ) \phi (\tau )d\tau .$
Then
\par $\phi (t)\le X \exp [\int_0^t\eta (\tau )d\tau ]+
\psi (t) + \int_0^t\exp [\int_{\tau }^t\eta (v)dv]
\psi (\tau )\eta (\tau )d\tau .$  \\
Use this lemma for $\phi (t):= |\int_{z\in
\{ \gamma (v): a\le v\le t \} } (f(z)-q(z)) dz|$,
$X:=C_1 \rho (f,q) V(\gamma )$, $\psi (t)=0$,
$\eta (t)={C'}_2R^{2^r+1}$, where $C_1>0$ and
${C'}_2>0$ are suitable constants independent of $f$, $q$ and $\gamma $,
since $\| ({\hat f}-{\hat q})(z) \| \le \rho (f,q) + \| Im \circ ({\hat f}
-{\hat q})(z) \| $ and $\| Im \circ ({\hat f} -{\hat q})(z) \|
\le \| {\hat f}-{\hat q}(z) \|$ and ${\hat f}(z).1=f$
for each $z\in \gamma ([a,b])$, \\
$|({\hat f}-{\hat q})(\gamma (x_{k+1})).(\gamma (x_{k+1})-\gamma (x_k))|
\le \| ({\hat f}-{\hat q})(\gamma (x_{k+1})) \|
|\gamma (x_{k+1})-\gamma (x_k)|$  \\
for each partitioning $P$: $a=x_0<...<x_k<x_{k+1}<...<x_w=b$,
where $Im (z):=(z-l(zl^*))/2$.
From Equations $2.5.(1,2)$ and Inequalities $(1-3)$ it follows, that
there exists $0<\epsilon <\infty $ such that
$$(4)\quad |I(f-q,\gamma ;P)|\le \rho (f,q)V(\gamma )C_1 \exp (C_2
R^{2^r+2})$$
for each partitioning $P$ of norm $|P|$ less than $\epsilon $,
where $C_1$ and $C_2$ are positive constants independent of $R$, $f$ and $q$.
In view of Formulas $2.6.(1,2)$ $\{ \int_{\gamma }
f^v(z)dz: v\in {\bf N} \} $ is a Cauchy sequence in ${\cal A}_r$
and the latter is complete as the metric space. Therefore,
there exists $\lim_v\lim_PI(f^v,\gamma ;P)=\lim_v\int_{\gamma }f^v(z)dz$,
which we denote by $\int_{\gamma }f(z)dz$.
As in \S 2.6 we get that all solutions $g$ differ on quaternion
constants on each connected component of $U$, consequently,
the functional $\int_{\gamma }$ is uniquely defined on $C^0_b(U,{\cal A}_r)$.
The functional $\int_{\gamma }: C^0_b(U,{\cal A}_r)\to {\cal A}_r$
is continuous due to Formula $(4)$ and evidently
it is $\bf R$-linear, since $\lambda z=z\lambda $
for each $\lambda \in \bf R$ and each $z\in {\cal A}_r$,
that is, $\int_{\gamma }(\lambda _1f_1(z)+
\lambda _2f_2(z))dz=$ $\int_{\gamma }(f_1(z)\lambda _1+
f_2(z)\lambda _2)dz=$ $\lambda _1\int_{\gamma }f_1(z)dz+
\lambda _2\int_{\gamma }f_2(z)dz$
for each $\lambda _1$ and $\lambda _2\in \bf R$,
$f_1$ and $f_2\in C^0_b(U,{\cal A}_r)$.
Moreover, it is left-${\cal A}_r$-linear,
that is, $\int_{\gamma }(\lambda _1f_1(z)+ \lambda _2f_2(z))dz=
\lambda _1\int_{\gamma }f_1(z)dz+ \lambda _2\int_{\gamma }f_2(z)dz$
for each $\lambda _1$ and $\lambda _2\in {\cal A}_r$,
$f_1$ and $f_2\in C^0_b(U,{\cal A}_r)$, since $I(f,\gamma ;P)$ is
left-${\cal A}_r$-linear.
If $g_k\in C^1(U,{\cal A}_r)$, then $g_k\lambda \in C^1(U,{\cal A}_r)$
for each $\lambda \in {\cal A}_r$ and $(Dg_k(z)\lambda ).h=
(Dg_k(z).h)\lambda $ for each $h\in {\cal A}_r$, since
$D\lambda =0$, in particular, for $g_k$
such that $(\partial g_k(z,{\tilde z})/\partial z).1=f_k(z)$
and satisfying $2.5.1.(1)$ by Lemma 2.5.1 on $U$, $k=1, 2$,
since $g(\mbox{ }_1z,\mbox{ }_2z)$ is $\mbox{ }_1z$-superdifferentiable.
Therefore,
$\int_{\gamma }(f_1(z)\lambda _1+f_2(z)\lambda _2)dz=
(\int_{\gamma }f_1(z)dz)\lambda _1+(\int_{\gamma }f_2(z)dz)\lambda _2$
for each constants $\lambda _1$ and $\lambda _2\in {\cal A}_r$,
consequently, $({\hat f}_k(z).h)\lambda =(f_k(z)\lambda )^{\hat .}.h$
for each $h\in {\cal A}_r$. But this certainly does not mean
that $(\int_{\gamma }f(z)dz)\lambda $ and $\lambda (\int_{\gamma }f(z)dz)$
are equal.
\par {\bf 2.8. Remark.} Let $\eta $ be a differential form
on open subset $U$ of the Euclidean space $\bf R^{2^rm}$ with values in
${\cal A}_r$, then it can be written as
$$(1)\quad \eta =\sum_{\Upsilon }\eta _{\Upsilon }db^{\wedge \Upsilon },$$
where $b=(\mbox{ }^1b,...,\mbox{ }^mb)\in \bf R^{2^rm}$,
$\mbox{ }^jb=(\mbox{ }^jb_1,...,\mbox{ }^jb_{2^r})$,
$\mbox{ }^jb_k\in \bf R$, $\eta _{\Upsilon }=\eta _{\Upsilon } (b):
{\bf R^{2^rm}}\to {\cal A}_r$ are $s$ times continuously differentiable
${\cal A}_r$-valued functions with $s\in \bf N$,
$\Upsilon =(\Upsilon (1),...,\Upsilon (m))$, $\Upsilon (j)=
(\Upsilon (j,1),...,\Upsilon (j,2^r))\in
\bf N^{2^r}$ for each $j$, $db^{\wedge \Upsilon }=
d\mbox{ }^1b^{\wedge \Upsilon (1)}\wedge ... \wedge d\mbox{ }^m
b^{\wedge \Upsilon (m)}$, $d\mbox{ }^jb^{\wedge \Upsilon (j)}=
d\mbox{ }^jb_1^{\Upsilon (j,1)}\wedge ... \wedge d\mbox{ }^j
b_{2^r}^{\Upsilon (j,2^r)}$, where $d\mbox{ }^jb_k^0:=1$, 
$d\mbox{ }^jb_k^1=d\mbox{ }^jb_k$, $d\mbox{ }^jb_k^v=0$
for each $v>1$. If $s\ge 1$, then there is defined
an (external) differential
$$d\eta =\sum_{\Upsilon ,(j,k)}(\partial \eta _{\Upsilon }/\partial
\mbox{ }^jb_k)(-1)^{\alpha (j,k)}db^{\wedge (\Upsilon +e(j,k))} ,$$
where $e(j,k)=(0,...,0,1,0,...,0)$ with $1$ on the
$2^r(j-1)+k$-th place, $\alpha (j,k)=(\sum_{p=1}^{j-1}\sum_{v=1}^{2^r}
\Upsilon (p,v))+\sum_{v=1}^{k-1}\Upsilon (j,v)$.
Now use the relations
\par $(2)$ $\mbox{ }^jb_1=(\mbox{ }^jz+
(2^r-2)^{-1} \{ - \mbox{ }^jz + \sum_{s\in \hat b} s(\mbox{ }^jzs^*) \})/2$
and \\
$\mbox{ }^jb_p=(i_p
(2^r-2)^{-1} \{ - \mbox{ }^jz + \sum_{s\in \hat b} s(\mbox{ }^jzs^*) \}
- \mbox{ }^jzi_p)/2$ for each $i_p\in \hat b$.
Then $\eta $ can be expressed in variables $z$. Consider basic elements
$S=(\mbox{ }^1S,...,\mbox{ }^mS)$
and their ordered product $S^{\to \Upsilon }:=
(...(\mbox{ }^1S^{\to \Upsilon (1)}\mbox{ }^2S^{\to \Upsilon (2)})
...)\mbox{ }^mS^{\to \Upsilon (m)}$,
where $\mbox{ }^jS=(\mbox{ }^jS_1,...,\mbox{ }^jS_{2^r})=
(1,i_1,i_2,...,i_{2^r-1})$,
$\mbox{ }^jS^{\to \Upsilon (j)}=(...(i_1^{\Upsilon (j,2)}
i_2^{\Upsilon (j,3)})...)i_{2^p-1}^{\Upsilon (j,2^r)}$, $S^0=1$.
Then Equation $(1)$ can be rewritten in the form:
$$(3)\quad \eta =\sum_{\Upsilon }\xi _{\Upsilon }
d(Sb)^{\wedge \Upsilon },$$
where $Sb=(\mbox{ }^1S_1\mbox{ }^1b_1,...,
\mbox{ }^1S_{2^r}\mbox{ }^1b_{2^r},...,
\mbox{ }^mS_1\mbox{ }^mb_1,...,\mbox{ }^mS_{2^r}
\mbox{ }^mb_{2^r})\in {\cal A}_r^{2^rm}$, \\
$d\mbox{ }^jS_k\mbox{ }^jb_k= \mbox{ }^jS_kd\mbox{ }^jb_k$, \\
$d(Sb)^{\wedge \Upsilon }:=(...((d\mbox{ }^1S\mbox{ }^1
b)^{\wedge \Upsilon (1)} \wedge (d\mbox{ }^2S\mbox{ }^2
b)^{\wedge \Upsilon (2)})\wedge ...) \wedge
(d\mbox{ }^mS\mbox{ }^mb)^{\wedge \Upsilon (m)}$, \\
$(d\mbox{ }^vS\mbox{ }^vb)^{\wedge \Upsilon (v)}:=
(...((d\mbox{ }^vS_1\mbox{ }^vb_1)^{\Upsilon (v,1)}
\wedge (d\mbox{ }^vS_2\mbox{ }^vb_2)^{\Upsilon (v,2)})
\wedge ...)\wedge (d\mbox{ }^vS_{2^r}\mbox{ }^vb_{2^r})^{\Upsilon (v,2^r)}$,
$\xi _{\Upsilon }:=
\eta _{\Upsilon }(S^{\to \Upsilon })^*$.
Relative to the external product $d\mbox{ }^jb_1$
anticommutes with others basic differential $1$-forms
$\mbox{ }^jS_kd\mbox{ }^jb_k$; for $k=2,...,2^r$
these $1$-forms commute with each other relative to the external product.
This means that the algebra of Cayley-Dickson algebra
differential forms is graded relative to the external product.
\par From Equation $(2)$ it follows, that
\par $(4)$ $db_1 = (dz+d (2^r-2)^{-1}
\{ - z + \sum_{s\in \hat b} s(zs^*) \} )/2,$ \\
$db_p = (di_p (2^r-2)^{-1} \{ - z + \sum_{s\in \hat b} s(zs^*) \} 
-dzi_p)/2$ for each $i_p\in \hat b$.
Therefore, the right side of
Equation $(3)$ can be rewritten with $d\mbox{ }^jzi_p$,
$di_p((s\mbox{ }^jz)s^*)$ on the right side, where $i_p\in
\{ 1,i_1,...,i_{2^r-1} \} $, $s\in \hat b$.
Here we can use also $di_p\mbox{ }^j{\tilde z}$ and
$d((2^r-2)^{-1} \{ - \mbox{ }^j{\tilde z} +
\sum_{s\in \hat b} s(\mbox{ }^j{\tilde z}s^*) \}
)i_p$ depending on the considered
representation either $z$ or $\tilde z$ or $(z,{\tilde z})$
of functions and differential forms.
These $1$-forms do neither commute nor anticommute,
since they are not pure elements of the graded algebra. For example,  \\
$d\mbox{ }^jz\wedge d\mbox{ }^jz=
2 \sum_{1 \le v<k\le 2^r-1} i_vi_k d\mbox{ }^jb_{v+1}\wedge d
\mbox{ }^jb_{k+1}$; \\
$(d\mbox{ }^jz)^{\wedge p}=p!\sum_{1\le v(1)<...<v(p)\le 2^r-1}
(...(i_{v(1)}i_{v(2)})...)i_{v(p)} d\mbox{ }^jb_{v(1)+1}\wedge ...
\wedge d\mbox{ }^jb_{v(p)+1}$ for each $3\le p\le 2^r$.
On the other hand Equation $(1)$ can be rewritten
using the identities $(2)$.
This shows, that the exterior differentiation operator
$\mbox{ }_{{\cal A}_r}d$ for ${\cal A}_r$-valued differential forms
over ${\cal A}_r$ and that of
for their real realization $\mbox{ }_{\bf R}d$
coincide and their common operator is denoted by $d$.
Consider the equality
$$(\partial \eta _{\Upsilon }/\partial \mbox{ }^jb^l)
\mbox{ }^jb^l \wedge db^{\Upsilon }=
[(\partial \eta _{\Upsilon }/\partial \mbox{ }^jz).
(\partial \mbox{ }^jz/\partial \mbox{ }^jb^l)]
\mbox{ }^jb^l \wedge db^{\Upsilon }$$ $$+
[(\partial \eta _{\Upsilon }/\partial \mbox{ }^j{\tilde z}).
(\partial \mbox{ }^j{\tilde z}/\partial \mbox{ }^jb_l)]
\mbox{ }^jb^l \wedge db^{\Upsilon }.$$
Applying it to $l=1,...,2^r$ and summing left anf right
parts of these equalities we get
$d\eta (z,\tilde z)=((\partial \eta /\partial z).
d\mbox{ }^jz)\wedge db^{\Upsilon }+ ((\partial \eta /\partial {\tilde z}).
d\mbox{ }^j{\tilde z})\wedge db^{\Upsilon }$,
hence the external differentiation can be presented in the form
$$(5)\quad d=\partial _z+\partial _{\tilde z},$$ where $\partial _z$
and $\partial _{\tilde z}$ are external differentiations
by variables $z$ and $\tilde z$ respectively.

\par Certainly for an external product $\eta _1\wedge \eta _2$
there is not (in general) a $\lambda \in {\cal A}_r$
such that $\lambda \eta _2\wedge \eta _1=
\eta _1\wedge \eta _2$, if $\eta _1$ and $\eta _2$ are not
pure elements (even or odd) of the graded algebra
of differential forms over ${\cal A}_r$.
\par {\bf 2.9. Definition.} A Hausdorff topological space
$X$ is said to be $n$-connected for $n\ge 0$
if each continuous map $f: S^k\to X$ from the $k$-dimensional
real unit sphere into $X$ has a continuous extension
over $\bf R^{k+1}$ for each $k\le n$. A $1$-connected space
is also said to be simply connected.
\par {\bf 2.10. Remark.} In accordance with Theorem
1.6.7 \cite{span} a space $X$ is $n$-connected
if and only if it is path connected and $\pi _k(X,x)$ is trivial
for every base point $x\in X$ and each $k$ such that $1\le k\le n$.
\par Denote by $Int (U)$ an interior of a subset $U$
in a topological space $X$, by $cl (U)=\bar U$ a closure
of $U$ in $X$. For a subset $U$ in ${\cal A}_r$, let
$\pi _{s,p,t}(U):= \{ u: z\in U, z=\sum_{v\in \bf b}w_vv,$
$u=w_ss+w_pp \} $ for each $s\ne p\in \bf b$,
where $t:=\sum_{v\in {\bf b}\setminus \{ s, p \} } w_vv
\in {\cal A}_{r,s,p}:= \{ z\in {\cal A}_r:$ $z=\sum_{v\in \bf b}
w_vv,$ $w_s=w_p=0 ,$ $w_v\in \bf R$ $\forall v\in {\bf b} \} $.
That is, geometrically $\pi _{s,p,t}(U)$ is the projections on
the complex plane ${\bf C}_{s,p}$ of the intersection
of $U$ with the plane ${\tilde \pi }_{s,p,t}\ni t$,
${\bf C}_{s,p}:=\{ as+bp:$ $a, b \in {\bf R} \} $,
since $sp^*\in {\hat b}$.
\par {\bf 2.11. Theorem.} {\it Let $U$ be a domain in ${\cal A}_r$,
$r\ge 3$, such that $\emptyset \ne Int (U)\subset U\subset cl (Int (U))$
and $U$ is $(2^r-1)$-connected; $\pi _{s,p,t}(U)$
is simply connected in $\bf C$ for each $k=0,1,...,2^{r-1}-1$,
$s:=i_{2k}$, $p:=i_{2k+1}$,
$t\in {\cal A}_{r,s,p}$ and $u\in {\bf C}_{s,p}$
for which there exists $z=u+t\in U$ (see \S 2.10).
Suppose $f\in C^0_b(U,{\cal A}_r)$ and $f$ is superdifferentiable
by $z\in U$ and $f$ has a continuous extension on an open domain $W$
such that $W\supset U$. Then for each
rectifiable closed path $\gamma $ in $U$ a Cayley-Dickson
line integral $\int_{\gamma }f(z)dz=0$ is equal to zero.}
\par {\bf Proof.} In view of Proposition 2.3 $f$ is $z$-represented
and $\partial _{\tilde z}f=0$ on $U$. Therefore, $\xi (\mbox{ }_1z,
\mbox{ }_2z)$ is independent from $\mbox{ }_2z$, where $\xi $
is the corresponding to $f$ function from \S 2.5, consequently,
$g(\mbox{ }_1z,\mbox{ }_2z)$ is also independent from $\mbox{ }_2z$
and we can write $g(z)$ shortly.
For a path $\gamma $ there exists
a compact canonical closed subset in ${\cal A}_r$:
$\quad W\subset Int (U)$
such that $\gamma ([0,1])\subset W$, since $\gamma $
is rectifiable and ${\cal A}_r$ is locally compact.
In view of Theorem 2.7 for each sequence of functions
$f_n\in C^1(U,{\cal A}_r)$ converging to $f$ in $C^0_b(U,{\cal A}_r)$
such that $f_n(z)=(\partial g_n(z)/\partial z).1$
with $g_n(z)\in C^2(U,{\cal A}_r)$ and satisfying conditions
of \S 2.5, since $\xi $ is independent from
$\mbox{ }_2z$, and each sequence of paths $\gamma _n: [0,1] \to U$
$C^3$-continuously differentiable and converging to $\gamma $
relative to the total variation $V(\gamma -\gamma _n)$
there exists $\lim_n\int_{\gamma _n}f_n(z,{\tilde z})dz=
\int_{\gamma }f(z)dz$.
Therefore, it is sufficient to consider the case of $f\in C^1(U,
{\cal A}_r)$ such that $f(z)=(\partial g(z)/\partial z).1$ on $U$,
and continuously differentiable $\gamma $,
where $g\in C^2(U,{\cal A}_r)$ satisfies conditions of \S 2.5.
Denote the integral $\int_{\gamma }f(z)dz$ by $Q$.
We can write this integral in the form $Q=\int_0^1 (\partial
g(z)/\partial z).\gamma '(t)dt$.
Write $f$ in the form:
\par $f=\sum_{s\in \bf b}f_ss=\sum_{\beta =0}^{2^{r-1}-1}
f_{2\beta ,2\beta +1}$, where $f_{\beta ,\nu }:=f_{i_{\beta }}i_{\beta }+
f_{i_{\nu }}i_{\nu }$, $f_s\in \bf R$ for each $s\in \bf b$. Therefore,
\par $d\gamma (t)=\gamma '(t)dt=\sum_{j=0}^{2^{r-1}-1}
{\gamma '}_{2j,2j+1}(t)dt$.
The condition $\gamma (0)=\gamma (1)$ is equivalent to
$\gamma _{2j,2j+1}(0)=\gamma _{2j,2j+1}(1)$ for each $j=0,1,...,2^{r-1}-1$.
We have $\gamma _{\beta ,\nu }\subset \pi _{\beta ,\nu ,t}(U)$
for each $\beta \ne \nu \in \bf b$. The multiplication in ${\cal A}_r$
is distrbutive, consequently,
\par $(\partial g(z)/\partial z).
(d\gamma (t))=\sum_{k=0}^{2^{r-1}-1}(\partial g(z)/\partial
z).(d\gamma _{2k,2k+1}(t))$. \\
In view of the Hurewicz isomorphism theorem (see \S 7.5.4
\cite{span}) $H_q(U,x)=0$ for each $x\in U$ and each
$q<2^r$, hence $H^l(U,x)=0$ for each $l\ge 1$.
\par If $f: Y\to V$ is continuous, then $r\circ f: Y\to
\Omega $ is continuous, if $f$ is onto $V$, then
$r\circ f$ is onto $\Omega $, where $r: V\to \Omega $ is a retraction,
$V$, $Y$ and $\Omega $ are topological spaces.
The topological space $U$ is metrizable, hence for each
closed subset $\Omega $ in $U$ there exists a canonical closed subset
$V\subset U$ such that $V\supset \Omega $ and $\Omega $
is a retraction of $V$, that is, there exists a continuous mapping
$r: V\to \Omega $, $r(z)=z$ for each $z\in \Omega $
(see \cite{eng} and Theorem 7.1 \cite{isbell}).
Therefore, if $V$ is a $(2^r-1)$-connected canonical closed subset
of $U$ and $\Omega $ is a two dimensional $C^0$-manifold
such that $\Omega $ is a retraction of $V$, then $\Omega $
is simply connected, since each continuous mapping
$f: S^k\to \Omega $ with $k\le 1$ has a continuous extension
$f: {\bf R^{k+1}}\to V$ and $r\circ f: {\bf R^{k+1}}\to \Omega $
is also a continuous extension of $f$ from $S^k$ on $\bf R^{k+1}$.
\par From $(2^r-1)$-connectedness of $U$ it follows, that there is
a two dimensional real differentiable manifold $\Omega $
contained in $U$ such that $\partial \Omega =\gamma $.
This may be lightly seen by considering partitions ${\sf Z}_n$ of $U$
by $S^n_{l,k}\cap U$ and taking $n\to \infty $, where
$S^n_{l,k}$ are parallelepipeds in ${\cal A}_r$ with ribs of length
$n^{-1}$, $l$, $k$ and $n\in {\bf N}$, two dimensional faces
$\mbox{ }_1S^n_l$ and $2^{r-1}$-dimensional faces
$\mbox{ }_2S^n_k$ of $S^n_{l,k}=
\mbox{ }_1S^n_l\times \mbox{ }_2S^n_k$ are parallel to
${\bf C}_{s,p}$ or ${\cal A}_{r,s,p}$ with $s=i_{2k}$
and $p=i_{2k+1}$ respectively such that there exists a sequence
of paths $\gamma _n$ converging to $\gamma $ relative to $|*|_{{\cal A}_r}$
and a sequence of (continuous) two dimensional $C^0$-manifolds
$\Omega ^n$ with $\partial \Omega ^n=\gamma ^n$,
$\Omega ^n\subset \bigcup_{l,k}[(\partial \mbox{ }_1S^n_l)\times (\partial
\mbox{ }_2S^n_k)]$.
Choose $\Omega $ orientable and of class $C^3$ as Riemann
manifolds such that taking their projections on ${\bf C}_{s,p}$ 
the corresponding paths $\gamma _{2k,2k+1}$ and regions $\Omega _{s,p}$ 
in ${\bf C}_{s,p}$ satisfy the conditions mentioned above in this proof.
\par To the appearing integrals the classical (generalized)
Stokes theorem can be applied (see Theorem V.1.1 \cite{weint}):
\par $\int_{\Omega _{2k,2k+1}}\eta (v)=\int_0^1(
\partial g(z)/\partial z).{\gamma '}_{2k,2k+1}(t)dt$, \\
where $\Omega _{2k,2k+1}$ is the simply connected domain in
${\bf C}_{i_{2k},i_{2k+1}}$
such that $\partial \Omega _{2k,2k+1}=\gamma _{2k,2k+1}$ for each $k$,
$\eta (v)=d[(\partial g(z)/\partial z).dv]$,
$v=z_{2k,2k+1}\in \Omega _{2k,2k+1}\subset {\bf C}_{i_{2k},i_{2k+1}}$.
The function $g$ is in $C^2(U,{\cal A}_r)$, hence
$(D^2g(z)).(h_1,h_2):=$ $(D^2g(z)).(h_2,h_1)$
for each $h_1$ and $h_2$ in ${\cal A}_r$. Therefore, due to conditions
of \S 2.5 imposed on $g$ we have
\par $\int_{\Omega _{2k,2k+1}}\eta (v)=
\int_{\Omega _{2k,2k+1}}d[(\partial g(z)/\partial z).dv]=
\int_{\Omega _{2k,2k+1}}d^2q(z_{2k,2k+1})=0$  \\
for each $k=0,1,...,2^{r-1}-1$, since \\
$(\partial g(z)/\partial z).{\gamma '}_{2k,2k+1}=(\partial g(z)/
\partial z).(s({\gamma '}_{2k}+ {\gamma '}_{2k+1}s^*p))$, such that \\
$(\partial g(z)/\partial z).{\gamma '}_{2k,2k+1}=
(\partial g(z)/\partial w_s).{\gamma '}_{2k} +
(\partial g(z)/\partial w_p).{\gamma '}_{2k+1}$ \\
$=(\partial g(z)/ \partial z_{2k,2k+1}).({\gamma '}_{2k}
+ {\gamma '}_{2k+1}s^*p)$, \\
$\partial y(z)/\partial w_s=(\partial y(z)/\partial z_{2k,2k+1}+
\partial y(z)/\partial {\bar z}_{2k,2k+1}).1$ and \\
$\partial y(z)/\partial w_p=(\partial y(z)/\partial z_{2k,2k+1}-
\partial y(z)/\partial {\bar z}_{2k,2k+1}).(s^*p)$
for each $z$-superdifferentiable function $y$ on $U$ and
$\partial g(z)/\partial {\bar z}_{2k,2k+1}=0$,
where $q$ corresponds to $g|_{\Omega _{2k,2k+1}}$, $s=i_{2k}$, $p=i_{2k+1}$.
\par {\bf 2.12. Definitions.} A continuous function
on an open domain $U$ in ${\cal A}_r$ such that $\emptyset \ne U$
and $\int_{\gamma }fdz=0$
for each rectifiable closed path $\gamma $ in $U$, then
$f$ is called ${\cal A}_r$-integral holomorphic on $U$ (see \S 2.5).
\par If $f$ is a $z$-superdifferentiable function
on $U$, then it is called ${\cal A}_r$-holomorphic on $U$.
\par {\bf 2.13. Corollary.} {\it Let $f$ be ${\cal A}_r$-holomorphic
function on an open $(2^r-1)$-connected domain $U$ in ${\cal A}_r$
such that $\pi _{s,p,t}(U)$ is simply connected
in ${\bf C}_{s,p}$ for each $t\in {\cal A}_{r,s,p}$ and $u\in
{\bf C}_{s,p}$, $s:=i_{2k}$, $p:=i_{2k+1}$ for which there exists
$z=t+u\in U$, then $f$ is ${\cal A}_r$-integral holomorphic.}
\par This follows immediately from Theorem 2.11.
\par {\bf 2.14. Definition.} Let $U$ be a subset of ${\cal A}_r$
and $\gamma _0: [0,1]\to {\cal A}_r$ and $\gamma _1: [0,1]\to {\cal A}_r$
be two continuous paths. Then $\gamma _0$ and $\gamma _1$ are called
homotopic relative to $U$, if there exists a continuous mapping
$\gamma : [0,1]^2\to U$ such that $\gamma ([0,1],[0,1])\subset U$
and $\gamma (t,0)=\gamma _0(t)$ and $\gamma (t,1)=\gamma _1(t)$ for each
$t\in [0,1]$.
\par {\bf 2.15. Theorem.} {\it Let $W$ be an open subset in ${\cal A}_r$,
$r\ge 3$, and $f$ be an ${\cal A}_r$-holomorphic function on $W$
with values in ${\cal A}_r$. Suppose that there are
two rectifiable paths $\gamma _0$ and $\gamma _1$ in $W$ with common
initial and final points ($\gamma _0(0)=\gamma _1(0)$ and $\gamma _0(1)
=\gamma _1(1)$) homotopic relative to $U$, where $U$ is a
$(2^r-1)$-connected subset in $W$ such that
$\pi _{s,p,t}(U)$ is simply connected in $\bf C$ for each
$t\in {\cal A}_{r,s,p}$ and $u\in {\bf C}_{s,p}$,
$s=i_{2k}$, $p=i_{2k+1}$, $k=0,1,...,2^{r-1}-1$,
for which there exists $z=u+t \in U$. 
Then $\int_{\gamma _0}fdz=\int_{\gamma _1}fdz$.}
\par {\bf Proof.} A homotopy of $\gamma _0$ with $\gamma _1$ realtive to $U$
implies homotopies of $(\gamma _0)_{2j,2j+1}$ with $(\gamma _1)_{2j,2j+1}$
relative to $\pi _{2j,2j+1,t}(U)$ in ${\bf C}_{s,p}$
with $s=i_{2j}$ and $p=i_{2j+1}$ for each
$j=0,1,2^{r-1}-1$ for each $t\in {\cal A}_{r,s,p}$ and $u\in {\bf C}_{s,p}$
for which there exists $z=t+u\in U$.
Consider a path $\zeta $ such that $\zeta (t)=
\gamma _0(2t)$ for each $0\le t\le 1/2$ and $\zeta (t)=\gamma _1(2-2t)$
for each $1/2 \le t\le 1$. Then $\zeta $ is a closed path
contained in a $U$. In view of Theorem 2.11 $\int_{\zeta }f(z)dz=0$.
On the other hand, $\int_{\zeta }f(z)dz=\int_{\gamma _0}f(z)dz-
\int_{\gamma _1}f(z)dz$, consequently, $\int_{\gamma _0}f(z)dz=
\int_{\gamma _1}f(z)dz$.
\par {\bf 2.15.1. Corollary.} {\it Let $f\in C^1$ satisfies conditions
of Theorem 2.15. Then for each $z\in U$ there exists $(\partial
(\int_{\gamma }f(\zeta )d\zeta )/ \partial z).h={\hat f}(z).h$ for each
$h\in {\cal A}_r$, where $\gamma (0)=z_0$, $\gamma (1)=z$,
$z_0$ is a marked point in $U$.}
\par {\bf 2.16. Theorem.} {\it Let $f$ be a ${\cal A}_r$ locally
$z$-analytic function on an open domain $U$ in ${\cal A}_r^n$, then $f$ is
${\cal A}_r$-holomorphic on $U$.}
\par {\bf Proof.} From the definition of the superdifferential
we get $(Dz^n).h=\sum_{k=0}^{n-1}z^khz^{n-k-1}$.
Using the formula of the superdifferential for a product of functions,
from \S 2.7 we obtain, that each $f$ of the form
$(2.5.i,ii)$ is superdifferentiable (by $z$) when $n_0\ge 0$ in
$(2.5.ii)$. Using the norm $\| *\|_{\omega }$-convergence
of series with respect to $z$ for a given $f\in C^{\omega }(U,{\cal A}_r)$
we obtain for each $a\in U$, that there exists its neighbourhood $W$,
where $f$ is ${\cal A}_r$-holomorphic, hence $f$ is
${\cal A}_r$-holomorphic on $U$.
\par {\bf 2.17. Note.} In the next section it is shown that
an octonion-holomorphic function is infinite differentiable;
furthermore, under suitable conditions equivalences between the
properties of octonion holomorphicity, octonion integral holomorphicity
and octonion local $z$-analyticity, will be proved there too.
Integral $(2.5.4)$ may be generalized for a continuous function
$q: U\to {\cal A}_r$ such that $V(q\circ \gamma )<\infty $.
Substituting $\Delta z_k$ on $q(z_{k+1})-q(z_k)=:\Delta q_k$
in Formula $(2.5.5)$ we get
$$(1)\quad \int_{\gamma }f(z,{\tilde z})dq(z):=
\lim_PI(f,q\circ \gamma ;P),\mbox{ where}$$
$$(2)\quad I(f,q\circ \gamma ;P)=\sum_{k=0}^{q-1}{\hat f}
(z_{k+1},{\tilde z}_{k+1}). (\Delta q_k).$$
In paticular, if $\gamma \in C^1$ and $q$ is ${\cal A}_r$-holomorphic
on $U$, also $f(z,{\tilde z})=(\partial g(z,{\tilde z})/
\partial z).1$, where $g\in C^1(U,{\cal A}_r)$, then
$$\int_{\gamma }f(z,{\tilde z})dq(z)=\int_0^1
(\partial g(z,{\tilde z})/\partial z).(
(\partial _zq(z)|_{z=\gamma (s)}).\gamma '(s))ds$$
and $V(\gamma )\le \int_0^1 |\gamma '(s)| ds$.
\par Let $f: U\to {\cal A}_r$ be an ${\cal A}_r$-holomorphic function
on $U$, where $U$ is an open subset of
${\cal A}_r^n$. If there exists an ${\cal A}_r$-holomorphic function
$g: U\to {\cal A}_r$ such that $g'(z).1=f(z)$
for each $z\in U$, then $g$ is called a primitive of $f$.
\par {\bf 2.18. Proposition.} {\it Let $U$ be an open connected subset
of ${\cal A}_m^n$, $m\ge 3$, and $g$ be a primitive of $f$ on $U$,
then a set of all primitives of $f$ is:
$\{ h: h=g+C, C=const \in {\cal A}_m \} $.}
\par {\bf Proof.} In view of Lemma $2.5.1$ for each two primitives
$g_1$ and $g_2$ of $f$ for each $z\in U$ there exists a ball
$B\subset U$, $z\in U$, such that $(g_1-g_2)|_U=const \in {\cal A}_m$
(see \S 2.7). Suppose $h'(z)=0$ for each $z\in U$, then
consider $q(s):=h((1-s)a+sz)$ for each $s\in [0,r]$, where
$a$ is a marked point in $U$ and $B(a,r,{\cal A}_m)$ is a ball
contained in $U$, $r>0$, $z\in B(a,r,{\cal A}_m)$. Then
$q$ is correctly defined and $q(0)=q(1)$. Therefore,
the set $V:=\{ z\in U: h(z)=h(a) \} $ is open in $U$, since
with each point $a$ it contains its neighbourhood.
On the other hand, it is closed due continuity of $h$,
hence $V=U$, since $U$ is connected, consequently, $h=const $ on $U$.
\section{Meromorphic functions and their residues.}
At first we define and describe the exponential and
the logarithmic functions of octonion variables
and then apply them to the investigation of octonionic residues.
Moreover, these studies are accomplished also for variables in
Cayley-Dickson algebras ${\cal A}_r$ for each $r\ge 4$.
\par {\bf 3.1. Note and Definition.} For a variable $z\in {\cal A}_r$
with $r\ge 3$ put
$$(3.1.)\quad \exp (z):=\sum_{n=0}^{\infty }z^n/n! .$$
In view of Note 2.1 $z^n$ and hence $\exp (z)$ are correctly
defined, since real numbers commute with each
element of ${\cal A}_r$, $n!\in {\bf N}\subset \bf R$. If $|z|\le R<\infty $,
then the series $(3.1)$ converges, since
$|exp (z)| \le \sum_{n=0}^{\infty } |z^n/n!|\le \exp (R)<\infty .$
Therefore, $\exp : {\cal A}_r\to {\cal A}_r$. The restriction of $\exp $ on
each of the subsets ${\bf Q}_s := \{ z: z\in {\cal A}_r, z=a+bs,$
$a, b\in {\bf R} \} $ is commutative, where $s\in \hat b$,
${\hat b}_r:={\bf b}_r\setminus \{ 1 \} $, ${\bf b}:={\bf b}_r$,
${\hat b}:={\hat b}_r$, but in general
two elements $z_1$ and $z_2\in {\cal A}_r$ do not commute and
the function $\exp (z_1+z_2)$ on ${\cal A}_r^2$ does not
coincide with $\exp (z_1)\exp (z_2)$.
\par {\bf 3.2. Proposition.} {\it Let $z\in {\cal A}_r$, $r\ge 3$, be written
in the form $z=v+M$, where $v\in \bf R$, $M\in {\cal I}_r$,
${\cal I}_r:= \{ \eta \in {\cal A}_r:$ $Re (\eta )=0 \} $, then
$$(3.2) \quad \exp (z)=\exp (v) \exp (M) ,\mbox{ where} $$
$$(3.3)\quad \exp (M)=(\cos |M|) + [(\sin |M|)/ |M|] M$$
for $M\ne 0$ and $\exp (0)=1$.}
\par {\bf Proof.} Consider $M\in {\cal I}_r$ and write it in the form
$M=a+bl$, where $l=i_{2^{r-1}}$ is the element of the doubling procedure
of ${\cal A}_r$ from ${\cal A}_{r-1}$, $a, b\in {\cal A}_{r-1}$,
$Re (a)=0$. Then \\
$M^2=a^2+a(bl)+a^*(bl)-(bl)(bl)^*=-(|a|^2+|b|^2)$,\\
since $a^*=-a$, $(bl)(bl)^*=(bl)^*(bl)=(lb^*)(bl^{-1})=b^*b=bb^*$,
consequently, \\
$M^{2n}=(-|M|^2)^n$, $M^{2n+1}=(-|M|^2)^nM$ for each $1\le n\in \bf Z$.
Therefore,
$$\exp (M)=1+\sum_{n=1}^{\infty } (-|M|^2)^n / (2n)!+
\sum_{n=0}^{\infty }(-|M|^2)^nM/(2n+1)! $$
$$=(\cos |M|) + [(\sin |M|)/|M|]M$$
for each $M\ne 0$, $\exp (0)=1$. Since $\lim_{0 \ne \phi \to 0}
\sin (\phi )/\phi =1$, where $\phi \in \bf R$,
then the limit taken in Formula
$(3.3)$ while $|M|\ne 0$ tends to $0$ gives the particular case
$\exp (0)=1$. Since $v\in \bf R$ commutes with $M$,
$[v,M]=0$, then $\exp (v+M)=\exp (v)\exp (M)$.
\par {\bf 3.3. Corollary.} {\it If $z\in {\cal A}_r$, $r\ge 3$,
is written in the form $z=\sum_{s\in \bf b}w_ss$ with real $w_s$
for each $s\in \bf b$, then $|\exp (z)|=\exp (v)$.}
\par {\bf Proof.} If $\sum_{s\ne 1}w_s^2=0$ this is evident.
Suppose $\sum_{s\ne 1}w_s^2\ne 0$. In view of Formulas $(3.2, 3.3)$
$$(3.4)\quad exp(z)=exp(v)A\mbox{, where }A=\cos |M|+[(\sin |M|)/|M|]M$$
Since $A\in {\cal A}_r$, then $A^*=\cos |M|-[(\sin |M|)/|M|]M$,
$M^*=-M$ and inevitably $|\exp (z)|=\exp (v)$.
\par {\bf 3.4. Corollary.} {\it The function $\exp (z)$
on the set ${\cal I}_r:=\{ z: z\in {\cal A}_r, Re (z)=0 \} $ is periodic
with $(2^r-1)$ generators of periods $s\in {\hat b}_r$
such that $\exp (z(1+2\pi n/|z|))=\exp (z)$ for each
$0\ne z\in {\cal I}_r$ and each integer number $n$.
If $z\in {\cal A}_r$ is written in the form $z=2\pi sM$, where $M
\in {\cal I}_r$, $|M|=1$,
then $\exp (z)=1$ if and only if $s\in {\bf Z}$.}
\par {\bf Proof.} In view of Formulas $(3.2, 3.3)$
$\exp (sM)=1$ for a given $z=sM\in {\cal I}_r$ with $|M|=1$
if and only if $\cos (s|M|)=1$ and $\sin (s|M|)=0$, that is equivalent to
$s\in \{ 2\pi n: n\in {\bf Z} \} $, since $|M|=1$
by the hypothesis of this corollary.
The particular cases of Formula $(3.3)$
are: $w_{s_0}\ne 0$ and $w_s=0$ for each $s\ne s_0$ in
${\hat b}:={\hat b}_r$,
hence $s\in \hat b$ are $(2^r-1)$ generators for the periods of $\exp $.
\par {\bf 3.5. Corollary.} {\it The function $\exp $
is the epimorphism from ${\cal I}_r$ on the $(2^r-1)$-dimensional
unit sphere $S^{2^r-1}(0,1,{\cal A}_r):=\{ z: z\in {\cal A}_r,
|z|=1 \} $. }
\par {\bf Proof.} In view of Corollary 3.3 the image
$\exp ({\cal I}_r)$ is contained in $S^{2^r-1}(0,1,{\cal A}_r)$.
The sphere $S^{2^r-1}(0,1,{\cal A}_r)$ is characterized by the condition
$\sum_{s\in \bf b} w_s^2=1$ or $w_1^2+|M_1|^2=1$, where
$M_1\in {\cal I}_r$. To show that $\exp ({\cal I}_r)=S^{2^r-1}(0,1,
{\cal A}_r)$ it is sufficient to find $z=v+M$, where $v\in \bf R$,
$M\in {\cal I}_r$, such that
$w_1=\cos |M|$, $M_1=[(\sin |M|)/|M|]M$. For this take
$|M|=\arccos w_1\in [0, \pi ]$, since $w_1\in [-1,1]$
and if $|w_1|\ne 1$ put $M=M_1(1-w_1^2)^{-1/2} \arccos w_1$.
In particular, for $w_1=1$ take $M=0$; for $w_1=-1$ take $M=\pi q$,
where $q\in \hat b$.
\par {\bf 3.6. Corollary.} {\it Each element of the Cayley-Dickson
algebra ${\cal A}_r$, $r\ge 3$, has a polar decomposition
$$(3.5)\quad z=\rho \exp (2\pi (\sum_{s\in \hat b}\phi _ss)),$$
where $\phi _s\in [-1,1]$ for each $s\in \hat b$,
$\sum_{s\in \hat b}\phi _s^2=1$, $\rho :=|z|.$}
\par {\bf Proof.} This follows from Formulas $(3.2, 3.3)$ and
Corollary 3.5.
\par {\bf 3.6.1. Definition.} Let ${\cal A}_{\infty }$ denotes
the family consisting of elements $z=\sum_{s\in \bf b}w_ss$
such that ${\tilde z}:=w_1-\sum_{s\in \hat b}w_ss$,
$z{\tilde z}=:|z|^2=\sum_{s\in \bf b}w_s^2<\infty $,
where ${\bf b}:={\bf b}_{\infty }:=\bigcup_{r=2}^{\infty }{\bf b}_r=
\{ 1, i_1, i_2,..., i_{2^r},... \} $,
${\hat b}:={\bf b}\setminus \{ 1 \} $, $w_s\in \bf R$ for each $s$.
\par {\bf 3.6.2. Theorem.} {\it The family ${\cal A}_{\infty }$
has the structure of the normed power-associative
left and right distributive algebra over $\bf R$
with the external involution of order two.}
\par {\bf Proof.}  Let ${\cal I}_{\infty }:= \{ z\in {\cal A}_{\infty } :
Re (z)=0 \} $. Then each $M\in {\cal I}_{\infty }$ is the limit
of the sequence $M_r\in {\cal I}_r$, also $|z|=:\rho $
is the limit of the sequence $\rho _r:=|z_r|$, where $z_r\in {\cal A}_r$.
Therefore, $z=\lim_{r\to \infty }z_r=\lim_{r\to \infty }
\rho _r \{ \cos |M_r| + [(\sin |M_r|)/|M_r|] M_r \} $
$= \rho \{ \cos |M| +[(\sin |M|)/|M|]M \} $
$=\rho \exp (M)$.
Using the polar coordinates $(\rho ,M)$ prove the power-associativity.
There exists the natural projection
$P_r$ from ${\cal A}_{\infty }$ onto ${\cal A}_r$  for each
$r\ge 2$ given by the formulas:
$M_r:=\{ \sum_{s\in {\hat b}_r}m_ss \} |M|
[\sum_{s\in {\hat b}_r}m_s^2]^{-1/2}$ for
$\sum_{s\in {\hat b}_r}m_s^2\ne 0$
and $M_r=0$ in the contrary case for each $M=\sum_{s\in {\hat b}} m_ss
\in {\cal I}_{\infty }$,
where $m_s\in \bf R$ for each $s\in {\hat b}$; then $z_r:=P_r(z)=
\rho _r \{ \cos |M_r| + [ (\sin |M_r|) / |M_r|] M_r \} $,
where $\rho _r=(\sum_{s\in {\bf b}_r} w_s^2)^{1/2}$,
$z=\sum_{s\in \bf b} w_ss$, $\lim_{0\ne \phi \to 0} \sin (\phi )/\phi =1$.
Then $\lim_{r\to \infty }z_r=z$ relative to the norm
$|z|$ in ${\cal A}_{\infty }$.
Therefore, for each $n\in \bf Z$ there exists
$\lim_{r\to \infty }(\rho _r)^n \exp (nM_r)=\rho ^n \exp (nM)=z^n$,
consequently, ${\cal A}_{\infty }$ is power-associative, since
each ${\cal A}_r$ is power-associative, $\cos $ and $\sin $ are continuous
functions.
Evidently, ${\cal A}_{\infty }$
is the $\bf R$-linear space. The continuity of multiplication
relative to the norm $|z|$ follows from the inequalities
$|\xi _r\eta _r -\psi _r\zeta _r|$ $\le
|\xi _r\eta _r-\xi _r\zeta _r| +|\xi _r\zeta _r-\psi _r \zeta _r|$
$\le |\xi | |\eta -\zeta |+ |\xi - \psi | |\zeta |$ and taking
the limit with $r$ tending to the infinity,
since $|\xi _r|\le |\xi |$ and $|\xi _r \eta _r|\le |\xi _r| |\eta _r|$
for each $\xi _r, \eta _r\in {\cal A}_r$ and for each
$r\in \bf N$. Left and right distributivity
$(\xi +\psi )\zeta =\xi \zeta + \psi \zeta $
and $\zeta (\xi +\psi )=\zeta \xi + \zeta \psi $ follow
from taking the limit with $r$ tending to the infinity
and such distributivity in each ${\cal A}_r$.
The involution $z\mapsto {\tilde z}=:z^*$ in ${\cal A}_r$
is of order two, since $(z^*)^*=z$. It is external, since there is not
any finite algebraic relation with constants in ${\cal A}_{\infty }$
transforming the variable $z\in {\cal A}_{\infty }$
into $z^*$. The relation $z^*=\lim_{r\to \infty }
(2^r-2)^{-1} \{ -z_r +\sum_{s\in {\hat b}_r}s(z_rs^*) \} $
$=\lim_{r\to \infty }
(2^r-2)^{-1} \{ -z +\sum_{s\in {\hat b}_r}s(zs^*) \} $
is of infinite order. The relations of the type
$z_r^*=l_{r+1}(z_rl_{r+1}^*)$ in ${\cal A}_r$ use external automorphism
with $l_{r+1}:=i_{2^r}\in {\cal A}_{r+1}\setminus {\cal A}_r$,
moreover, the latter relation is untrue for $z^*$ and
$z\in {\cal A}_{\infty }$ instead of $z_r\in {\cal A}_r$.
\par No any finite set of non-zero constants $a_1,...,a_n\in
{\cal A}_{\infty }$
can provide the automorphism $z\mapsto \tilde z$ of ${\cal A}_{\infty }$.
To prove this consider an $\bf R$-subalgebra $\Upsilon _{M_1,...,M_n}$
of ${\cal A}_{\infty }$ generated by $ \{ M_1,...,M_n \} $, where
$a_j=|a_j|e^{M_j}$, $M_j\in {\cal I}_{\infty }$.
Since $a_1a_1^*=|a_1|^2>0$, then ${\bf R}|a_1|^2={\bf R}\subset
\Upsilon _{M_1,...,M_n}$, hence $1\in \Upsilon _{M_1,...,M_n}$.
If $\Upsilon _{M_1,...,M_n}=\bf R$, then it certainly can not provide
the automorphism $z\mapsto \tilde z$ of ${\cal A}_{\infty }$.
Consider $\Upsilon _{M_1,...,M_n}\ne \bf R$, without loss of generality
suppose $a_1\notin \bf R$. There is the scalar product $Re (z{\tilde y})$
in ${\cal A}_{\infty }$ for each $z, y\in {\cal A}_{\infty }$.
Let $b_1$ be the projection of $a_1$ in a subspace of ${\cal A}_{\infty }$
orthogonal to ${\bf R}1$, then by our supposition $b_1\ne 0$ and $b_1
\in {\cal I}_{\infty }$. Therefore, $b_1^2/|b_1|^2=-1$, consequently,
$\Upsilon _{M_1}$ is isomorphic to $\bf C$. Certainly, no any
${\cal A}_r$, $r\in \bf N$, can provide the automorphism $z\mapsto \tilde z$
of ${\cal A}_{\infty }$. Therefore, without loss of generality suppose,
that $\Upsilon _{M_1,...,M_n}$ is not isomorphic to $\bf C$
and $a_2\notin \bf C$. If $M, N\in {\cal I}_r$ and $Re(MN^*)=0$, then
$MN\in {\cal I}_r$ and hence $(MN)^*=NM=-MN$, since
$A^*=-A$ for each $A\in {\cal I}_r$.
Let $b_2$ be the projection of $M_2$ in a susbspace of ${\cal A}_{\infty }$
orthogonal to $\Upsilon _{M_1}$ relative to the scalar product
$Re (z{\tilde y})$.
Then $b_2\ne 0$ by our supposition and $b_2\in {\cal I}_{\infty }$,
$b_2^2/|b_2|^2=-1$, hence after the doubling procedure with $b_2/|b_2|$
we get, that $\Upsilon _{M_1,M_2}$ is a subalgebra of ${\cal A}_4$.
Then proceed by induction, suppose $\Upsilon _{M_1,...,M_k}$
is a subalgebra of ${\cal A}_{2^k}$, $k\in \bf N$, $k<n$.
Since ${\cal A}_{2^k}$ can not provide the automorphism $z\mapsto \tilde z$
of ${\cal A}_{\infty }$, suppose without loss of generality, that
$a_{k+1}\notin \Upsilon _{M_1,...,M_k}$ and
consider the orthogonal projection $b_{k+1}$ of $M_{k+1}$ in a subspace of
${\cal A}_{\infty }$ orthogonal to $\Upsilon _{M_1,...,M_k}$ relative
to the scalar product $Re (z{\tilde y})$. Then $b_{k+1}\ne 0$ and
$b_{k+1}\in {\cal I}_{\infty }$, $b_{k+1}^2/|b_{k+1}|^2=-1$. Then
the doubling procedure with $b_{k+1}/|b_{k+1}|$ gives the algebra
$\Upsilon _{M_1,...,M_{k+1}}$ which is the subalgebra of
${\cal A}_{2^{k+1}}$, etc.
As the result $\Upsilon _{M_1,...,M_n}$ is the subalgebra
of ${\cal A}_{2^n}$ and can not provide the automorphism
$z\mapsto \tilde z$ of ${\cal A}_{\infty }$, where $a_1,...,
a_n\in {\Upsilon }_{M_1,...,M_n}$ due to the formula of polar
decomposition (see (3.2,3)) of Cayley-Dickson numbers.
\par {\bf 3.6.3. Note.} From Theorem $3.6.2$ it follows, that
${\cal A}_{\infty }$ together with ${\bf C}={\cal A}_1$
are two extreme cases, where the conjugation $z\to z^*$ is the external
automorphism, though the complex field $\bf C$ is easier to work
due to its commutativity and associativity, than ${\cal A}_{\infty }$
which is neither commutative nor associative.
In view of Definition $3.6.1$ and
Theorem $3.6.2$ preceding results can be transferred from the
case ${\cal A}_r$ with $r\ge 4$ on ${\cal A}_{\infty }$.
Definitions $2.1-2.2$ are transferrable on ${\cal A}_{\infty }$,
moreover, in view of algebraic independence
of $z$ and $z^*$ in ${\cal A}_{\infty }$ we have that
$z$ and $z^*$ are automatically independent variables and
having projections $P_r: {\cal A}_{\infty }\to {\cal A}_r$
for each $r\ge 2$ justify our definitions of $C^{\omega }_{z,{\tilde z}}
(U,{\cal A}_r)$ and $C^n_{z,{\tilde z}}(U,{\cal A}_r)$,
also it justifies axioms of superdifferentiations in \S 2.2.
Evidently Propositions $2.2.1, 2.3, 2.6$ and Corollary $2.4$, Lemma $2.5.1$
are true in the case ${\cal A}_{\infty }$ with ${\bf b}={\bf b}_{\infty }$
instead of ${\bf b}={\bf b}_r$. Definition $2.5$ is valuable for
${\cal A}_{\infty }$ also. Theorem $2.7$ is true also for
${\cal A}_{\infty }$, since with the help of projections $P_r$
we have $\gamma =\lim_{r\to \infty }P_r(\gamma )$,
$P_r(\gamma )\subset U_r$, $ \{ P_r(\gamma ): r \in {\bf N} \} $ converges
to $\gamma $ uniformly on a compact segment $[a,b]\subset \bf R$,
$\gamma : [a,b]\to {\cal A}_{\infty }$, hence the line
integral has the continuous unique extension on
$C^0_b(U,{\cal A}_{\infty })$, where $U_r=P_r(U)$.
In Remark $2.8$ use $l_2({\bf R})^m$ instead of $\bf R^{2^rm}$
and representing differential forms $\eta $ over ${\cal A}_{\infty }$ as
pointwise limits (or uniform convergence on compact subsets)
of differential forms over ${\cal A}_r$ while $r$ tends to the infinity,
since $z_r\to z$ while $r$ tends to infinity, where
$z\in {\cal A}_{\infty }$, $z_r:=P_r(z)$. In general: \\
$(i)\quad \eta (z,{\tilde z})=\sum_{I,J} \eta _{I,J}
\{ (d\mbox{ }^{p_1}z^{\wedge I_1}\alpha _1 \wedge
... \wedge d\mbox{ }^{p_n}z^{\wedge I_n}\alpha _n \wedge
d\mbox{ }^{t_1}{\tilde z}^{\wedge J_1}\beta _1\wedge ... \wedge
d\mbox{ }^{t_n}{\tilde z}^{\wedge J_n}\beta _n \} _{q(|I|+|J|+2n)}$ \\
is the differential form over ${\cal A}_{\infty }$,
where each $\eta _{I,J}(z,{\tilde z})$ is a continuous function
on open $U_n$ in ${\cal A}_{\infty }^n$ with values in ${\cal A}_{\infty }$,
$I=(I_1,...,I_n)$,  $J=(J_1,...,J_n)$, $|I|:=I_1+...+I_n$,
$1\le p_1\le p_2\le ... \le p_n\in \bf N$,
$1\le t_1\le t_2\le ... \le t_n\in \bf N$,
$0\le I_k \in \bf Z$, $0\le J_k\in \bf Z$, 
$\alpha _k, \beta _k\in {\cal A}_{\infty }$ are constants
for each $k=1,...,n$, $d\mbox{ }^pz^0:=1$, $d\mbox{ }^p{\tilde z}^0:=1$,
$n\in \bf N$, $\pi ^l_n(U_l)\subset U_n$ for each $l\ge n$,
where $\pi ^l_n: {\cal A}_{\infty }^l\to {\cal A}_{\infty }^n$
is the natural projection for each $l\ge n$.
The convergence on the right of
Formula $(i)$ in the case of infinite series by $I$ or $J$
is supposed relative to the $C^0_b(W,{\cal A}_{\infty }^{\wedge *})$
topology of uniform convergence on $W$, where $W=pr-\lim \{
U_n, \pi ^l_n, {\bf N} \} $, ${\cal A}_{\infty }^{\wedge *}$ is
supplied with the norm topology inherited from the adjoint space of all
poly $\bf R$-homogeneous ${\cal A}_{\infty }$-additive functionals.
\par In Remark $2.10$ define ${\cal A}_{\infty ,s,p}$ and use projections
$\pi _{s,p,t}$ for each $s\ne p\in \bf b$. Theorems $2.11, 2.15$ and
Corollaries $2.13, 2.15.1$ are transferrable
onto ${\cal A}_{\infty }$ by imposing condition of $(2^r-1)$-connectedness
of $P_r(U)=:U_r$ for each $r\ge 3$, considering $\pi _{s,p,t}(U)$
for each $s=i_{2k}$, $p=i_{2k+1}$, $0\le k\in \bf Z$.
Then Definitions $2.12, 2.14$ and Theorem $2.16$, Notes $2.17, 3.1$
are valid for ${\cal A}_{\infty }$ also. The validity of Corollary
$3.3$ follows from the proof of Theorem $3.6.2$.
We also have instead of $3.4$ and $3.5$ the following.
\par {\bf 3.4'. Corollary.} {\it The function $\exp (z)$
on the set ${\cal I}_{\infty }:= \{ z\in {\cal A}_{\infty }:
Re (z)=0 \} $ is periodic
with the infinite family of generators of periods $s\in {\hat b}_{\infty }$
such that $\exp (z(1+2\pi n/|z|))=\exp (z)$ for each
$0\ne z\in {\cal I}_{\infty }$ and each integer number $n$.
If $z\in {\cal A}_{\infty }$ is written in the form $z=2\pi sM$, where $M
\in {\cal I}_{\infty }$, $|M|=1$,
then $\exp (z)=1$ if and only if $s\in {\bf Z}$.}
\par {\bf 3.5'. Corollary.} {\it The function $\exp $
is the epimorphism from ${\cal I}_{\infty }$ on the infinite dimensional
unit sphere $S^{\infty }(0,1,{\cal A}_{\infty }):=
\{ z: z\in {\cal A}_{\infty }, |z|=1 \} $. }
\par {\bf 3.7. Note.} In the noncommutative ${\cal A}_r$,
$2\le r\le \infty $, case
there is the following relation for $\exp $ and its (right) derivative:
$$(3.6) \quad \exp (z)'.h=\sum_{n=1}^{\infty }\sum_{k=0}^{n-1}
((z^k)h)z^{n-k-1}/n!,$$
where $z$ and $h\in {\cal A}_r$. In particular,
$$(3.7)\quad \exp (z)'.v=v \exp (z)$$ for each $v\in \bf R$,
but generally not for all $h\in {\cal A}_r$.
The function $\exp $ is periodic on ${\cal A}_r$, hence
the inverse function denoted by $Ln$ is defined only locally.
\par Let at first $2\le r\in \bf N$.
Consider the space $\bf R^{2^r}$ of all variables $(w_s: s\in {\bf b})$
for which $\exp $ is periodic on ${\cal A}_r$. The condition
$\sum_{s\in {\hat b}}w_s^2=1$ defines in $\bf R^{2^r}$
the unit sphere $S^{2^r-2}$.
The latter has a central symmetry element $C$ for the transformation
$C(w_s: s\in {\bf b})=(-w_s: s\in {\bf b})$. Consider a subset
$P=P_0\cup \bigcup_{q,j_1,...,j_q}P_{j_1,...,j_q}$ of $S^{2^r-2}$,
where $1\le q\le 2^{r-1}-1$, of all points characterized by the conditions:
\par $P_0:=\{ (w_s: s\in {\hat b})\in S^{2^r-2}: w_s\le 0, \forall s\in
{\hat b} \} ,$
\par $P_{j_1,...,j_q}:=\{ (w_s: s\in {\hat b})\in S^{2^r-2}: w_s\le 0,
\forall s\in {\hat b}\setminus \{ j_1,...,j_q \} ,$
$ w_s\ge 0, \forall s\in \{ j_1,...,j_q \} \} ,$
then $P\cup CP=S^{2^r-2}$ and the intersection $P\cap CP$
is $(2^r-3)$-dimensional over $\bf R$.
This sphere $S^{2^r-2}$ corresponds to
the embedding $\theta _1: (w_s: s\in {\hat b})\hookrightarrow (0,w_s:
s\in {\hat b})\in \bf R^{2^r}$.
Consider the embedding of $\bf R^{2^r}$ into ${\cal A}_r$
given by $\theta _2: (w_s: s\in {\bf b})\hookrightarrow
\sum_{s\in {\bf b}} w_ss \in {\cal A}_r$.
This yields the embedding $\theta :=\theta _2\circ \theta _1$
of $S^{2^r-2}$ into ${\cal A}_r$.
Each unit circle with the center $0$ in ${\cal A}_r$ intersects
the equator $\theta (S^{2^r-2})$ of the unit sphere $S^{2^r-1}(0,1,
{\cal A}_r)$.
Join each point $\sum_{s\in \hat b}w_ss$ on $\theta (S^{2^r-2})$
with the zero point in ${\cal A}_r$ by a line $\{ a\sum_{s\in \hat b}
w_ss: s\in {\bar {\bf R}}_+ \} $,
where ${\bar {\bf R}}_+:=\{ a\in {\bf R}: a\ge 0 \} $.
This line crosses a circle embedded into $S^{2^r-1}(0,1,{\cal A}_r)$,
which is a trace of a circle $\{ \exp (2\pi aM): a\in [0,1] \} $
of radius $1$
in ${\cal A}_r$, where $M=\sum_{s\in \hat b}w_ss$.
Therefore, $\psi (s):= \exp (v + 2\pi aM)$
as a function of $(v,a)$ for fixed $(w_s: s\in {\hat b})\in
S^{2^r-2}$ defines a bijection
of the domain $X\setminus \{ aM: a\in {\bar {\bf R}}_+ \} $
onto its image,
where $X$ is $\bf R^2$ embedded as $(v,a)\hookrightarrow
(v+aM) \in {\cal A}_r$, where $v\in \bf R$. This means,
that $Ln (z)$ is correctly
defined on each subset $X\setminus \{ aM: a\in {\bar {\bf R}}_+ \} $
in ${\cal A}_r$. The union $\bigcup_{(w_s: s\in {\hat b})\in P}
\{ a\sum_{s\in \hat b}w_ss: a\in {\bar {\bf R}}_+ \} $
produces the $(2^r-1)$-dimensional
(over $\bf R$) subset $Q:=Q_0\cup \bigcup_{q,j_1,...,j_q}Q_{j_1,...,j_q}$,
where $Q_0=\theta (S_0)$, $Q_{j_1,...,j_q}:=\theta (S_{j_1,...,j_q})$,
$S_0:={\bar {\bf R}}_+P_0$, $S_{j_1,...,j_q}:={\bar {\bf R}}_+
P_{j_1,...,j_q}$, $1\le q\le 2^{r-1}-1$.
Then, on the domain ${\cal A}_r\setminus Q$, the function $\exp (z)$
defines a bijection with image $\exp ({\cal A}_r\setminus Q)$
and its inverse function $Ln (z)$ is correctly defined
on ${\cal A}_r\setminus \exp (Q)$.
By rotating ${\cal A}_r\setminus Q$ one may produce other
domains on which $Ln$ can be defined as
the univalued function (that is, $Ln(z)$ is one point in ${\cal A}_r$),
but not on the entire ${\cal A}_r$. This means that $Ln(z)$ is a
locally bijective function.
We have elementary identities $\cos (2\pi -\phi )=\cos (\phi )$
and $\sin (2\pi -\phi )=-\sin (\phi )$ for each $\phi \in \bf R$.
If $0<\phi <2\pi $, then $w_1\sin (\phi )/\phi  =w_2\sin (2\pi -\phi )/
(2\pi -\phi )$ if and only if $w_1=-\phi w_2/(2\pi -\phi ).$
To exclude this ambiguity we put in Formula $(3.3)$
$\phi =|M|\ge 0$ and $w_{i_1}\ge 0$.
Therefore, $Ln(\exp (z))=z$ on ${\cal A}_r\setminus Q$, hence
using Formulas $(3.4,3.5)$ we obtain the multivalued function
$$(3.8)\quad Ln (z)=ln (|z|)+Arg(z), \mbox{ where }
Arg(z):=arg (z)+2\pi aM$$
on ${\cal A}_r\setminus \{ 0 \} $, where $ln$ is the usual
real logarithm on $(0,\infty )$, $a\in \bf Z$,
$$|z|\exp (2\pi arg (z))=z,\quad arg (z):=\sum_{s\in \hat b}w_{s,z}s,\quad
(w_{s,z}: s\in {\hat b})\in {\bf R^{2^r-1}},$$
$\sum_{s\in \hat b}w_{s,z}^2\le 1$, $w_{i_1,z}\ge 0$,
$M=\sum_{s\in \hat b}w_ss$ is any unit vector (that is, $|M|=1$) in
${\cal I}_r$ commuting with
$arg (z)\in {\cal A}_r$, $arg (z)$ is uniquely defined by such
restriction on $(w_s: s\in {\hat b})$, for example, $M=\zeta arg (z)$
for any $\zeta \in \bf R$, when $arg (z)\ne 0$.
\par For each fixed
$M\in {\cal I}_r$ $\quad \exp (aM)$ is a one-parameter family of
special transformations of ${\cal A}_r$, that is, $\exp (aM)\eta \in
{\cal A}_r$ for each $\eta \in {\cal A}_r$ and $| \exp (aM) | =1$,
where ${\cal A}_r$ as the linear space
over $\bf R$ is isomorphic with $\bf R^{2^r}$. On the other hand, there
are special transformations of ${\cal A}_r$ for which
$a=\pi /2+\pi k$, but $M$ is variable with $|M|=1$, where
$k\in \bf Z$, then $\exp (z)=(-1)^k M.$
To each closed curve $\gamma $ in ${\cal A}_r$ there corresponds
a closed curve $P_{\xi }(\gamma )$ in a $\bf R$-linear subspace
$\xi \ni 0$, where $P_{\xi }$ is a projection on $\xi $, for example,
\par $P_{{\bf R}\oplus {\bf R}s}(z)=(z-s(zs))/2=w_1+w_ss$
for $\xi ={\bf R}\oplus {\bf R}s$,
\par $P_{{\bf R}s\oplus {\bf R}p}(z)=
sP_{{\bf R}\oplus {\bf R}s^*p}(s^*z)=[z-p(s^*z(s^*p))]/2=w_ss+w_pp$ \\
for each $s\ne p\in {\hat b}$.
Particular cases of such special transformations also correspond to
$w_s=0$ for some $s\in \hat b$ for $M\ne 0$. To each closed curve
$\gamma $ in ${\cal A}_r$ and each $a$ and $b$ in ${\cal A}_r$ with
$ab\ne 0$ there corresponds a closed curve $(a\gamma )b$ in ${\cal A}_r$.
\par Instead of the Riemann two dimensional surface of the
complex logarithm function we get the $2^r$-dimensional
manifold $W$, that is, a subset of $Y^{\aleph _0}:=
\prod_{i\in \bf Z}Y_i,$ where $Y_i= Y$ for each $i$, such that
each $Y$ is a copy of ${\cal A}_r$ embedded into ${\cal A}_r
\times \bf R^{2^r-1}$ and cut by a $(2^r-1)$-dimensional submanifold
$Q$ and with diffeomorphic bending of a neighbourhood of $Q$
such that two $(2^r-1)$-dimensional edges $\mbox{ }_1Q$ and
$\mbox{ }_2Q$ of $Y$ diffeomorphic
to $Q$ do not intersect outside zero, $\mbox{ }_1Q\cap \mbox{ }_2Q
= \{ 0 \} ,$
that is, the boundary $\partial Q$ is also cut everywhere outside
zero. We have $\partial Q= \bigcup_{j\in \hat b} \partial Q^j$,
where $\partial Q^j:=\{ \theta (w_s: s\in {\hat b}): w_j=0,
(w_s: s\in {\hat b})\in S_0\cup \bigcup_{q,j_1,...,j_q}S_{j_1,...,j_q} \} $,
for each $j\in \hat b$.
To exclude rotations in each subspace $v+a\xi $
isomorphic with $\bf R^2$ and embedded into ${\bf R}+\partial Q^j$,
where $\xi \in \partial Q^j$, $\xi \ne 0$, we have cut $\partial Q$.
Then in ${\cal A}_r\times \bf R^{2^r-1}$ two copies
$Y_i$ and $Y_{i+1}$ are glued by the equivalence relation
of $\mbox{ }_2Q_i$ with $\mbox{ }_1Q_{i+1}$ via the segments
$\{ a_{l,i}M: a_{l,i} \in {\bar {\bf R}}_+  \} $
such that $a_{1,i+1}=a_{2,i}$ for each $a_{l,i}\in {\bar {\bf R}}_+$
and each given real $(w_s: s\in {\hat b})\in P$ with $M=\sum_{s\in \hat b}
w_ss$, $|M|=1$.
This defines the $2^r$-dimensional manifold $W$ embedded into ${\cal A}_r
\times \bf R^{2^r-1}$ and $Ln: {\cal A}_r\setminus \{ 0 \} \to W$ is the
univalued function, that is, $Ln(z)$ is a singleton
in $W$ for each $z\in {\cal A}_r\setminus \{ 0 \}$.
\par In the case of ${\cal A}_{\infty }$ consider a family
$\Upsilon $ of subsets of ${\hat b}={\hat b}_{\infty }$
such that if $A\in \Upsilon $,
then ${\hat b}\setminus A\notin \Upsilon $, put $Q:= \{ M\in
{\cal I}_{\infty } :$ $w_s\le 0$ $\forall s\in A,$ $w_s\ge 0$ $\forall
s\in {\bf b}\setminus A,$ $A\in \Upsilon \} $,
where $M=\sum_{s\in \hat b}w_ss$, $w_s\in \bf R$ for each $s\in \bf b$.
Then $\partial Q=\bigcup_{j\in \hat b}\partial Q^j$,
where $\partial Q^j:= \{ z\in Q:$ $w_j=0 \} $.
Consider $Y^{\aleph _0}=\prod_{k\in \bf N}Y_k$, where $Y_k=Y$
is a copy of ${\cal A}_{\infty }$ embedded into ${\cal A}_{\infty }\times
l_2({\bf R})$ and cut by the infinite dimensional submanifold $Q$
and with diffeomorphic bending of a neighbourhood of $Q$ such that
infinite dimensional edges $\mbox{ }_1Q$  and $\mbox{ }_2Q$
do not intersect outside zero, $\mbox{ }_1Q\cap \mbox{ }_2Q= \{ 0 \} $.
Then in ${\cal A}_{\infty }\times l_2({\bf R})$ two copies
$Y_k$ and $Y_{k+1}$ are glued by the equivalence relation
of $\mbox{ }_2Q_k$ with $\mbox{ }_1Q_{k+1}$ via the segments
$\{ a_{l,k}M: a_{l,k} \in {\bar {\bf R}}_+  \} $
such that $a_{1,k+1}=a_{2,k}$ for each $a_{l,k}\in {\bar {\bf R}}_+$
and each given $M=\sum_{s\in \hat b} w_ss \in Q\cap
S(0,1,{\cal I}_{\infty })$ with $|M|=1$, where
$S(z_0,\rho ,{\cal I}_{\infty }):=$
$\{ z\in {\cal I}_{\infty }:$ $|z-z_0|=\rho \} $, $\rho >0$.
This defines the infinite dimensional manifold $W$ embedded
into ${\cal A}_{\infty } \times l_2({\bf R})$ and
$Ln: {\cal A}_{\infty }\setminus \{ 0 \} \to W$ is the
univalued function.
\par {\bf 3.8. Remark.} In the real case trigonometric
and hyperbolic functions are different, but defined as functions
of the ${\cal A}_r$-variable they are related. Put 
\par $\cos (v):=[\exp (vM) + \exp (-vM)]/2$,
\par $\sin (v):=[\exp (vM) - \exp (-vM)]M^*/2$,
\par $\cosh (v) := [\exp (v) + \exp (-v) ]/2$,
\par $\sinh (v) := [\exp (v) - \exp (-v) ]/2$ \\
for each $v\in \bf R$, $M\in {\cal I}_r$, $|M|=1$,
$2\le r \le \infty $ and
\par $\cos (z) := [\exp (z) + \exp (-z) ]/2$ and
\par $\sin (z) : = [ \exp (z) - \exp (-z) ] (z^* - z )/[2 |z-z^*|]$
for each $z\in {\cal A}_r$, where $z\ne z^*$ in the latter case, then
\par $\cos (v+yM)=\cos (v) \cosh (y) - \sin (v) \sinh (y) M$,
\par $\sin (v+yM)= \sin (v) \cosh (y) + \cos (v) \sinh (y) M$ \\
for each $v, y \in \bf R$ and $M$ as above.
\par {\bf 3.8.1. Proposition.} {\it Let suppositions of Proposition $2.2.1$
be satisfied with $k=n=m$ and $f\circ g(z)=z$ for each $z\in U$,
where $g: U\to W$ is a bijective surjective mapping, that is,
there exists the inverse mapping $f=g^{-1}$. Then
\par $(Df(g)|_{g=g(z)}).\eta
=(Dg(z))^{-1}.\eta $ for each $\eta \in {\cal A}_r^n$, $2\le r\le \infty $.
If $f$ is either $z$ or $\tilde z$ or $(z,{\tilde z})$-superdifferentiable,
then $g$ is either $z$ or $\tilde z$ or
$(z,{\tilde z})$-superdifferentiable respectively.}
\par {\bf Proof.} In view of Proposition $2.2.1$ $(Df(g).
((Dg(z)).h)=h$ for each $z\in U$ and each $h\in {\cal A}_r^n$.
Since $Dg(z)$ is $\bf R$-homogeneous and ${\cal A}_r^n$-additive,
then $DG$ has an inverse operator, where $G=g\circ \sigma $ is
a diffeomorphism of a real domain $P$ in $\bf R^{2^rn}$ corresponding
to $U$, $\sigma (P)=U$. Therefore, there exists an inverse
$\bf R$-homogeneous ${\cal A}_r^n$-additive operator
$(Dg(z))^{-1}$ on ${\cal A}_r^n$. Putting $\eta =(Dg(z)).h$ we
get the statement of this proposition.  The latter statement of this
proposition follows from Proposition $2.3$, since $g$ is either
$z$ or $\tilde z$ or $(\mbox{ }_1z,\mbox{ }_2z)
|_{\mbox{ }_1z=z,\mbox{ }_2z=\tilde z}$-represented together with $f$.
\par {\bf 3.8.2. Proposition.} {\it Let $g: U\to {\cal A}_r$,
$\infty \ge r\ge 3$, be an ${\cal A}_r$-holomorphic function on $U$,
where $U$ is open in ${\cal A}_r$. If $0\notin g(U)$, then
there exists the ${\cal A}_r$-holomorphic function $f=1/g$
on $U$ such that
\par $(Df(z)).h=[D(1/g)|_{g=g(z)}].((Dg(z)).h)$ for each
$h\in {\cal A}_r$. \\
In particular, $(Df(z)).h=- \xi [((Dg(z)).h)\xi ]$ for $r=3$,
${\cal A}_3=\bf K$, where $\xi :=1/g$.}
\par {\bf Proof.} In view of Formulas $(3.2, 3.3)$ there exists
$f=1/g$, that is $f(z)g(z)=g(z)f(z)=1$ for each $z\in U$, but
it does not mean an existence of a solution of the equation
$(ab)z=a$ in general in ${\cal A}_r$ for $r\ge 4$.
Using Proposition $2.2.1$ we get
\par $(Df(z)).h=(D(1/g)|_{g=g(z)}).((Dg(z)).h)$, \\
such that $g((D(1/g)).h)=-((Dg).h)/g=-h/g$ for each $h\in {\cal A}_r$.
That is $D(1/g)$ exists. If $r=3$, ${\cal A}_3=\bf K$, then
$(D(1/g)).h=-\xi (h\xi )$, where $\xi :=1/g$, since $\bf K$
is alternative.
\par {\bf 3.8.3. Theorem.} {\it The function $Ln$ is
${\cal A}_r$-holomorphic on any domain $U$ in ${\cal A}_r$ obtained by an
${\cal A}_r$-holomorphic diffeomorphism of ${\cal A}_r\setminus Q$
onto $U$, where $2\le r\le \infty $.
Each path $\gamma $ in ${\cal A}_r$ such that $\gamma (t)=\rho
\exp (2\pi tnM)$ with $t\in [0,1]$, $n\in {\bar {\bf R}}_+$,
$M\in {\cal I}_r$, $|M|=1$ is closed in ${\cal A}_r$ if and only if
$n\in \bf N$, where $\rho >0$. In this case
$$(3.9)\quad \int_{\gamma }z^{-1}dz=
\int_{\gamma }d(Ln z)=2\pi nM.$$}
\par {\bf Proof.} In view of Formulas $(3.2, 3.3)$ each $0\ne z\in
{\cal A}_r$ has $z^{-1}$. If $U$ and $V$ are two open subsets in
${\cal A}_r$ and $g: V\to U$ is an ${\cal A}_r$-holomorphic
diffeomorphism of $V$ onto $U$
and $f$ is an ${\cal A}_r$-holomorphic function on $V$,
then $f\circ g^{-1}$ is ${\cal A}_r$-holomorphic function
on $U$, since $(f\circ g^{-1})'(z).h=(f'(\zeta )|_{\zeta =
g^{-1}(z)}.(g^{-1}(z))'.h$ for each $z\in U$ and each
$h\in {\cal A}_r$ (see Propositions $2.2.1$ and $3.8.1$).
Since $\exp $ is the diffeomorphism
from ${\cal A}_r\setminus Q$ onto ${\cal A}_r\setminus \exp (Q)$,
we have that $Ln $ is ${\cal A}_r$-holomorphic on  ${\cal A}_r\setminus Q$
due to Propositon $3.8.1$ 
and on each of its ${\cal A}_r$-holomorphic images after choosing
a definite branch of the multivalued function $Ln (z)$ (see Formula $(3.8)$).
\par A path $\gamma $ is defined for each $t\in \bf R$
not only for $t\in [0,1]$ due to the existence of $\exp $.
In view of Formulas $(3.2, 3.3)$ a path $\gamma $ is closed
(that is, $\gamma (t_0)=\gamma (t_0+1)$ for each $t_0\in \bf R$)
if and only if $\cos (2\pi n)=\cos (0)=1$ and $\sin (2\pi n)=0$, that is,
$n\in \bf N$.
\par From the definition of the line integral
we get the equality: $\int_{\gamma }d(Lnz)=
\int _0^1(Ln z)'.(\gamma '(t)dt)$. Considering integral sums
by partitions $P$ of $[0,1]$ and taking the limit by the family
of all $P$ we get, that $\int_{\gamma }d(Ln z)=Arg(\gamma (1))-
Arg (\gamma (0))$ for a chosen branch of the function $Arg (z)$
(see Formula $(3.8)$). Therefore, $\int_{\gamma }d(Ln z)=2\pi nM.$
\par Since ${\cal A}_r$ is power-associative, then
$z$ commutes with itself and we have:
$\exp (z)'.z=\exp (z)z.$ Therefore,
$\exp (Ln(z))'.1=(d \exp (\eta )/d \eta )|_{\eta =Ln(z)}
.(Ln(z))'.1=\exp (Ln(z))(Ln(z))'.1$, consequently,
$(Ln (z))'.1=\exp (-Ln(z))=z^{-1}$ and inevitably
$$\lim_PI(z^{-1},\gamma ;P)=\lim_P\sum_l{\hat z}_l^{-1}\Delta z_l
=\lim_P\Delta Ln(z_l)=\int_{\gamma }dLn(z),$$
hence $\int_{\gamma }z^{-1}dz=\int_{\gamma }dLn(z)$.
That is, $\int_{\gamma }d Ln (z)$ can be considered
as the definition of $\int_{\gamma }z^{-1}dz$.
\par {\bf 3.8.4. Notation.} Denote an ordered compositions of functions
$ \{ f_1\circ f_2\circ ... f_m \} _{q(m)}$, where $q(m):=
(q_m,...,q_3)$, $q_m\in \bf N$ means that the first (the most
internal bracket of) composition is $f_{q_m}\circ f_{q_{m+1}}=
f_{q_m}(f_{q_{m+1}})=:[f_{q_m}\circ f_{q_{m+1}}]$ such that
to the situation \\
$(f_1\circ ... \circ [f_t\circ f_{t+1}]\circ ...
\circ [f_w\circ f_{w+1}]\circ ...\circ f_m)$ with two simultaneous
independent compositions, but $t<w$ by our definition of ordering
there corresponds $q_m=w$ (apart from the multiplication). After the first
composition we get the composition of $(\mbox{ }'f_1\circ ...\circ
\mbox{ }'f_{m-1})$, where not less than $(m-2)$ elements here are the same
as in the first composition, then $q_{m-1}$ corresponds to the first
composition in this new ordered family of functions, and so on by
induction from $j$ to $j-1$, $j=m,m-1,...,3$.
Since $q_2$ and $q_1$ are unique, we omit them.
After steps $q_m,...,q_{m-j+1}$ let the corresponding composition
be denoted by \\
$ \{ \mbox{ }^jg_1\circ ...\circ \mbox{ }^jg_{m-j} \} _{q(m-j)}$, where \\
$\mbox{ }^jg_1,..., \mbox{ }^jg_{m-j}$ are resulting composites
on preceding steps, some of them may belong to the set
$ \{ f_1,...,f_m \} $. If $f_l$ is in the composite
on the $k=k(l)$ step, then there may be two variants:
$f_l\circ \mbox{ }^lg_p$ or $\mbox{ }^lg_p\circ f_l$, where
$p=p(l)$, $j=j(l)=k(l)-1$, $\mbox{ }^0g_p:=f_p$. In the first case
suppose that $dom (f_l)\subset \mbox{ }^jg_p(dom \mbox{ }^jg_p)$,
in the second case let $dom (\mbox{ }^jg_p)\subset f_l(dom (f_l))$
for each $l=1,...,m$.
\par {\bf 3.8.5. Proposition.} {\it Let $f_1$,...,$f_m$ be a family
of ${\cal A}_r$ superdifferentiable functions either all by $z$ or
all by $\tilde z$ or all by $(z,{\tilde z})$, $f_j: U_j\to
{\cal A}_r^{t(j)}$, $U_j$ is open in ${\cal A}_r^{t(j+1)}$, $t(j)
\in \bf N$ for each $j=1,...,m$ such that their domains satisfy
conditions above, where $2\le r\le \infty $. Then
$$(i)\quad (D \{ f_1\circ f_2\circ ... \circ f_m \} _{q(m)}(z)).h=
\{ Df_1(\mbox{ }^{j(1)}g_{p(1)}).Df_2(\mbox{ }^{j(2)}g_{p(2)}).
...(Df_m(z)).h \} _{q(m)}$$
for each $h\in {\cal A}_r^{t(m+1)}$, where $Df_l(\mbox{ }^{j(l)}
g_{p(l)}).\xi =(Df_l(\eta )|_{\eta =\mbox{ }^{j(l)}g_{p(l)}(z)}).\xi $
for each $\xi \in {\cal A}_r^{t(l+1)}$.
Moreover, $\{ f_1\circ f_2\circ ... \circ f_m \} _{q(m)}$
is superdifferentiable either by $z$ or by $\tilde z$ or by
$(z,{\tilde z})$ correspondingly.
If $r=2$, ${\cal A}_2=\bf H$, then the composition in $(i)$ is associative,
for each $r\ge 3$ it can be in general nonassociative. At a marked point
$z=a\in U_m$ it takes the form:
$$(ii)\quad (D \{ f_1\circ f_2\circ ... \circ f_m \} _{q(m)}(z)).h=
\{ Df_1(\eta _1).Df_2(\eta _2). ...(Df_m(z)).h \} _{q(m)},$$
where $\eta _l:=f_l(f_{l+1}(...f_{m-1}(f_m(z))...))$ for each
$l=1,...,m-1.$}
\par {\bf Proof.} For $m=2$ this proposition was proved in \S 2.2.1.
Prove this proposition by induction and applying Proposition 2.2.1
to pairs of functions in appearing compositions.
At first mention that the order of composition for the differential
is essential for ${\cal A}_r$ with $r\ge 3$. The particular case of
all right superdifferentiable functions shows, that in general
$(Df_1.Df_2).Df_3$ is not equal to $Df_1.(Df_2.Df_3)$, since
these operators are right superlinear, but multiplication of matrices
with entries in the Cayley-Dickson algebra ${\cal A}_r$
with $r\ge 3$ is not associative. Moreover, these expressions may
be different, when $Df_j$ are not right superlinear, but only
${\cal A}_r^{t(j+1)}$-additive.
\par Let proposition be proved for all $n\le m$, consider
$\{ f_1\circ ... \circ f_{m+1})_{q(m+1)}$. In it $f_1$ is in the
composition on the $k(1)$ step with $\mbox{ }^{j(1)}g_{p(1)}$
of the type $f_1\circ \mbox{ }^{j(1)}g_{p(1)}$, since $f_1$
is in the extreme left position.
If $k(1)=1$, then $\mbox{ }^{j(1)}g_{p(1)}=f_2$ and \\
$\{ f_1\circ ...\circ f_{m+1} \} _{q(m+1)}= \{ [f_1\circ f_2]\circ ...
\circ f_{m+1} \} _{q(m)}$.  \\
Applying supposition of induction to composition of functions
$[f_1\circ f_2]$, $f_3$,...,$f_{m+1}$ and then substituting to it
expression of $D(f_1\circ f_2)$ by Proposition $2.2.1$ we get
the statement of this proposition in the case $k(1)=1$.
If $k(1)>1$, then on the $k(1)$ step the composition of
$f_1$ and $\mbox{ }^{j(1)}g_2$,...,$\mbox{ }^{j(1)}g_{m+1-j(1)}$
is considered, where $j(1)=k(1)-1\ge 1$, hence $m+1-j(1)\le m$.
Applying supposition of induction to $\{ f_1\circ 
\mbox{ }^{j(1)}g_2\circ ...\circ \mbox{ }^{j(1)}g_{m+1-j(1)}
\} _{q(m+1-j(1))}$ and then to each $\mbox{ }^{j(1)}g_p$ while being
a nontrivial composition we get the statement of proposition
in the case $k(1)>1$. The last statement also follows by the
considered above induction from Proposition $2.2.1$.
\par Set-theoretical composition of functions is independent
from brackets, but it depends only on order of functions:
$(f_1\circ f_2)\circ f_3(z)=(f_1\circ f_2)(f_3(z))=
f_1(f_2(f_3(z)))=f_1\circ (f_2\circ f_3(z))$, etc. by induction.
A non-associativity in general appears after superdifferentiation
over ${\cal A}_r$ with $r\ge 3$. Since $\{ f_1\circ f_2\circ
...\circ f_m \} _{q(m)}(z)=f_1(f_2(...(f_m(z))...)$, then
for a marked point $z=a\in U_m$ Formula $(i)$ takes the form $(ii)$.
\par In the case of $r=2$, ${\cal A}_2=\bf H$, each $\bf R$-homogeneous
$\bf H$-additive operator $A$ has the form
\par $A.h=\sum_jB_jhC_j$ for each $h\in \bf H^n$, where the sum by
$j$ is finite, $1\le j\le 4$ (see formulas in the proof of Theorem 3.28
\cite{luoyst}), where $B_j$ is the $n\times n$ matrix,
$h$ is the $n\times 1$ matrix, $A_j$ is the $1\times 1$ matrix
with entries in $\bf H$. Let $A_k$ be an operator corresponding to
$Df_k(\eta _k)$ for a given marked point $z=a\in U_3$,
$A_k.h=\sum_jB_{j,k}hC_{j,k}$ for each $k=1, 2, 3$,
then $(A_1.A_2).A_3.h=\sum_{j_1,j_2,j_3}(B_{j_1,1}B_{j_2,2})(B_{j_3,3}h
C_{j_3,3})(C_{j_2,2}C_{j_1,1})$
$=\sum_{j_1,j_2,j_3}B_{j_1,1}(B_{j_2,2}B_{j_3,3})h
(C_{j_3,3}C_{j_2,2})C_{j_1,1}$ $=(A_1.(A_2.A_3)).h$,
since the matrix multiplication over $\bf H$ is associative.
Applying the latter formula by induction, we get, that
in the case of $\bf H$ the composition in Formula $(i)$ is associative.
\par {\bf 3.9. Theorem.} {\it Let $f$ be an ${\cal A}_r$-holomorphic
function on an open domain $U$ in ${\cal A}_r$, $\infty \ge r\ge 3$.
If $(\gamma +z_0)$ and $\psi $ are presented as piecewise unions of paths
$\gamma _j+z_0$ and $\psi _j$ with respect to parameter
$\theta \in [a_j,b_j]$ and $\theta \in [c_j,d_j]$
respectively with $a_j<b_j$ and $c_j<d_j$ for each $j=1,...,n$
and $\bigcup_j[a_j,b_j]=\bigcup_j[c_j,d_j]=[0,1]$ 
homotopic relative to $U_j\setminus \{ z_0 \} $, where
$U_j\setminus \{ z_0 \} $ is a $(2^r-1)$-connected open domain in
${\cal A}_r$ such that $\pi _{s,p,t}(U_j\setminus \{ z_0 \})$
is simply connected in $\bf C$ for each $s=i_{2k}$, $p=i_{2k+1}$,
$k=0,1,...,2^{r-1}-1$ ($\forall 0\le k\in \bf Z$ and $P_m(U_j\setminus
\{ z_0 \} )$ is $(2^m-1)$-connected for each $4\le m\in \bf N$
if $r=\infty $), each $t\in {\cal A}_{r,s,p}$ and $u\in {\bf C}_{s,p}$
for which there exists $z=t+u\in {\cal A}_r$.
If $(\gamma +z_0)$ and $\psi $ are closed rectifiable paths (loops)
in $U$ such that
$\gamma (\theta )=\rho \exp (2\pi \theta M)$ with $\theta \in [0,1]$
and a marked $M\in {\cal I}_r$, $|M|=1$ and $z_0\notin \psi $. Then
$$(3.10)\quad (2\pi )f(z)M=\int_{\psi }f(\zeta )(\zeta -z)^{-1}d\zeta $$
for each $z\in U$ such that $|z-z_0|< \inf_{\zeta \in \psi ([0,1])}
|\zeta -z_0|$. If ${\cal A}_r$ is alternative, that is, $r=3$, ${\cal A}_3=
\bf K$, or $f(z)\in \bf R$, then
$$(3.11)\quad f(z)=(2\pi )^{-1}(\int_{\psi }f(\zeta )(\zeta -z)^{-1}d\zeta )
M^*.$$}
\par {\bf Proof.} Join $\gamma $ and $\psi $ by a rectifiable path
$\omega $ such that $z_0\notin \omega $,
which is going in one direction and the opposite direction,
denoted $\omega ^-$, 
such that $\omega _j\cup \psi _j\cup \gamma _j\cup \omega _{j+1}$
is homotopic to a point relative to $U_j\setminus \{ z_0 \}$
for suitable $\omega _j$ and $\omega _{j+1}$, where $\omega _j$ joins
$\gamma (a_j)$ with $\psi (c_j)$ and $\omega _{j+1}$ joins
$\psi (d_j)$ with $\gamma (b_j)$ such that $z$ and $z_0\notin
\omega _j$ for each $j$. Then
$\int_{\omega _j} f(\zeta )(\zeta -z)^{-1}d\zeta = -
\int_{\omega _j^-}f(\zeta )(\zeta -z)^{-1}d\zeta $
for each $j$. In view of Theorem 2.15 there is the equality
$-\int_{\gamma ^-+z}f(\zeta )(\zeta -z)^{-1}d\zeta =
\int_{\psi }f(\zeta )(\zeta -z)^{-1}d\zeta .$
Since $\gamma +z$ is a circle around $z$ its radius $\rho >0$ can be
chosen so small, that $f(\zeta )=f(z)+ \alpha (\zeta ,z)$,
where $\alpha $ is a continuous function on $U^2$ such that
$\lim_{\zeta \to z}\alpha (\zeta ,z)=0$, then
$\int_{\gamma +z}f(\zeta )(\zeta -z)^{-1}d\zeta =$
$\int_{\gamma +z}f(z)(\zeta -z)^{-1}d\zeta +\delta (\rho )=$
$2\pi f(z) M + \delta (\rho ),$
where $|\delta (r)|\le |\int_{\gamma }\alpha (\zeta ,z)
(\zeta -z)^{-1}d\zeta |$ $\le 2\pi \sup_{\zeta \in \gamma }
|\alpha (\zeta ,z)| C_1 \exp (C_2 \rho ^6)$,
where $C_1$ and $C_2$ are positive constants (see Inequality
$(2.7.4)$), hence there exists $\lim_{\rho \to 0,
\rho >0}\delta (\rho )=0$. Taking the limit while $\rho >0$ tends to zero
yields the conclusion of this theorem. If $r=3$, that is,
${\cal A}_r=\bf K$, or $f(z)\in \bf R$, then $((2\pi )f(z)M)M^*=2\pi f(z)$.
\par {\bf 3.9.1. Corollary.} {\it Let $f$, $U$, $\psi $,
$z$ and $z_0$ be as in Theorem 3.9, then
$$|f(z)|\le \sup_{(\zeta
\in \psi , h\in {{\cal A}_r}, |h|\le 1)} |{\hat f} (\zeta ).h|.$$}
\par {\bf 3.9.2. Theorem.} {\it Let $\{ f_n: n\in {\bf N} \} $
be a sequence of ${\cal A}_r$-holomorphic functions on a neighbourhood
$W$ of bounded canonical closed subset $U$ in ${\cal A}_r$,
$2\le r\le \infty $,
such that $Int (U)$ satisfies conditions of Theorem $3.9$
and there exists $\delta >0$ for which $\{ f_n: n \in {\bf N} \} $
converges uniformly on the $\delta $-neighbourhood of the topological
boundary $Fr (U)^{\delta }$ of $U$. Then $\{ f_n: n\in {\bf N} \} $
converges uniformly on $U$ to an ${\cal A}_r$-holomorphic function
on $Int (U)$.}
\par {\bf Proof.} In view of \S 2.7 the sequence ${\hat f}_n$
is uniformly converging on $Fr (U)^{\delta }\cap W$.
Consider sections of $W$ by planes $i_{2m}{\bf R}\oplus
i_{2m+1}\bf R$ for each $m=0,1,2,...$ and rectifiable loops $\gamma $ in
$Fr (U)^{\delta }\cap W\cap (i_{2m}{\bf R}
\oplus i_{2m+1}{\bf R})$. For $r>3$ consider
all possible embeddings of $\bf K$ into ${\cal A}_r$
and such copies of $\bf K$ contain all corresponding loops $\gamma $.
In view of Estimate $2.7.(4)$ and Theorem $3.9$
the sequence $\{ f_n: n\in {\bf N} \} $ converges uniformly on $U$
to a holomorphic function, since $U$ is of finite diameter.
\par {\bf 3.10. Theorem.} {\it Let $f$ be a continuous function
on an open subset $U$ of ${\cal A}_r$, $3\le r\le \infty $.
If $f$ is ${\cal A}_r$-integral holomorphic on $U$, then
$f$ is ${\cal A}_r$ locally $z$-analytic on $U$.}
\par {\bf Proof.} Let $z_0\in U$ be a marked point and let
$\Gamma $ denotes the family of all rectifiable paths
$\gamma : [0,1]\to U$ such that $\gamma (0)=z_0$, then
$U_0=\{ \gamma (1): \gamma \in \Gamma \} $ is a connected component
of $z_0$ in $U$. Therefore, $g= \{ \gamma (1), \int_{\gamma }
f(z)dz \} $ is the function with the domain $U_0$.
From Formulas $(3.2, 3.3)$ it follows, that each $0\ne z\in
{\cal A}_r$ has an inverse $z^{-1}z=zz^{-1}=1$ (this does not
mean an existence of a solution of the equation $(ab)z=a$ with $b\ne 0$
in general in ${\cal A}_r$ for $r\ge 4$). In view of Proposition
$2.2.1$ $\partial _z(\zeta -z)^{-k}.1=k(\zeta -z)^{-k-1}$, since
$0=(\partial _z1).h=$  $(\partial _z z^kz^{-k}).h=$
$((\partial _zz^k).h)z^{-k}+z^k((\partial _zz^{-k}).h)$ for each
$h\in {\cal A}_r$ and $z\ne 0$, $z\in {\cal A}_r$.
As in \S 2.15 it can be proved, that
$F(z):=\int_{\gamma }f(z)dz$, for each rectifiable $\gamma $ in
$U$, depends only on initial and final points.
This integral is finite, since $\gamma ([0,1])$ is contained
in a compact canonical closed subset $W\subset U$
on which $f$ is bounded.
Therefore, $(\partial \int_{z_0}^zf(\zeta )d\zeta /\partial z).h=
{\hat f}(z).h$ for each $z\in U$ and $h\in {\cal A}_r$,
$(\partial \int_{z_0}^zf(\zeta )d\zeta /\partial {\tilde z})=0$
for each $z\in U$ and $h\in {\cal A}_r$, where
$z_0$ is a marked point in $U$ such that $z$ and $z_0$ are in one
connected component of $U$, since $\int_{z_0}^{z+\Delta z}f(\zeta )d\zeta-
\int_{z_0}^zf(\zeta )d\zeta =$  ${\hat f}(z).\Delta z+\epsilon (\Delta z)
|\Delta z|$, where $\lim_{\Delta z\to 0}\epsilon (\Delta z)=0$
(see \S 2.5). In particular, ${\hat f}(z).1=f(z)$ for each $z\in U$.
Here $\hat f$ is correctly defined for each $f\in C^{1,0}(U,{\cal A}_r)$
(see Corollary $2.15.1$)
by continuity of the differentiable integral functional
on $C^0(U,{\cal A}_r)$. For a given $z\in \bf U$ choose
a neighbourhood $W$ satisfying the conditions of Theorem 3.9.
Then there exists a rectifiable path $\psi \subset W$ such that $f(z)$
is presented by Formula $(3.10)$. The latter integral is
infinite differentiable by $z$ such that
$$(3.12)\quad 2\pi ((\partial ^k F(z)/\partial z^k).h)M=
(\int_{\psi }F(\zeta )(\partial ^k (\zeta -z)^{-1}/
\partial z^k).h) d\zeta ),$$
where $h\in {\cal A}_r^k$, particularly, for $h=(1,...,1)=:1^{\otimes k}$:
$$(3.13)\quad 2\pi ((\partial ^k F(z)/\partial z^k).1^{\otimes k})M=
k! (\int_{\psi }F(\zeta ) (\zeta -z)^{-k-1} d\zeta ) ,$$
where $M\in {\cal I}_r$, $|M|=1$. For simplicity of notation
we can omit $1^{\otimes k}$ on the left of $(3.13)$.
In particular, we may choose a ball $W=B(a,R,{\cal A}_r):=
\{ \xi \in {\cal A}_r: |\xi -a|\le R \} \subset U$
for a sufficiently small $R>0$ and $\psi =\gamma +a$, where
$\gamma (s)=\rho \exp (2\pi tM)$ with $t\in [0,1]$, $0<\rho <R$.
If we prove that $F(z)$ is ${\cal A}_r$ locally $z$-analytic, then
evidently its $z$-derivative $f(z)$ will also be  
${\cal A}_r$ locally $z$-analytic.
If $r\ge 4$, then there exist different embeddings of
the octonion algebra $\bf K$ into ${\cal A}_r$ (see \S 3.6.2).
Suppose there is a series $f(z) := \sum_j a_j (z-z_0)^j b_j$
converging for $|z-z_0| < \rho $, where $a_j, b_j$ for each $j$
belong to the same $\Upsilon _a\hookrightarrow \bf C$,
$0\ne a\in {\cal A}_r$,
expansion coefficients do not depend on the embedding of
$\Upsilon _{a,z-z_0}$ into ${\cal A}_r$, while $z-z_0$ is varying
within the same copy of $\bf K$. Then this expansion
is valuable for each $z\in {\cal A}_r$ with $|z-z_0|<\rho $,
since different embeddings of $\bf K$ (with generators
$\{ 1,M_1,...,M_7 \} \hookrightarrow {\cal A}_r$ such that
$|M_i|=1$, $Re (M_iM_j)=0$ and $|M_iM_j|=|M_i| |M_j|$ for
each $i\ne j$) into ${\cal A}_r$ give all possible values of
$z\in {\cal A}_r$ with $|z-z_0|<\rho $.
Consider $\psi $ such that $\Upsilon _{\zeta -a,z-a}$ has an
embedding into $\bf K$ for each $\zeta \in \psi $, which is not
finally restrictive due to Theorem 2.15.
In view of the latter statement using the
alternative property of $\bf K$ we
consider $z\in B(a,\rho ',{\cal A}_r)$  with $0<\rho '<\rho $,
then $|z-a|<|\zeta -a|$
for each $\zeta \in \psi $ and $(\zeta -a -(z-a))^{-1}=
(1-(\zeta -a)^{-1}(z-a))^{-1}(\zeta -a)^{-1}=$ $\sum_{k=0}^{\infty }
((\zeta -a)^{-1}(z-a))^k(\zeta -a)^{-1}$, 
where $0\notin \psi $. Therefore,
$$(3.14)\quad 2\pi F(z)M= \sum_{k=0}^{\infty }\phi _k(z),$$
$$\mbox{where }\phi _k(z):=
(\int_{\psi }F(\zeta )((\zeta -a)^{-1}(z-a))^k(\zeta -a)^{-1}d\zeta ).$$
Thus $|\phi _k(z)|\le \sup_{\zeta \in \psi }|F(\zeta )|(\rho '/\rho )^{-k}$
for each $z\in B(a,\rho ',{\cal A}_r)$ and series $(3.14)$ converges
uniformly on $B(a,\rho ',{\cal A}_r)$. Each function $\phi _k(z)$
is evidently ${\cal A}_r$ locally $z$-analytic on $B(a,\rho ',{\cal A}_r)$,
hence $F(z)M$ is such too. Since for each $a\in U$
there is an $\rho '>0$, for which the foregoing holds,
it follows that $F(z)M$ is the ${\cal A}_r$ locally $z$-analytic function.
Now write $f(z)=\sum_{s\in \bf b}f_ss$, where $f_s\in \bf R$ for
each $s\in \bf b$. If the loop $\gamma $ is nondegenerate, then
loops $s\gamma $, $\gamma s$, $(2^r-2)^{-1} \{ -\gamma +\sum_{s \in \hat b}
s(\gamma s^*) \} ={\tilde {\gamma }} $ are nondegenerate.
If $(i)\quad \int_{\gamma }f(\zeta )d\zeta =0$, then $\int_{\gamma }
sf(\zeta )d\zeta =0$ and  $\int_{\gamma }f(\zeta )s d\zeta =0$
(see Theorem $2.7$). In view of Formulas $2.8.(2)$: \\
$f_1=(f+(2^r-2)^{-1} \{ -f +\sum_{s\in {\hat b}_r} s(fs^*) \} )/2$ and \\
$f_p=(i_p(2^r-2)^{-1} \{ -f +\sum_{s\in {\hat b}_r} s(fs^*) \}
-fi_p)/2$ \\
for each $i_p\in {\hat b}_r$ (for $r=\infty $ use $\lim_{r\to \infty }$
on the right of the latter formulas) we have
that Condition $(i)$ is equivalent to:
$\int_{\gamma }f_s(\zeta )d\zeta =0$ for each $s\in \bf b$.
Therefore, the proof above shows, that each function
$2\pi f_s(z)M$ is ${\cal A}_r$ locally $z$-analytic, where
$M$ is arbitrary in ${\cal I}_r$, $|M|=1$. But $f_s(z)\in \bf R$
for each $z\in U$, hence $f_s(z)=(2\pi f_s(z)M)(2\pi )^{-1}M^*=
f_s(z)$ is ${\cal A}_r$ locally $z$-analytic for each $s\in \bf b$,
hence $f$ is also ${\cal A}_r$ locally $z$-analytic on $U$.
\par {\bf 3.11. Note.} Theorems 2.11, 2.15, 2.16, 3.10 and Corollary 2.13
establish the equivalence of notions of ${\cal A}_r$-holomorphic,
${\cal A}_r$-integral holomorphic
and ${\cal A}_r$ locally $z$-analytic classes of functions on domains
satisfying definite conditions.
\par {\bf 3.11.1. Definitions.} Let $U$ be an open subset in ${\cal A}_r$,
$2\le r\le \infty $ and
$f\in C^0(U,{\cal A}_r)$, then we say that $f$ possesses a primitive
$g\in C^1(U,{\cal A}_r)$ if $g'(z).1=f(z)$ for each $z\in U$.
A region $U$ in ${\cal A}_r$ is said to be
${\cal A}_r$-holomorphically simply connected
if every function ${\cal A}_r$-holomorphic on it possesses a primitive.
\par From \S 3.10 we get.
\par {\bf 3.11.2. Theorem.} {\it If $f\in C^{\omega }_z(U,{\cal A}_r)$,
$2\le r\le \infty $, 
where $U$ is $(2^r-1)$-connected ($P_m(U)$ is $(2^m-1)$-connected for
each $4\le m\in \bf N$ for $r=\infty $); $\pi _{s,p,t}(U)$
is simply connected in $\bf C$ for each $s=i_{2k}$, $p=i_{2k+1}$ in
$\bf b$, $t\in {\cal A}_{r,s,p}$ and $u\in
{\bf C}_{s,p}$ for which there exists $z=t+u\in U$, $U$ is an open subset in
${\cal A}_r$, then there exists $g\in C^{\omega }_z(U,{\cal A}_r)$ such that
$g'(z).1=f(z)$ for each $z\in U$.}
\par {\bf 3.11.3. Theorem.} {\it Let $U$ and $V$ be
${\cal A}_r$-holomoprhically simply connected regions in ${\cal A}_r$,
$2\le r\le \infty $ with  $U\cap V\ne \emptyset $ connected.
Then $U\cup V$ is ${\cal A}_r$-holomorphically simply connected.}
\par {\bf 3.12. Corollary.} {\it Let $U$ be an open subset in
${\cal A}_r^n$, $2\le r\le \infty $, then the family of all
${\cal A}_r$-holomorphic functions
$f: U\to {\cal A}_r$ has a structure of an ${\cal A}_r$-algebra.}
\par {\bf Proof.} If $f_1(z)=\alpha g(z)\beta +\gamma h(z)\delta $
or $f_2(z)=g(z)h(z)$ for each $z\in U$, where $\alpha $, $\beta $,
$\gamma $ and $\delta \in {\cal A}_r$ are constants, $g$ and $h$ are
${\cal A}_r$-holomorphic functions on $U$,
then $F_1$ and $F_2$ are Frech\'et
differentiable on $U$ by $(w_s: s\in {\bf b})$ (see \S 2.1 and \S 2.2)
and $\partial _{\tilde z}f_1(z)=\alpha (\partial _{\tilde z}g)\beta +
\gamma (\partial _{\tilde z}h)\delta =0$ and $\partial _{\tilde z}f_2(z)=
(\partial _{\tilde z}g)h+g(\partial _{\tilde z}h)=0$,
hence $f_1$ and $f_2$ are also ${\cal A}_r$-holomorphic on $U$.
\par {\bf 3.13. Proposition.} {\it For each complex holomorphic function
$f$ in a neighbourhood $Int (B(q_0,\rho ,{\bf C}))$, $\infty > \rho >0$,
of a point $q_0\in \bf C$ and each $2\le r\le \infty $,
$s\ne p\in {\bf b}_r$, there exists
an ${\cal A}_r$ $z$-analytic function $g$ on a neighbourhood
$Int (B(a,\rho ,{\cal A}_r))$ of $a\in {\cal A}_r$
such that $s^*g_{s,p}(u,t_0)=f(v)$ on $Int (B(u_0,\rho ,{\bf C}))$,
$u=s(Re (v)+s^*(p Im_{\bf C}(v)))$, where 
$B(x,\rho ,X):=\{ y\in X: d_X(x,y)\le \rho \}$
is the ball in a space $X$ with a metric $d_X$, $a=u_0+t_0$,
$u_0\in {\bf C}_{s,p}$, $t_0\in {\cal A}_{r,s,p}$,
$u_0=s(Re (q_0)+s^*( p Im_{\bf C}(q_0)))$, $Im_{\bf C}(v):=
(v-{\tilde v})/(2i)$.}
\par {\bf Proof.} Among Conditions $(2.3.1)$ there are independent: \\
$(3.15)\quad \partial F_1/ \partial \mbox{ }^jw_p =
\partial F_{pq^*} / \partial \mbox{ }^jw_q $, \\
$\quad \partial F_1 / \partial \mbox{ }^jw_q =
- \partial F_{pq^*} / \partial \mbox{ }^jw_p $ \\
for each $p=i_{m}$, $q=i_{m+1}\in {\bf b}_r$, $0\le m\in \bf Z$.
Consider a homogeneous polynomial function
on an open ball $Int (B)$ in ${\cal A}_r$ such that
$P_n(\lambda z)=\lambda ^n P_n(z)$ for each $\lambda \in \bf R$,
$P_n: Int (B) \to {\cal A}_r$, then $P_{n+1}$ can be written in the form \\
$(3.16)\quad P_{n+1} (z) = \sum_{s\in {\bf b}_r; k; j_1,...,j_k}
C_{s; k; j_1,...,j_k} v_1^{j_1}...v_k^{j_k}s$; \\
where $v_l := w_{i_l}$ is the real variable for each $l$,
$z= \sum_{s\in {\bf b}_r} w_ss$, $0\le j_1\in {\bf Z},...,0\le
j_k \in {\bf Z}$, $k\in \bf N$, $j_1+...+j_k=n+1$,
$C_{s; k; j_1,...,j_k} \in \bf R$ is the real expansion coefficient
for each $s, k, j_1,...,j_k$. In view of $(2.3.1)$ a function $f$ is
right superlinearly ${\cal A}_r$-superdifferentiable at a point $z_0$
if and only if $sf$ is right superlinearly
${\cal A}_r$-superdifferentiable at $z_0$ for each $s\in {\hat b}_r$. 
Then $(3.15)$ applied to $(3.16)$
gives conditions on coefficients of homogeneous polynomial
providing its right superlinear superdifferentiability: \\
$(3.17) \quad C_{1; k; j_1,...,j_m+1,...,j_k}(j_m+1)=
C_{pq^*; k; j_1,...,j_{m+1}+1,...,j_k} (j_{m+1}+1)$ and \\
$\quad C_{1; k; j_1,...,j_{m+1}+1,...,j_k}(j_{m+1}+1)=
-C_{pq^*; k; j_1,...,j_m+1,...,j_k} (j_m+1)$, \\
for each $p=i_m$, $q=i_{m+1}$ in ${\bf b}_r$.
Since $j_l+1\ge 1$ for each $l$, then coefficients on the right
sides are expressible through coefficients on the left, which can be
taken as free variables. Therefore, for each $n\ge 0$
there exists nontrivial (nonzero) $P_{n+1}$ satisfying $(3.15)$.
Evidently, the $\bf R$-linear space of all right superlinearly
${\cal A}_r$-superdifferentiable functions is infinite dimensional,
since for each $n$ there exists a nontrivial solution of
the homogeneous system of linear equations $(3.17)$.
\par Consider first an extension in the class
of ${\cal A}_r$-holomorphic functions with a superdifferential
not necessarily right superlinear on the superalgebra ${\cal A}_r$.
Since $f$ is holomorphic in $Int (B(q_0,\rho ,{\bf C}))$, it has
a decomposition $f(t)=\sum_{n=0}^{\infty }f_n(q-q_0)^n$,
where $f_n\in \bf C$ are expansion coefficients, $q\in
Int (B(q_0,\rho ,{\bf C}))$. Consider its extension in
$Int (B(z_0,\rho ,{\cal A}_r))$
such that $f(z)=\sum_{n=0}^{\infty } f_n (z-z_0)^n$,
$z_0=u_0+t_0$, $u_0=s(Re (q_0)+s^*(p Im_{\bf C}(q_0)))$,
$t_0\in {\cal A}_{r,s,p}$.  Evidently this series converges for each $z\in
Int (B(z_0,\rho ,{\cal A}))$
and this extension of $f$ is ${\cal A}_r$-holomorphic, since
$f_n\in {\cal A}_r$ for each $n$ and $\partial f/\partial {\tilde z}=0$.
\par Consider now more narrow class of quaternion holomorphic functions
with a right superlinear superdifferential on the superalgebra ${\cal A}_r$.
The Cauchy-Riemann conditions for complex holomorphic functions
are particular cases (part) of Conditions $(3.15)$.
Having grouped the series for complex holomorphic function $f$
in the series by homogeneous polynomials and applying $(3.17)$
we get expansion coefficients for the right superlinear
superdifferentiable extension of $f$ on $Int (B(z_0,\rho ,{\cal A}))$.
\par {\bf 3.14. Proposition.} {\it If $f$ is an ${\cal A}_r$-holomorphic
function on an open subset $U$ in ${\cal A}_r$, $2\le r\le \infty $,
and $ker (f'(z_0))= \{ 0 \} $ and $f'(z_0)$ is right superlinear,
then it is a conformal mapping at a marked point $z_0\in U$,
that is preserving angles between differentiable curves.
If $r=3$, ${\cal A}_3=\bf K$, then $ker (f'(z_0))= \{ 0 \} $
if and only if $f'(z_0)\ne 0$.}
\par {\bf Proof.} Let $z_0\in U$ and $f$ be right
superlinearly superdifferentiable at $z_0$.
Consider a scalar product $(*,*)_r$ in ${\bf C}^m$ as in the introduction,
$m=2^{r-1}$ for $r<\infty $ or in $l_2({\bf C})$
inherited from ${\cal A}_{\infty }$, where
$(x,y)_{l_2}=\lim_{r\to \infty }(x_r,y_r)_r$ (see \S 3.6.2 also).
To the right superlinear operator $f'(z_0)$ there corresponds
the unique bounded operator $A$ on ${\bf C}^m$ or $l_2({\bf C})$ respectively
such that $ker (A)= \{ 0 \} $, that is, $Ah=0$ if and only if
$h=0$. If $r=3$, then $A$ is invertible if and only if
$f'(z_0)\ne 0$, since $f'(z_0)\in {\cal A}_r$ and $\bf K$ is
alternative: $(ab)y=a$ has the unique solution $y$ for each
$b\ne 0$ in $\bf K$. In view of the polar decomposition $(3.5)$
$f'(z_0)=\rho \exp (M)$, where $\rho >0$ and $M\in {\cal I}_r$.
Then the adjoint operator $A^*$ corresponds to $\rho \exp (-M)$,
but $\exp (-M) \exp (M) = 1$, hence $A$ is the unitary operator:
$A\in U(m)$ or $A\in U(\infty )$ respectively.
Since the unitary group preserves scalar product, then
$f(z)$ preserves an angle $\alpha $ between each two
differentiable curves in $U$ crossing at the marked point $z_0$:
if $\psi $ and $\phi : (-1,1)\to U$ are two differentiable curves
crossing at a point $z_0\in U$, then
$f(\psi (\theta ))'=f'(z)|_{z=\psi (\theta )}.\psi '(\theta )$,
where $\cos (\alpha )=Re (\psi '(0),\phi '(0))/
(|\psi '(0)| |\phi '(0)|)$ for $\psi '(0)\ne 0$ and $\phi '(0)\ne 0$.
\par {\bf 3.14.1. Remark.} For each $r\ge 4$ the Cayley-Dickson
algebra is not the division algebra, hence the condition
$ker (f'(z_0)) = \{ 0 \}  $ is essential in Proposition $3.14$.
\par {\bf 3.15. Theorem.} {\it Let $f$ be an ${\cal A}_r$-holomorphic
function, $2\le r\le \infty $, on an open subset $U$
in ${\cal A}_r$ such that
$\sup_{z\in U, h\in B(0,1,{\cal A}_r)} |[f(z)(\zeta -z)^{-2}]^{\hat .}.h|
\le C/|\zeta -z|^2$ for each $\zeta \in {\cal A}_r\setminus cl (U)$.
Then $|f'(z)|\le C/d(z)$ for each $z\in U$,
where $d(z) := \inf_{\zeta \in {\cal A}_r\setminus U}|\zeta -z|$;
$|f(\xi )-f(z)|/|\xi -z|\le 2C/\rho $ for each $\xi $ and $z\in
B(a,\rho /2,{\cal A}_r)\subset Int (B(a,\rho ,{\cal A}_r))\subset U$, where
$\rho >0$. In particular, if $f$ is an ${\cal A}_r$-holomorphic 
function with bounded $[f(z)(\zeta -z)^{-2}]^{\hat .}.h|\zeta -z|^2$ on
${\cal A}_r^2\times B(0,1,{\cal A}_r)$ with $|\zeta |\ge 2|z|$,
that is, $\sup_{\zeta ,z\in {\cal A}_r, |\zeta |\ge 2|z|,
h\in B(0,1,{\cal A}_r)} |[f(z)(\zeta -z)^{-2}]^{\hat .}.h| |\zeta -z|^2
<\infty $, then $f$ is constant.}
\par {\bf Proof.} In view of Theorem $3.9$ there exists
a rectifiable path $\gamma $ in $U$ such that
$$(3.18)\quad (\partial ^kf(z)M/\partial z^k).h=
(2\pi )^{-1} (\int_{\gamma +z_0}f(\zeta ) (\partial ^k
(\zeta -z)^{-1}/\partial z^k).h d\zeta )$$
for each $h\in {\cal A}_r^k$, where $\gamma (t)=\rho ' \exp (2 \pi tM)$
with $t\in [0,1]$, $0<\rho '$. Then, in particular, for
$h=1^{\otimes k}$ and omitting it for short, we get
$$(3.19)\quad (\partial ^kf(z)M/\partial z^k)=
k! (2\pi )^{-1} (\int_{\gamma +z_0}f(\zeta ) 
(\zeta -z)^{-k-1} d\zeta ).$$
Therefore, $|f'(z)|\le C/d(z)$, since $ |(\partial (\zeta -z)^{-1}/
\partial z).s|=|\zeta -z|^{-2}$ for each $s\in \bf b$.
Since $\int_{\zeta }^zdf(z)=f(z)-f(\zeta )$,
then $|f(\xi )-f(z)|/|\xi -z|\le
\sup_{z\in B(a,\rho /2,{\cal A}_r)}[C/d(z)]\le 2C/\rho $, where $\rho '<
\rho /2$, $\xi $ and $z\in B(a,\rho /2,{\cal A}_r)\subset
Int (B(a,\rho ,{\cal A}_r))\subset U$.
Taking $\rho $ tending to infinity,
if $f$ is ${\cal A}_r$-holomorphic with bounded
$[f(z)(\zeta -z)^{-2}]^{\hat .}.h |\zeta -z|^2$
on ${\cal A}_r^2\times B(0,1,{\cal A}_r)$ for $|\zeta |\ge 2|z|$, then
$f'(z)=0$ for each $z\in {\cal A}_r$, since $f$ is locally $z$-analytic
and $\sup_{\zeta ,z\in U, |\zeta |\ge 2|z|,
h\in B(0,1,{\cal A}_r)} |[f(z)(\zeta -z)^{-2}]^{\hat .}.h|
|\zeta -z|^2 <\infty $ is bounded, hence $f$ is constant on ${\cal A}_r$.
\par {\bf 3.16. Remark.} Theorems 3.9, 3.10 and 3.15 are
the ${\cal A}_r$-analogs of the Cauchy, Morera and Liouville theorems
correspondingly. Evidently, Theorem 3.15 is also true for
right superlinear ${\hat f}(z)$ on ${\cal A}_r$ for each $z\in U$ and
with bounded ${\hat f}(z).h$ on $U\times B(0,1,{\cal A}_r)$ instead of
$[f(z)(\zeta -z)^{-2}]^{\hat .}.h |\zeta -z|^2$.
\par {\bf 3.17. Theorem.} {\it Let $P(z)$ be a polynomial
on ${\cal A}_r$, $2\le r\le \infty $,
such that $P(z)=z^{n+1}+\sum_{\eta (k)=0}^n(A_k,z^k)$,
where $A_k=(a_{1,k},...,a_{m,k}),$ $a_{j,l}\in {\cal A}_r$,
$k=(k_1,...,k_m)$, $0\le k_j\in \bf Z$, $\eta (k)=k_1+...+k_m$,
$0\le m=m(k)\in \bf Z$, $m(k)\le \eta (k)+1$,
$(A_k,z^k):=\{ a_{1,k}z^{k_1}...a_{m,k}z^{k_s} \} _{q(m+\eta (k))}$,
$z^0:=1$. Then $P(z)$ has a root in ${\cal A}_r$.}
\par {\bf Proof.} Consider at first $r<\infty $.
Suppose that $P(z)\ne 0$ for each $z\in {\cal A}_r$.
Consider a rectifiable path $\gamma _R$ in ${\cal A}_r$ such that
$\gamma _R ([0,1]) \cap {\cal A}_r = [-R,R]$ and outside $[-R,R]$:
$\quad \gamma _R(t) = R\exp (2\pi tM)$, where $M$ is a vector
in ${\cal I}_r$ with $|M|=1$, $0\le t\le 1/2$. Express $\tilde P$ through
variable $z$ also using $z^* = (2^r-2)^{-1} \{ -z
+ \sum_{s\in {\hat b}_r} s(zs^*) \} $.
Since $\lim_{|z|\to \infty } P(z) z^{-n-1}=1$, then due to Theorem 2.11
$\lim_{R\to \infty } \int_{\gamma _R}(P{\tilde P})^{-1}(z)dz=$
$\int_{-R}^R(P{\tilde P})^{-1}(v)dv$ $=\int_{-R}^R |P(v)|^{-2}dv\ge 0.$
But $\lim_{R\to \infty } \int_{\gamma _R} (P{\tilde P})^{-1}(z)dz=
\lim_{R\to \infty }\pi R^{-2n-1}=0$.
On the other hand, $\int_{-R}^R|P(v)|^{-2}dv=0$ if and only if
$|P(v)|^{-2}=0$ for each $v\in \bf R$.
This contradicts our supposition, hence
there exists a root $z_0\in {\cal A}_r$, that is, $P(z_0)=0$.
In the case $r=\infty $ use that $z=\lim_{r\to \infty }z_r$.
\par {\bf 3.18. Theorem.} {\it Let $f$ be an ${\cal A}_r$-holomorphic
function on an open subset $U$ in ${\cal A}_r$, $2\le r\le \infty $.
Suppose that $\epsilon >0$ and $\sf K$ is a compact subset of $U$.
Then for each $M\in {\cal I}_r$, $|M|=1$, there exists a function
$g_M(z)=P_{\infty }(z)+\sum_{k=1}^{\nu }P_k[(z-a_k)^{-1}]$,
$z\in {\cal A}_r\setminus \{ a_1,...,a_{\nu } \} $, $\nu \in \bf N$,
where $P_{\infty }$ and $P_j$ are polynomials, $a_j\in Fr (U)$,
$Fr (U)$ denotes a topological boundary of $U$ in ${\cal A}_r$,
such that $|f(z)M-g_M(z)|<\epsilon $ for each $z\in \sf K$.
If $r=2$ or $r=3$, then this statement is true for $M=1$ also.}
\par {\bf Proof.} Consider
cubes $S_j$ with ribs parallel to the basic axes and of length
$n^{-1}$ in ${\cal A}_r$ and putting $S := \cup_j S_j$ such that
${\sf K}\subset Int (S)$, where $n\in \bf N$ tends to infinity.
Since $f$ is ${\cal A}_r$-holomorphic and $\bf  K$ is compact,
we may apply Formula $(3.10)$ to each $\gamma \subset Fr (S_j)$
defined by directing vector $M\in {\cal I}_r$, $|M|=1$.
It can be seen, that $f$ can be approximated uniformly on $\sf K$
by a sum of the form $\sum_{k=1}^{\mu }\{ (a_{1,k}(\zeta _k-z)^{-1}
a_{2,k}) \} _{q(3)}$, where $a_{j,k}\in {\cal A}_r$, $\zeta _k\in Fr (S_j)$. 
For a given $n\in \bf N$ if $b\in Fr (S_j)$, then there exists
$a\in Fr (U_j)\cup \partial B(0,\rho ,{\cal A}_r)$ such that
$|b-a|\le n^{-1}$. If $z\in \sf K$ and $|z-a|\ge n^{-1}$, then
the series $(z-b)^{-1}=(\sum_{k=0}^{\infty }[(z-a)^{-1}(b-a)]^k)
(z-a)^{-1}$ converges uniformly on $\sf K$ and it is clear that
$fM$ can be approximated uniformly on $\sf K$ by a function of
the indicated form. In particular, for $r=2$ or $r=3$ the equation
$(ab)y=a$ has a solution for each $b\ne 0$ in ${\cal A}_r$,
that gives approximation of $f$ by $g_MM^*$.
\par {\bf 3.19. Note and Definitions.} Consider a one-point
(Alexandroff) compactification ${\hat {\cal A}}_r$ of the locally compact
topological space ${\cal A}_r$ for $2\le r\in \bf N$.
It is homeomorphic to a unit $2^r$-dimensional sphere $S^{2^r}$
in the Euclidean space $\bf R^{2^r+1}$.
In case $r=\infty $ consider a unit sphere $S^{\infty }$
with centre in $0$ in ${\bf R}\oplus {\cal A}_{\infty }$
such that ${\cal A}_{\infty }$
is topologically homeomorphic to $S^{\infty }\setminus \{ (1,0,0,...) \} $.
If $\zeta $ is a point in $S^{2^r}$ or $S^{\infty }$ different
from $(1,0,0,...)$, then the straight line containing
$(1,0,0,...)$ and $\zeta $ crosses $\pi _S$ in a finite point
$z$, where $\pi _S$ is the $2^r$-dimensional or $\infty $-dimensional
respectively plane orthogonal to the vector
$(1,0,0,...)$ and tangent to $S^{2^r}$ or $S^{\infty }$ at the south pole
$(-1,0,0,...)$. This defines the bijective continuous mapping from
$S^{2^r}\setminus \{ (1,0,0,...) \} $ or $S^{\infty }\setminus
\{ (1,0,0,...) \} $ onto $\pi _S$ such that $(1,0,0,...)$ corresponds to
the point of infinity. Therefore each function on a subset $U$
of ${\cal A}_r$ as a topological space can be considered on the
homeomorphic subset $V$ in $S^{2^r}$ or $S^{\infty }$.
The infinite dimensional sphere $S^{\infty }$ is not (locally)
compact relative to the norm topology, but it is compact relative to
the weak topology inherited from $l_2$ (see \cite{eng,nari}).
\par Let $z_0\in {\hat {\cal A}}_r$ be a marked point.
If a function $f$ is defined and ${\cal A}_r$-holomorphic 
on $V\setminus \{ z_0 \} $,
where $V$ is a neighbourhood of $z_0$, then $z_0$ is called
a point of an isolated singularity of $f$.
\par Suppose that $f$ is an ${\cal A}_r$-holomorphic function
in $B(a,0,\rho ,{\cal A}_r)\setminus \{ a \} $ for some $\rho >0$.
Then we say that $f$ has an isolated singularity at $a$.
Let $B(\infty ,\rho ,{\cal A}_r):=\{ z\in {\hat {\cal A}}_r$ such that
$\rho ^{-1}< |z| \le \infty \} $. Then we say that $f$ has an isolated
singularity at $\infty $ if it is ${\cal A}_r$-holomorphic
in some $B(\infty ,\rho ,{\cal A}_r)$.
\par Let $f: U\to {\cal A}_r$ be a function, where $U$ is a neighbourhood
of $z\in {\hat {\cal A}}_r$. Then $f$ is said to be meromorphic at
$z$ if $f$ has an isolated singularity at $z$.
If $U$ is an open subset in ${\hat {\cal A}}_r$, then $f$ is called
meromorphic in $U$ if $f$ is meromorphic at each point $z\in U$.
If $U$ is a domain of $f$ and $f$ is meromorphic in $U$, then
$f$ is called meromorphic on $U$. Denote by ${\bf M}(U)$ the set
of all meromorphic functions on $U$. Let $f$ be meromorphic on
a region $U$ in ${\hat {\cal A}}_r$. A point
$$c \in \bigcap_{V\subset U, V \mbox{ is closed and bounded }}
cl (f(U\setminus V)) $$
is called a cluster value of $f$.
\par {\bf 3.20. Proposition.} {\it Let $f$ be an
${\cal A}_r$-holomorphic function, $2\le r\le \infty $, with a right
${\cal A}_r$-superlinear superdifferential on
an open connected subset $U\subset {\hat {\cal A}}_r$ and
suppose that there exists a sequence of points $z_n\in U$
having a cluster point $z\in U$ such that $f(z_n)=0$
for each $n\in \bf N$, then $f=0$ everywhere on $U$.}
\par {\bf Proof} follows from the local $z$-analyticity of
$f$ and the fact $f^{(k)}(z)=0$ for each $0\le k\in \bf Z$
(see Definiton $2.2$, Theorems 2.11 and 3.10), when $f'(z)$ is right
${\cal A}_r$-superlinear on $U$. Therefore, $f$ is equal to
zero on a neighbourhood of $z$. The maximal subset of $U$
on which $f$ is equal to zero is open in $U$. On the other hand it is closed,
since $f$ is continuous, hence $f$ is equal to zero on $U$, since
$U$ is connected.
\par {\bf 3.21. Theorem.} {\it Let $\bf A$ denote the family of
all functions $f$ such that $f$ is ${\cal A}_r$-holomorphic
on $U:=Int (B(a,\rho ,R, {\cal A}_r))$, $3\le r\le \infty $,
where $a$ is a marked point in ${\cal A}_r$, $0\le \rho <R<\infty $
are fixed. Let $\bf S$ denote a subset
of ${\bf Z}^{\bf N}$ such that for each $k\in \bf S$ there exists
$m(k):=\max \{ j: $ $k_j\ne 0, k_i=0 $ $\mbox{for each }
i>j \} \in \bf N$ and let ${\bf B}$ be a family of finite sequences
$b_k=(b_{k,1},...,b_{k,m(k)}; q(m(k)+\eta (k)))$
such that $b_{k,j}\in {\cal A}_r$ for each $j=1,...,n$,
$n\in \bf N$, $q(m(k)+\eta (k))$ as in \S 2.1.
Then there exists a bijective correspondence between
$\bf A$ and $\zeta \in {\bf B^S}$ such that
$$(3.20) \quad \lim_{m+\eta
\to \infty }\sup_{z\in B(a,\rho _1,R_1,{\cal A}_r)}
\sum_{k, m(k)=m, \eta (k)=\eta } |\{ (b_k,z^k) \} _{q(m(k)+\eta (k))}|=0$$
for each $\rho _1$ and $R_1$ such that $\rho <\rho _1<R_1<R$,
where $\eta (k):=k_1+...+k_{m(k)}$, $\zeta (k)=:b_k
=(b_{k,1},...,b_{k,m(k)};q(m(k)+ \eta (k)))$,
$\{ (b_k,z^k) \} _{q(m(k)+\eta (k))}=
\{ b_{k,1}z^{k_1}...b_{k,m(k)}z^{k_{m(k)}} \} _{q(m(k)+\eta (k))}$
for each $k\in \bf S$, that is,
$f\in \bf A$ can be presented by a convergent series
$$(3.21)\quad f(z)=\sum_{b\in \zeta }\{ (b_k,z^k) \} _{q(m(k)+\eta (k))}.$$}
\par {\bf Proof.} If Condition $(3.20)$ is satisfied, then
the series $(3.21)$ converges on $B(a,\rho ',R',{\cal A}_r)$ for each
$\rho '$ and $R'$ such that $\rho <\rho '<R'<R$, since $\rho _1$
and $R_1$ are arbitrary such that $\rho <\rho _1<R_1<R$ and
$\sum_{n=0}^{\infty }p^n$
converges for each $|p|<1$. In particular taking
$\rho _1<\rho '<R'<R_1$ for $p=R'/R_1$ or
$p=\rho _1/\rho '$. Therefore, from $(3.20)$ and
$(3.21)$ it follows, that $f$ presented by the series $(3.21)$
is ${\cal A}_r$-holomorphic on $U$.
\par Vice versa let $f$ be in $\bf A$. In view of Theorems 2.11 and 3.9
there are two rectifiable closed paths $\gamma _1$ and $\gamma _2$
such that $\gamma _2(t)=a+\rho ' \exp (2\pi t M_2)$ and
$\gamma _1(t)=a+R' \exp (2\pi t M_1)$, where $t\in [0,1]$,
$M_1$ and $M_2\in {\cal A}_r$ with $|M_1|=1$ and $|M_2|=1$,
where $\rho <\rho '<R'<R$, because as in \S 3.9 $U$ can be presented
as a finite union of regions $U_j$ each of which satisfies
the conditions of Theorem 2.11. Using a finite number
of rectifiable paths $w_j$ (joining $\gamma _1$ and $\gamma _2$
within $U_j$) going twice in one and the
opposite directions leads to the conclusion
that for each $z\in Int (B(a,\rho ',R',{\cal A}_r))$
the function $f(z)M$ with $M=M_1=M_2$ is presented by the integral formula:
$$(3.22)\quad f(z)M=(2\pi )^{-1}\{
(\int_{\gamma _1}f(\zeta )(\zeta -z)^{-1} d\zeta ) -
(\int_{\gamma _2}f(\zeta )(\zeta -z)^{-1} d\zeta ) \} .$$
On $\gamma _1$ we have the inequality:
$|(\zeta -a)^{-1}(z-a)|<1$, on $\gamma _2$ another
inequality holds: $|(\zeta -a)(z-a)^{-1}|<1$.
In view of \S 3.6.2 as in \S 3.10 considering
different possible embeddings $\Upsilon _{\zeta -a,z-a}\subset \bf K$
into ${\cal A}_r$ and using alternative property of the octonion algebra
$\bf K$ we get, that for $\gamma _1$ the series
$$(\zeta -z)^{-1}=(\sum_{k=0}^{\infty }((\zeta -a)^{-1}(z-a))^k)
(\zeta -a)^{-1}$$
converges uniformly by $\zeta \in B(a,R_2+\epsilon,R_1,{\cal A}_r)$ and
$z\in B(a,\rho _2,R_2,{\cal A}_r)$ for $\zeta $ and $z$ such that
$\Upsilon _{\zeta -a,z-a}\hookrightarrow \bf K$,
while for $\gamma _2$ the series
$$(\zeta -z)^{-1}=-(z-a)^{-1}(\sum_{k=0}^{\infty }((\zeta -a)(z-a)^{-1})^k)$$
converges uniformly by $\zeta \in B(a,\rho _1,\rho _2-\epsilon,{\cal A}_r)$
and $z\in B(a,\rho _2,R_2,{\cal A}_r)$
for $\zeta $ and $z$ such that
$\Upsilon _{\zeta -a,z-a}\hookrightarrow \bf K$,
for each $\rho '<\rho _2<R_2<R'$
and each $0<\epsilon <\min (\rho _2-\rho _1,R_1-R_2)$,
since expansion coefficients by $((\zeta -a)^{-1} (z-a))$
in the first series and by $((\zeta -a)(z-a)^{-1})$ in the second series
are independent from the type of embedding.
Consider corresponding $\gamma _1$ and $\gamma _2$ such that
$\zeta \in \gamma _1$ or $\zeta \in \gamma _2$ respectively and
$a, \zeta , z$ are subjected to the condition $\Upsilon _{\zeta -a,
z-a}\hookrightarrow \bf K$, which is not finally restrictive
due to Theorem 2.15. Consequently,
$$(3.23)\quad f(z)M=\sum_{k=0}^{\infty }(\phi _k(z)+\psi _k(z)),
\mbox{ where}$$
$$\phi _k(z):=(2\pi )^{-1} \{ \int_{\gamma _1}
f(\zeta )(((\zeta -a)^{-1}(z-a))^k(\zeta -a)^{-1}) d\zeta ) \} ,$$
$$\psi _k(z):=(2\pi )^{-1}\{ \int_{\gamma _2}
f(\zeta )((z-a)^{-1}((\zeta -a)(z-a)^{-1})^k) d\zeta ) \} ,$$
and where $\phi _k(z)$ and $\psi _k(z)$ are ${\cal A}_r$-holomorphic
functions, hence $fM$ has decomposition $(3.21)$ in $U$,
since due to \S 2.15
and \S 3.9 there exists $\delta >0$ such that integrals for $\phi _k$
and $\psi _k$ by $\gamma _1$ and $\gamma _2$ are the same for each
$\rho '\in (\rho _1,\rho _1+\delta )$, $R'\in (R_1-\delta ,R_1)$.
Using the definition of the ${\cal A}_r$ line integral we get $(3.21)$
converging on $U$. Varying $z\in U$ by $|z|$ and $Arg (z)$ we get
that $(3.21)$ converges absolutely on $U$, consequently,
$(3.20)$ is satisfied for $fM$. Since $M\in {\cal A}_r$, $|M|=1$,
is arbitrary, then as in the proof of Theorem $3.10$ we get
the statement of this theorem for $f$.
\par {\bf 3.22. Notes and Definitions.} Let $\gamma $ be a closed
curve in ${\cal A}_r$. There are natural projections from ${\cal A}_r$
on complex planes: $\pi _s(z)=w_1+w_ss$ for each $s\in {\hat b}_r$,
where $2\le r\le \infty $, $z=\sum_{s\in {\bf b}_r}w_ss$ with real $w_s$
for each $s\in {\bf b}_r$. Therefore, $\pi _s(\gamma )=:\gamma _s$ are curves
in complex planes $\bf C_s$ isomorphic to ${\bf R}\oplus {\bf R}s$.
A curve $\gamma $ in ${\cal A}_r$ is closed
(a loop, in another words) if and only if $\gamma _s$ is
closed for each $s\in {\hat b}_r$, that is, $\gamma (0)=\gamma (1)$
and $\gamma _s(0)=\gamma _s(1)$ correspondingly.
In each complex plane there is the standard complex notion of
a topological index $In (a_s,\gamma _s)$ of a curve $\gamma _s$ at $a_s=
\pi _s(a)$.
Therefore, there exists a vector $In (a,\gamma ):=\{ In (a_s,\gamma _s):
s\in {\hat b}_r \} $ which we call
the topological index of $\gamma $ at a point $a\in {\cal A}_r$.
This topological index is invariant relative to homotopies satisfying
conditions of Theorem 3.9.
Consider now a standard closed curve $\gamma (s)=a+\rho \exp (2\pi tnM)$,
where $M\in {\cal I}_r$ with $|M|=1$, $n\in \bf Z$, $\rho >0$, $t\in [0,1]$.
Then ${\hat I}n (a,\gamma ):=(2\pi )^{-1}(\int_{\gamma }d Ln (z-a))=nM$
is called the ${\cal A}_r$-index of $\gamma $ at a point $a$.
It is also invariant relative to homotopies satisfying the
conditions of Theorem
3.9. Moreover, ${\hat I}n (h_1(ah_2),h_1(\gamma h_2))={\hat I}n (a,\gamma )$
for each $h_1$ and $h_2\in {\cal A}_r\setminus \{ 0 \} $
such that $h_1(Mh_2)=M$. For $M=\sum_{s\in {\hat b}_r} m_ss$
there is the equality ${\hat I}n (a,\gamma )=
\sum_{s\in {\hat b}_r}In (a_s,\gamma _s)m_ss$
(adopting the corresponding convention for signs of indexes
in each $\bf C_s$ and the convention of positive directions of going along
curves).
In view of the properties of $Ln$ for each curve $\psi $ in ${\cal A}_r$
there exists $\int_{\gamma }d Ln (z-a)=2\pi qM$ for some
$q\in \bf R$ and $M\in {\cal I}_r$ with $|M|=1$.
For a closed  curve $\psi $ up to a composition of homotopies
each of which is charaterized by homotopies in $\bf C_s$
for $s\in {\hat b}_r$ there exists a standard $\gamma $ with a generator $M$
for which ${\hat I}n (a,\gamma )=qM$, where $q\in \bf Z$.
Therefore, we can take as a definiton
${\hat I}n (a,\psi )={\hat I}n (a,\gamma )$.
Define also the residue of a meromorphic function with an isolated
singularity at a point $a\in {\cal A}_r$ as \\
$(i)\quad res (a,f)M:=(2\pi )^{-1}(\int_{\gamma }f(z)dz)$,
where $\gamma (t)=a+\rho \exp (2\pi tM)\subset V,$ $\rho >0$,
$|M|=1$, $M\in {\cal I}_r$, $t\in [0,1]$,
$f$ is ${\cal A}_r$-holomorphic on $V\setminus \{ a \} $.
Extend $res (a,f)M$ by Formula $(i)$ on ${\cal I}_r$ as \\
$(ii)\quad res (a,f)M:=[res (a,f)(M/|M|)] |M|$, $\forall M\ne 0$;
$res (a,f)0:=0$, \\
when $res (a,f)M$ is finite for each $M\in {\cal I}_r$, $|M|=1$.
For $r=2$ or $r=3$ the equation for $res (a,f)$ can be resolved
for each $M\ne 0$.
\par If $f$ has an isolated singularity at $a\in {\hat {\cal A}}_r$,
then coefficients $b_k$ of its Laurent series (see \S 3.21)
are independent of $\rho >0$. The common series is called
the $a$-Laurent series. If $a=\infty $, then $g(z):=f(z^{-1})$
has a $0$-Laurent series $c_k$ such that $c_{-k}=b_k$.
Let $\beta := \sup_{b_k\ne 0}
\eta (k)$, where $\eta (k)=k_1+...+k_m$, $m=m(k)$
for $a=\infty $; $\beta =\inf_{b_k\ne 0} \eta (k)$ for $a\ne \infty $.
We say that $f$ has a removable singularity, pole,
essential singularity at $\infty $ according as
$\beta \le 0$, $0<\beta <\infty $, $\beta =+\infty $.
In the second case $\beta $ is called the order of the pole
at $\infty $. For a finite $a$ the corresponding cases
are: $\beta \ge 0$, $-\infty <\beta <0$, $\beta =-\infty $.
If $f$ has a pole at $a$, then $|\beta |$ is called
the order of the pole at $a$.
\par A value of a function $\partial _f(a):=\inf \{ \eta (k): b_k\ne 0 \} $
is called a divisor of $f$ at $a\ne \infty $,
$\partial _f(a):=\inf \{ - \eta (k): b_k\ne 0 \} $ for $a=\infty $,
where $b_k\ne 0$ means that $b_{k,1}\ne 0$,....,$b_{k,m(k)}\ne 0$.
Then $\partial _{f+g}(a)\ge \min \{ \partial _f(a), \partial _g(a) \} $
for each $a\in dom (f)\cap dom (g)$
and $\partial _{fg}(a)=\partial _f(a)+\partial _g(a).$
For a function $f$ meromorphic
on an open subset $U$ in ${\hat {\cal A}}_r$ the function $\partial _f(p)$
by the variable $p\in U$ is called the divisor of $f$.
\par {\bf 3.23. Theorem.} {\it Let $U$ be an open region in
${\hat {\cal A}}_r$, $2\le r\le \infty $, with $n$ distinct marked
points $p_1,...,p_n$, and let $f$ be an ${\cal A}_r$-holomorphic
function on $U\setminus \{ p_1,...,p_n \} =:U_0$
and $\psi $ be a rectifiable closed curve lying in
$U_0$ such that $U_0$ satisfies the conditions of Theorem 3.9
for each $z_0\in \{ p_1,...,p_n \} $. Then
$$\int_{\psi }f(z)dz= 2 \pi \sum_{j=1}^n res (p_j,f)
{\hat I}n (p_j,\psi )$$
and $res (p_j,f)M$ is the $\bf R$-homogeneous ${\cal I}_r$-additive
(of the variable $M$ in ${\cal I}_r$)
${\cal A}_r$-valued functional for each $j$.}
\par {\bf Proof.} For each $p_j$ consider the principal part $T_j$
of a Laurent series for $f$ in a neighbourhood of $p_j$, that is, \\
$T_j(z)=\sum_{k, \eta (k)<0} \{ (b_k,(z-p_j)^k) \} _{q(m(k)+ \eta (k))}$,
where $\eta (k)=k_1+...+k_n$ for $k=(k_1,...,k_n)$
(see Theorem 3.21).
Therefore, $h(z):=f(z)-\sum_jT_j(z)$ is a function having
an ${\cal A}_r$-holomorphic extension on $U$.
In view of Theorem 3.9 for an ${\cal A}_r$-holomorphic
function $g$ in a neighbourhood $V$ of a point $p$ and a rectifiable
closed curve $\zeta $ we have
$$g(p) {\hat I}n (p,\zeta )=(2\pi )^{-1}(\int_{\zeta }g(z)(z-p)^{-1}dz)$$
(see \S 3.22).
We may consider small $\zeta _j$ around each $p_j$ with
${\hat I}n (p_j,\zeta _j)={\hat I}n (p_j,\gamma )$
for each $j=1,...,n$.
Then $\int_{\zeta _j}f(z)dz=\int_{\zeta _j}T_j(z)dz$ for each $j$.
Representing $U_0$ as a finite union of open regions $U_j$
and joining $\zeta _j$ with $\gamma $ by paths $\omega _j$ going
in one and the opposite direction as in Theorem 3.9 we get
$$\int_{\gamma }f(z)dz+\sum_j\int_{\zeta _j^-}f(z)dz=0,$$
consequently,
$$\int_{\gamma }f(z)dz=\sum_j\int_{\zeta _j}f(z)dz=
\sum_j 2\pi res (p_j,f) {\hat I}n (p_j,\gamma ) ,$$
where ${\hat I}n (p_j,\gamma )$ and $res (p_j,f)$ are
invariant relative to homotopies satisfying conditions of Theorem 3.9.
Since $\int_{\zeta _j} g(z)d Ln (z-p_j)$ is $\bf R$-homogeneous
and ${\cal I}_r$-additive relative to a directing vector
$M\in {\cal I}_r$ of a loop $\zeta _j$, then
$res (p_j,f)M$ defined by Formulas
$2.22.(i,ii)$ is $\bf R$-homogeneous ${\cal I}_r$-additive
of the argument $M$ in ${\cal I}_r$.
\par {\bf 3.24. Corollary.} {\it Let $f$ and $T$
be the same as in \S 3.23, then $res (p_j,f)M=res (p_j,T_j)M=
res (p_j,\sum_{k,\eta (k)=-1} \{ (b_k,(z-p_j)^k) \} _{q(m(k)+\eta (k))} M$,
in particular, $res (p_j, \{ b (z-p_j)^{-1}c
\} _{q(3)} ) M= \{ bMc \} _{q(3)} $ for each $b, c \in {\cal A}_r$.}
\par {\bf Proof.} The first statement follows from \S 3.23,
the second statement follows from left and right-${\cal A}_r$-linearity
of the line integral, though it is not the superlinear functional
(see Theorem $2.7$).
\par {\bf 3.25. Corollary.} {\it Let $U$ be an open region
in ${\hat {\cal A}}_r$, $2\le r\le \infty $,
with $n$ distinct points $p_1,...,p_n$,
let also $f$ be an ${\cal A}_r$-holomorphic function
on $U\setminus \{ p_1,...,p_n \} =:U_0$, $p_n=\infty $,
and $U_0$ satisfies conditions of Theorem 3.9 with at least
one $\psi $, $\gamma $ and each $z_0\in \{ p_1,...,p_n \} $.
Then $\sum_{p_j\in U} res (p_j,f)M =0$.}
\par {\bf Proof.} If $\gamma $ is a closed curve encompassing
$p_1$,...,$p_{n-1}$, then $\gamma ^-(t):=\gamma (1-t)$, where $t\in [0,1]$,
encompasses $p_n=\infty $ with positive going by $\gamma ^-$ relative to
$p_n$. Since $\int_{\gamma }f(z)dz+\int_{\gamma ^-}f(z)dz=0$,
we get from Theorem 3.23, that $\sum_{p_j\in U} res (p_j,f)M=0$
for each $M\in {\cal I}_r$, hence
$\sum_{p_j\in U} res (p_j,f)M =0$ is the zero $\bf R$-homogeoneous
${\cal I}_r$-additive ${\cal A}_r$-valued functional on ${\cal I}_r$.
\par {\bf 3.26. Definitions.} Let $f$ be an ${\cal A}_r$-holomorphic
function, $2\le r\le \infty $, on a neighbourhood $V$
of a point $z\in {\cal A}_r$.
Then the infimum: $\eta (z;f):=\inf \{ k: k\in {\bf N},
f^{(k)}(z)\ne 0 \} $ is called a multiplicity of $f$ at $z$.
Let $f$ be an ${\cal A}_r$-holomorphic function on an open subset $U$
in ${\hat {\cal A}}_r$, $2\le r\le \infty $. Suppose
$w\in {\hat {\cal A}}_r$, then the valence $\nu _f(w)$ of $f$
at $w$ is by the definition $\nu _f(w):=\infty $, when the set
$\{ z: f(z)=w \} $ is infinite, and otherwise
$\nu _f(w):=\sum_{z, f(z)=w}\eta (z;f)$.
\par {\bf 3.26.1. Theorem.} {\it Let $f$ be an ${\cal A}_r$-meromorphic
right superlinearly superdifferentiable
function on a region $U\subset {\hat {\cal A}}_r$. If $b\in {\hat
{\cal A}}_r$ and $\nu _f(b)<\infty $, then $b$ is not a cluster value of
$f$ and the set $ \{ z: \nu _f(z)=\nu _f(b) \} $ is a neighbourhood
of $b$. If $U\ne {\hat {\cal A}}_r$ or $f$ is not constant, then
the converse statement holds. Nevertheless, it is false, when
$f=const $ on ${\hat {\cal A}}_r$.}
\par {\bf 3.26.2. Theorem.} {\it Let $U$ be a proper open subset
of ${\hat {\cal A}}_r$, $2\le r\le \infty $, let also $f$ and $g$
be two continuous functions
from ${\bar U}:=cl (U)$ into ${\hat {\cal A}}_r$ such that on a
topological boundary $Fr (U)$ of $U$ they satisfy the inequality
$|f(z)|<|g(z)|$ for each $z\in Fr (U)$. Suppose $f$ and $g$ are
${\cal A}_r$-meromorphic functions in $U$ and $h$ be a unique
continuous map from $\bar U$ into ${\hat {\cal A}}_r$ such that
$h|_E=f|_E +g|_E$, where $E:= \{ z: f(z)\ne \infty , g(z)\ne \infty \} $,
$[\partial Ln (h(z))/\partial z]$ is right superlinear in $U_{z_0}$
for each zero $z_0$ and in $U_{z_0}\setminus \{ z_0 \} $ for each pole
$z_0$, where $U_{z_0}$ is a neighbourhood of $z_0$, $z\in U_{z_0}$ or
$z\in U_{z_0}\setminus \{ z_0 \} $ respectively. Then
$\nu _{g|_U}(0)-\nu _{g|_U}(\infty )=\nu _{h|_U}(0)-\nu _{h|_U}(\infty )$.}
\par {\bf Proofs} of these two theorems are analogous to that of
Theorems VI.4.1, 4.2 \cite{heins}. To prove Theorem 3.26.2 consider
the function $\zeta (z,t) := tf(z)+g(z)$ for each $z\in {\bar U}$
and each $t\in [0,1]\subset \bf R$. If $z_0$ is a pole of $h(z)$, then
$z_0$ is a zero of $1/h(z)$. By the supposition of Theorem 3.26.1,
Proposition 2.3 and Theorem 3.9 $res (z_0, Ln (h))$ is the right
superlinear operator for each zero or pole $z_0$.
On the other hand there exists $\delta >0$ such that
$\Delta _{\gamma } Arg (1 + t g^{-1}f) =0$ for each $t\in [-1,1]$, when
no any pole or zero of $g$ or $f$ belongs to a rectifiable loop
$\gamma $ in $U$ with $dist (\gamma , Fr (U))< \delta $,
where $dist (A,B):=\sup_{z\in A}
(\inf_{\xi \in B} |\xi -z|) + (\sup_{\xi \in B} \inf_{z\in A} |\xi -z|)$.
Then $\int_{z\in \gamma } d Ln
\zeta (z,t)$ is the continuous function by $t\in [0,1]$ taking
values $2 \pi nM$, where $M\in {\cal I}_r$ charterizes a rectifiable loop
$\gamma $ contained in $U$, $|M|=1$, $n\in \bf Z$,
$M$ is independent of $t$. For each $\delta >0$ it is possible to choose
a rectifiable loop $\gamma $ in $U$ such that
$dist (\gamma , Fr (U))<\delta $. Then apply Theorem 3.23
to suitable pieces of $U$ whose boundaries do not contain zeros
and poles of $f$ and $g$.
\par From the proof of Theorem 3.26.2 we get.
\par {\bf 3.26.3. Corollary.} {\it Let suppositions of Theorem 3.26
be satisfied may be besides the condition of right supelinearity
of $[\partial Ln (h(z))/\partial z]$, then $\Delta _{\partial U}
Arg (f) = \Delta _{\partial U} Arg (g)=\int_{\gamma }d Ln (f(z))$,
where $\gamma $ is as in \S 3.26.2.}
\par {\bf 3.27. Theorem.} {\it Let $U$ be an open subset
in ${\cal A}_r^n$, $2\le r\le \infty $, then there exists a represenation
of the $\bf R$-linear space
$C^{\omega }_{z,\tilde z}(U,{\cal A}_r)$ of locally
$(z,\tilde z)$-analytic functions on $U$ such that it is isomorphic to the
$\bf R$-linear space $C^{\omega }_z(U,{\cal A}_r)$
of ${\cal A}_r$-holomorphic functions on $U$.}
\par {\bf Proof.} Evidently, the proof can be reduced to the case
$n=1$ by induction considering local $(z,\tilde z)$-series
decompositions by $(\mbox{ }^nz,\mbox{ }^n{\tilde z})$
with coefficients being convergent series of $(\mbox{ }^1z,
\mbox{ }^1{\tilde z},...,\mbox{ }^{n-1}z,\mbox{ }^{n-1}{\tilde z})$.
We have \\
${\tilde z}=(2^r-2)^{-1} \{ -z +\sum_{s\in {\hat b}_r} s(zs^*) \} $ 
for each $2\le r<\infty $, \\
${\tilde z} = \lim_{r\to \infty } (2^r-2)^{-1}
\{ -z + \sum_{s\in {\hat b}_r} s(zs^*) \} $ in ${\cal A}_{\infty }$. \\
Consequently, each polynomial in $(z,\tilde z)$ is also a polynomial
in $z$ only, moreover, each polynomial locally $(z,\tilde z)$
analytic function on $U$ is polynomial locally $z$-analytic on $U$.
Then if a series by $(z,{\tilde z})$ converges in a ball
$B(z_0,\rho ,{\cal A}_r^n)$, then its series in the $z$-representation
converges in a ball $B(z_0,\rho ',{\cal A}_r^n)$, where
$\rho '=2^r (2^r-2)^{-1} \rho $. Considering basic polynomials
of any polynomial basis in $C^{\omega }_{z,\tilde z}(U,{\cal A}_r)$
we get (due to infinite dimensionality of this space) a polynomial base of
$C^{\omega }_z(U,{\cal A}_r)$. This establishes the 
$\bf R$-linear isomorphism between these two spaces.
Moreover, in such representation of the space
$C^{\omega }_{z,\tilde z}(U,{\cal A}_r)$ we can put $D_{\tilde z}=0$,
yielding for differential forms $\partial _{\tilde z}=0$,
this leads to differential calculus and integration with
respect to $D_z$ and $dz$ only.
\par {\bf 3.28. Theorem} (Argument principle). {\it
Let $f$ be an ${\cal A}_r$-holomorphic function on an open region
$U$ satisfying conditions of \S 3.9, $2\le r\le \infty $,
and let $\gamma $ be a closed curve contained in $U$,
where $[\partial Ln (f(z))/\partial z]$ is right supelinear
in some neighbourhood $U_{z_0}$ for each zero $z_0$ of $f(z)$. Then
${\hat I}n (0; f\circ \gamma )=\sum_{\partial _f(a)\ne 0}
{\hat I}n (a; \gamma )\partial _f(a)$.}
\par {\bf Proof.} There is the equality ${\hat I}n(0;f\circ \gamma )
=\int_{\zeta \in \gamma }d Ln (f(\zeta ))=$
$\int_0^1d Ln (f\circ \gamma (s))=$ $\int_{\gamma }f^{-1}(\zeta )df(\zeta )$.
Let $\partial _f(a)=n\in \bf N$, then
$$f^{-1}(a)f'(a).s=
\sum_{l,k; n_1+...+n_k=\partial _f(a), 0\le n_j\in {\bf Z}, j=1,...,k}
\{ (z-a)^{n_1} g_{s,l,k,1;n_1,...,n_k}(z)$$
$$(z-a)^{n_2}g_{s,l,k,2;n_1,...,n_k}(z)...
(z-a)^{n_k} g_{s,l,k,k;n_1,...,n_k}(z) \} _{q(n_1+...+n_k+k)} ,$$
where $g_{s,l,p,k;n_1,...,n_k}(z)$ are
${\cal A}_r$-holomorphic functions of $z$ on $U$ for each $s\in {\bf b}_r$
such that $g_{s,l,p,k;n_1,...,n_k}(a)\ne 0$, 
where $l=1,...,m$, $1\le m\le 2^{r \partial _f(a)}$
for finite $r$ and each $m\in \bf N$ for $r=\infty $ (see \S \S 2.8,
3.7, 3.21, 3.27), since each term
$\xi (z) \prod_{s\in {\bf b}_r} (w_s-w_{s,0})^{n_s}$
with $\sum_{s\in {\bf b}_r}n_s\ge \partial _f(a)$, $n_j\ge 0$,
has such decomposition, where $2\le r< \infty $,
$\xi (z)$ is an ${\cal A}_r$-holomorphic function on a neighbourhood
of $a$ such that $\xi (a)\ne 0$. When $r=\infty $ use
$z = \lim_{r\to \infty } z_r$.
Suppose $\psi $ is a closed curve such that ${\hat I}n(p,\psi )=
2\pi nM$, $|M|=1$, $M\in {\cal I}_r$, $0\ne n\in \bf Z$.
Then we can define a curve $\psi ^{1/n}=:\omega $ as a closed curve
for which ${\hat I}n(p,\omega )=2\pi M$ and $\omega ([0,1])\subset
\psi ([0,1])$. Then we call $\omega ^n=\psi $. That is,
${\hat I}n(p,\psi ^{1/n})={\hat I}(p,\psi )/n$. The latter formula
allows an interpretation also when ${\hat I}n (p,\psi )/n$ is
equal to $2\pi qM$, where $0\ne q\in \bf Q$. That is, a curve
$\psi ^{1/n}$ can be defined for each $0\ne n\in \bf Z$.
This means that $\gamma $ can be presented
as union of curves $\omega _j$ for each of which there exists
$n_j\in \bf N$ such that $\omega _j^{n_j}$ is a closed curve.
Using Theorem 3.9 for each $a\in U$ with $\partial _f(a)\ne 0$, also
using the series given above we can find a finite family of
$\omega _j$ for which one of the terms in the series is not less,
than any other term. We may also use small homotopic
deformations of $\omega _j$ satisfying the conditions of Theorem 3.9
such that in the series one of the terms is greater
than any other for almost all points on $\omega _j$.
Such deformation is permitted, since otherwise two terms would coincide
on an open subset of $U$, that is impossible.
Considering such series, Formulas $2.5.(4,5)$ and using Theorem 3.26.2
we get the statement of this theorem.
\par {\bf 3.29. Theorem.} {\it If $f$ has an essential
singularity at $a$, then $cl (f(V))={\hat {\cal A}}_r$
for each $V\subset dom (f)$, $V=U\setminus \{ a \} $,
where $U$ is a neighbourhood of $a$.}
\par {\bf Proof.} Suppose that the statement of this theorem
is false, then there would exist $\rho >0$ and $m>0$ and an
element $A\in {\cal A}_r$ such that $f$ is $z$-analytic in
$B(a,0,\rho ,{\cal A}_r)\setminus
\{ a \} $ and $|f(z)-A|\ge m$ for each $z$ such that
$0<|z-a|<\rho $. If $\infty \notin cl (f(V))$, then there exists
$R>0$ such that $A\notin cl (f(V))$ for each $|A|>R$.
Therefore, the function $[f(z)-A]^{-1}$ is ${\cal A}_r$-holomorphic
in $B(a,0,\rho ,{\cal A}_r)\setminus \{ a \} $.
Hence $[f(z)-A]^{-1}=\sum_k \{ (p_k,(z-a)^k) \} _{q(\eta (k)+m(k))}$, where
in this sum $k=(k_1,...,k_{m(k)})$ with $k_j\ge 0$
for each $j=1,...,m(k)\in \bf N$, $p_k$ are finite sequences
of coefficients for $[f(z)-A]^{-1}$ as in \S 3.21.
If $D_z^n([f(z)-A]^{-1})|_{z=a}=0$ for each $n\ge 0$, then
$[f(z)-A]^{-1}=0$ in a neighbourhood of $a$.
Therefore, $[f(z)-A]^{-1}=\sum_{n_1+...+n_l=n} \{ g_1z^{n_1}...
g_lz^{n_l} \} _{q(\eta (n)+l)} $ for some $n$ such that $0\le n\in \bf N$,
$n_j\ge 0$ for each $j=1,...,l\in \bf N$,
each $g_j$ is an ${\cal A}_r$-holomorphic function (of $z$).
Consequently, taking inverses of both sides
$[f(z)-A]$ and $(\sum_{n_1+...+n_l=n} \{ g_1z^{n_1}...
g_lz^{n_l} \} _{q(\eta (n)+l)} )^{-1}$ and comparing their expansion
series we see that finite sequences $b_k$ of expansion
coefficients for $f$ have the property $b_k=0$ for each
$\eta (k)<-n$. This contradicts the hypothesis
and proves the theorem.
\par {\bf 3.30. Definition.} Let $a$ and $b$ be two points in ${\cal A}_r$
and $\theta $ be a stereographic mapping of the unit 
real sphere $S^{2^r}$ for $2\le r<\infty $ or $S^{\infty }$ for $r=\infty $
on ${\hat {\cal A}}_r$. Then $\chi (a,b):=
|\phi (a)-\phi (b)|_Y$ is called the chordal metric,
where $\phi :=\theta ^{-1}: {\hat {\cal A}}_r\to S^m$,
$S^m$ is embedded in $Y:=\bf R^{m+1}$ for $m:=2^r$ with $r<\infty $
or in ${\bf R}\oplus l_2({\bf R})$ for $m=\infty $ with
$r=\infty $, $|*|_Y$ is the Euclidean or Hilbert norm
in $Y$ respectively.
\par {\bf 3.30.1. Theorem.} {\it Let $U$ be an open
region in ${\hat {\cal A}}_r$, $\{ f_n: n \in {\bf N} \} $ be a sequence
of functions meromorphic on $U$ tending uniformly in $U$
to $f$ relative to the chordal metric. Then either $f$ is the constant
$\infty $ or else $f$ is meromorphic on $U$.}
\par {\bf 3.30.2. Theorem.} {\it Let $\{ f_k: k\in {\bf N} \} $
be a sequence of meromorphic functions on an open subset $U$ in
${\hat {\cal A}}_r$,
which tends uniformly in the sence of the chordal metric in $U$
to $f$, $f\ne const $. If $f(a)=b$ and $\rho >0$ are such that
$B(a,\rho ,{\cal A}_r)\subset U$ and $f(z)\ne b$ for each
$z\in B(a,\rho ,{\cal A}_r)
\setminus \{ a \} $, then there exists $m\in \bf N$ such that
the value of the valence of $f_k|_{B(a,\rho ,{\cal A}_r)}$ at $b$
is $\eta (b;f)=\eta (a;f)$ for each $k\ge m$.}
\par {\bf 3.30.3. Note.} The proofs of these theorems are
formally similar to the proofs of VI.4.3 and 4.4 \cite{heins}.
Theorems 3.26.2 and 3.30.2 are the ${\cal A}_r$ analogs of the
Rouch\'e and Hurwitz theorems respectively.
There are also the following ${\cal A}_r$ analogs of
the Mittag-Leffler and Weierstrass theorems.
Their proofs are similar to those for Theorems VIII.1.1 and
1.2 \cite{heins} respectively. Nevertheless the second part
of the Weierstrass theorem is not true in general because of
noncommutativity of ${\cal A}_r$, that is, a function $h\in {\bf M}(U)$
with $\partial _h=\partial $ is not necessarily
representable as $h=fg$, where $g$ is an ${\cal A}_r$-holomorphic on $U$
and $f$ is another marked function $f\in {\bf M}(U)$ such that
$\partial _f=\partial $.
\par {\bf 3.31. Theorem.} {\it Let $U$ be a nonempty proper open
subset of ${\hat {\cal A}}_r$, $2\le r\le \infty $,
let $A\subset U$ not containing
any cluster point in $U$. Let there be a function
$g_b\in {\bf M}({\hat {\cal A}}_r)$
for each $b\in A$ having a pole at $b$ and no other.
Then there exists $f\in {\bf M}(U)$ ${\cal A}_r$-holomorphic on
$U\setminus B$ and having the same principal part at $b$ as $g_b$.
If $f$ is such a function, then each other such function
is the function $f+g$, where $g$ is ${\cal A}_r$-holomorphic on $U$.}
\par {\bf 3.32. Theorem.} {\it Let $U$ be a proper nonempty
open subset of ${\hat {\cal A}}_r$, $2\le r\le \infty $.
Let $\partial : U\to \bf Z$ be a
function such that $ \{ \partial (z)\ne 0 \} $ does not have
a cluster point in $U$.
Then there exists $f\in {\bf M}(U)$ such that $\partial _f=\partial $.}
\par {\bf Proof.} If $a_j$ is a zero, that is, $\partial (a_j)\ne 0$,
then take a circle of radius $\delta _j>0$ with centre at $a_j$.
There are possible two cases: $|a_j| \delta _j\ge 1$ and
$|a_j| \delta _j<1$. At $a_j$ of the first type construct in $U$ a
meromorphic function $g(z)$ with the principal part
$g_j(z)=n_j(z-a_j)^{-1}$, $n_j=\partial (a_j)$,
$g(z):=\sum_{j=1}^{\infty } (g_j(z)-h_j(z))$, where
$h_j(z):=-n_j \sum_{p=1}^{k_j}(a_j^{-1}z)^{p-1}a_j^{-1} $.
Choose $k_j$ such that in each bounded canonical closed $V$
in ${\cal A}_r$, $V\subset U$, the series for $g$ is unifomly convergent.
At $a_j$ of the second type construct \\
$g_j(z):=n_j \sum_{p=0}^{\infty }(z-b_j)^{-1}[(a_j-b_j)(z-b_j)^{-1}]^p$
and \\
$h_j(z):=n_j \sum_{p=0}^{k_j}(z-b_j)^{-1}[(a_j-b_j)(z-b_j)^{-1}]^p$, \\
where $|a_j-b_j| < |z-b_j|$, $a_j, b_j, z$ are subjected to the condition
$\Upsilon _{a_j-b_j,z-b_j}\hookrightarrow \bf K$, which is not
finally restrictive due to Theorem 2.15.
Choose $k_j$ such that for each $z$ with $|z-b_j|\ge R>\delta _j$
we have $|g_j(z)-h_j(z)| = |n_j \sum_{p=k_j+1}^{\infty }
(z-b_j)^{-1}[(a_j-b_j)(z-b_j)^{-1}]^p |\le n_j \sum_{p=k_j+1}^{\infty }
|(a_j-b_j)^pR_j^{-p-1} | <\epsilon _j$, where $\sum_{j=1}^{\infty }
\epsilon _j<\infty $ converges, $R_j>0$ is a constant for each $j$.
These series by $(a_j^{-1}z)^pa_j^{-1} $ or by
$(z-b_j)^{-1}[(a_j-b_j)(z-b_j)^{-1}]^p$, $p=0,1,2,...$,
are with real coefficients independent of the type of embedding of
$\bf K$ into ${\cal A}_r$, where $z\in \Upsilon _{a_j,b_j,z_1}\hookrightarrow
{\bf K}\hookrightarrow {\cal A}_r$ for each given $z_1\in {\cal A}_r$
with the variable $z$ within a given copy of $\bf K$,
hence they can be extended on the corresponding balls
in ${\cal A}_r$. Now integrate $g(\zeta )$ along a rectifiable path
$\gamma $ in $U$ which does not contain any $a_j$,
$\gamma (0)=z_0$, $\gamma (1)=z$. Then $\int_{z_0}^z
g(\zeta )d\zeta =\sum_{j=1}^{\infty }[w_j(z)-w_j(z_0)]$
such that $f_1(z):=\exp (\int_{z_0}^z g(\zeta )d\zeta )$ is independent
of the path (see 
Theorems $2.15$ and $3.8.3$ and Corollary $3.4$ above), where \\
$\exp (w_j(z))=  ((1-a_j^{-1}z) \exp (\sum_{p=1}^{\infty }
(a_j^{-1}z)^p/p]))^{n_j}$ in the first case,  \\
$\exp (w_j(z))=([(z-a_j)(z-b_j)^{-1}] \exp (\sum_{p=1}^{\infty }
[(a_j-b_j)(z-b_j)^{-1}]^p/p))^{n_j}$ in the second case,
but $w_j$ and $w_l$ generally do not commute for $j\ne l$.
The convergence of series is analogous to that of
Satz 25 \cite{behnsomm} in the complex case.
Two functions satisfying Theorem $3.32$ need not differ on
an ${\cal A}_r$-holomorphic multiplier apart from the complex case, since,
for example, $\zeta _j:=\{ f_1...f_jgf_{j+1}...f_n \} _{q(n)}$
with different $j=a$ and $j=b$ in $\{ 1,...,n \} \subset \bf N$
do not correlate: $\zeta _a\ne \{ h\zeta _b k \} _{q(3)}$
in general for any ${\cal A}_r$-holomorphic $h$ and $k$ functions
on $U$, where each $f_l$ has a zero of order $n_l>0$ at $a_l$,
$n\ge 2$. Moreover, for $r\ge 3$ there is also dependence on
the order of multiplication $\{ * \} _{q(n)}$.
\par {\bf 3.33. Theorem.} {\it Let $U$ be an open region
in ${\cal A}_r$, $2\le r\le \infty $, and $f$ be a function
${\cal A}_r$-holomorphic on $U$
with a right superlinear superdifferential on $U$.
Suppose $f$ is not constant and $B(a,\rho ,{\cal A}_r)\subset U$,
where $0<\rho <\infty $. Then $f(B(a,\rho ,{\cal A}_r))$ is a neighbourhood
of $f(a)$ in ${\cal A}_r$.}
\par {\bf 3.34. Remarks.} For several ${\cal A}_r$
variables a multiple ${\cal A}_r$ line integral
${\bf I}:=
\{ \int_{\gamma _n}...\int_{\gamma _1} f(\mbox{ }^1z,...,\mbox{ }^nz
d\mbox{ }^1z...d\mbox{ }^nz \} _{q(n)}$ may be naturally considered
for rectifiable curves $\gamma _1$,...,$\gamma _n$ in ${\cal A}_r$,
$3\le r\le \infty $, where $\{ * \} _{q(n)}$ denotes the order
of brackets in the order or iterated integrations (see also
\S 2.1). Generally, the order of integration is essential, since
the existence of the partial derivative
$\partial ^n g(\mbox{ }^1z,...,\mbox{ }^nz)/
\partial \mbox{ }^1z ... \partial \mbox{ }^nz $ does not imply
an existence of a continuous $g^{(n)}$ and as Proposition
3.8.5 shows the order of differentiation is essential,
for example, even in the case of $g$ corresponding to
$f := \{ f_1\circ ... \circ f_n \} _{q(n)}$
with $f_n(\mbox{ }^1z,...,\mbox{ }^nz)$ with values in ${\cal A}_r^{n-1}$,
$f_{n-1}(\mbox{ }^1z,...,\mbox{ }^{n-1}z)$ with values in ${\cal A}_r^{n-2}$,
...,$f_2(\mbox{ }^1z,\mbox{ }^2z)$ and $f_1(\mbox{ }^1z)$
with values in ${\cal A}_r$, $ \{ \partial ^n g(\mbox{ }^1z,...,\mbox{ }^nz)/
\partial \mbox{ }^1z ... \partial \mbox{ }^nz \} _{q(n)} .
1^{\otimes n} =f(\mbox{ }^1z,...,\mbox{ }^nz)$ (see also \S 2.7).
Therefore, there is the natural generalization
of Theorem 3.9 for several ${\cal A}_r$ variables:
$$(i)\quad (2\pi )^nf(z_0) \{ M_1...M_n \} _{q(n)} =$$
$$ \{ \int_{\psi _n} ...
\int_{\psi _1}f(\mbox{ }^1\zeta ,...,\mbox{ }^n\zeta )(\mbox{ }^1\zeta -
\mbox{ }^1z_0)^{-1}d\mbox{ }^1\zeta )...
(\mbox{ }^n\zeta - \mbox{ }^nz_0)^{-1}d\mbox{ }^n\zeta ) \} _{q(n)} $$
for the corresponding $U=\mbox{ }^1U\times ...\times \mbox{ }^nU,$
where $\psi _j$ and $\mbox{ }^jU$ satisfy conditions of Theorem 3.9
for each $j$ and $f$ is a continuous ${\cal A}_r$-holomorphic
function on $U$.

\thanks{
Address: Sergey V. Ludkovsky, Mathematical Department, \\
Brussels University, V.U.B.,
Pleinlaan 2, Brussels 1050, Belgium.\\
{\underline {Acknowledgment}}. The author thanks the Flemish
Science Foundation for support through the Noncommutative Geometry
from Algebra to Physics project and Prof. Fred van Oystaeyen
for his interest and fruitful discussions.}
\end{document}